\documentclass{article}
\usepackage{arxiv}
\usepackage{comment}
\usepackage[utf8]{inputenc} 
\usepackage[T1]{fontenc}    
\usepackage{hyperref}       
\usepackage{url}            
\usepackage{booktabs}       
\usepackage{amsmath}
\usepackage{amssymb}
\usepackage{amsfonts}       
\usepackage{nicefrac}       
\usepackage{microtype}      
\usepackage{xcolor}
\usepackage{graphicx}
\usepackage{stmaryrd}
\SetSymbolFont{stmry}{bold}{U}{stmry}{m}{n}
\usepackage{bm}
\usepackage[normalem]{ulem}
\usepackage[tickmarkheight=0.1cm,textsize=tiny]{todonotes}
\usepackage{algorithmicx}
\usepackage{algorithm}
\usepackage{algpseudocode}
\usepackage{caption}
\usepackage{subcaption}
\usepackage{lineno}
\usepackage{orcidlink}
\usepackage[inline,final]{showlabels}

\usepackage{xpatch}
\makeatletter
\xpatchcmd{\algorithmic}
  {\ALG@tlm\z@}{\leftmargin\z@\ALG@tlm\z@}
  {}{}
\makeatother

\hypersetup{
    colorlinks=true, 
    linkcolor=blue, 
    urlcolor=red, 
    citecolor=magenta,
    linktoc=all 
}

\newlength\myindent 
\newcommand{\re}{\mathbb{R}} 

\newcommand{\poly}{\mathbb{P}} 

\newcommand{\emh}{{e-\frac{1}{2}}}

\newcommand{\eph}{{e+\frac{1}{2}}}

\newcommand{\pph}{p + \frac{1}{2}}
\newcommand{\pmh}{p - \frac{1}{2}}

\newcommand{\Nmh}{{N-\frac{1}{2}}}
\newcommand{\Nph}{{N+\frac{1}{2}}}
\newcommand{\half}{\frac{1}{2}}
\newcommand{\ud}{\textrm{d}}
\newcommand{\pd}[2]{\frac{\partial #1}{\partial #2}}

\newcommand{\avg}[1]{\overline{#1}}
\newcommand{\au}{\avg{\uu}}
\newcommand{\tu}{\tilde{\uu}}
\newcommand{\bL}{\bs{L}}

\newcommand{\myvector}[1]{\mathsf{#1}}
\newcommand{\vu}{\myvector{u}}

\newcommand{\vD}{\myvector{D}}

\newcommand{\vR}{\myvector{R}}

\newcommand{\bs}{\boldsymbol}

\newcommand{\mathLaplace}{\Delta} 
\newcommand{\extrapolate}{\textbf{AE}} 
\newcommand{\evaluate}{\textbf{EA}}

\newcommand{\Uad}{\mathcal{U}_{\textrm{ad}}}

\newtheorem{definition}{Definition}

\usepackage{graphicx,latexsym,tabularx,xcolor}

\newcommand{\tmmathbf}[1]{\ensuremath{\boldsymbol{#1}}}
\newcommand{\tmop}[1]{\text{#1}}

\newtheorem{theorem}{Theorem}

\newcommand{\bzero}{\boldsymbol{0}}
\newcommand{\uu}{\boldsymbol{u}}
\newcommand{\uU}{\boldsymbol{U}}
\newcommand{\bv}{\boldsymbol{v}}
\newcommand{\ff}{\boldsymbol{f}}
\newcommand{\F}{\boldsymbol{F}}
\newcommand{\pf}{\ff}

\newcommand{\pg}{\boldsymbol{g}}

\newcommand{\fg}{\boldsymbol{g}}
\newcommand{\bD}{\boldsymbol{D}}
\newcommand{\fG}{\boldsymbol{G}}
\newcommand{\bss}{\boldsymbol{s}}
\newcommand{\bS}{\boldsymbol{S}}

\newcommand{\pdx}{\partial_x}
\newcommand{\ddgx}{\partial_x^\text{DG}}
\newcommand{\dfrx}{\partial_x^\text{FR}}
\newcommand{\dlocx}{\partial_x^\text{loc}}
\newcommand{\nc}{M}
\newcommand{\nel}{\ensuremath{N_e}}
\newcommand{\ad}{P}

\newcommand{\uep}{\uu_{e,p}}

\newcommand{\uez}{\tmmathbf{u}_{e, 0}}

\newcommand{\uepoz}{\tmmathbf{u}_{e + 1, 0}}
\newcommand{\ueN}{\tmmathbf{u}_{e, N}}
\newcommand{\xep}{x^e_p}

\newcommand{\wg}{w}

\newcommand{\discf}{\ensuremath{\pf_h^\delta}}
\newcommand{\discfref}{\ensuremath{\hat{\pf}_{h,e}^\delta}}
\newcommand{\discfrefe}{\ensuremath{\hat{\pf}_h^\delta}}
\newcommand{\discF}{\ensuremath{\F_h^\delta}}
\newcommand{\discFref}{\ensuremath{\hat{\F}_{h,e}^\delta}}
\newcommand{\discFrefe}{\ensuremath{\hat{\F}_h^\delta}}

\newcommand{\cE}{E}

\newcommand{\Eonetwo}{E_{12}}
\newcommand{\Etwotwo}{E_{22}}
\newcommand{\vone}{v_1}
\newcommand{\vtwo}{v_2}
\newcommand{\ccE}{\tmmathbf{E}}

\newcommand{\cR}{p}
\newcommand{\ccP}{\tmmathbf{p}}
\newcommand{\ccR}{\ensuremath{\tmmathbf{p}}}
\newcommand{\Poneone}{p_{11}}
\newcommand{\Ponetwo}{p_{12}}
\newcommand{\Ptwotwo}{p_{22}}

\newcommand{\correction}[1]{{#1}}

\title{Compact Runge-Kutta Flux Reconstruction for Hyperbolic Conservation Laws with admissibility preservation}

\author{
Arpit~Babbar \orcidlink{0000-0002-9453-370X} \\
Institute of Mathematics\\
Johannes Gutenberg University Mainz\\
Staudingerweg 9, 55122 Mainz, Germany\\
\texttt{ababbar@uni-mainz.de} \\
\And
Qifan Chen \orcidlink{0009-0003-0790-6972}\\
Department of Mathematics\\
The Ohio State University Columbus\\
OH 43210 USA\\
\texttt{chen.11010@osu.edu}
}
\begin{document}
\maketitle
\begin{abstract}
Compact Runge-Kutta (cRK) Discontinuous Galerkin (DG)  methods, recently introduced in [Q. Chen, Z. Sun, and Y. Xing,
SIAM J. Sci. Comput., 46: A1327–A1351, 2024], are a variant of RKDG methods for solving hyperbolic conservation laws and are characterized by their compact stencil including only immediate neighboring finite elements. This article proposes a cRK Flux Reconstruction (FR) method by interpreting cRK as a procedure to approximate time-averaged fluxes, which requires computing only a single numerical flux for each time step and further reduces data communication. The numerical flux is carefully constructed to maintain the same Courant–Friedrichs–Lewy (CFL) numbers as cRKDG methods and achieve optimal accuracy uniformly across all polynomial degrees, even for problems with sonic points. A subcell-based blending limiter is then applied for problems with nonsmooth solutions, which uses Gauss-Legendre solution points and performs MUSCL-Hancock reconstruction on subcells to mitigate the additional dissipation errors. Additionally, to achieve a fully admissibility preserving cRKFR scheme, a flux limiter is applied to the time-averaged numerical flux to ensure admissibility preservation in the means, combined with a positivity preserving scaling limiter. The method is further extended to handle source terms by incorporating their contributions as additional time averages. Numerical experiments including Euler equations and the ten-moment problem are provided to validate the claims regarding the method's accuracy, robustness, and admissibility preservation.
\end{abstract}
\section{Introduction}

The compact Runge-Kutta Discontinuous Galerkin (cRKDG) method was introduced in~\cite{chen2024} as a high order method to solve hyperbolic conservation laws using a more compact stencil compared to the original RKDG method proposed in~\cite{Cockburn1991,Cockburn1989a,Cockburn1989}. The goal of this article is to develop the compact Runge-Kutta method in a Flux Reconstruction (FR) framework~\cite{Huynh2007} that requires only one numerical flux computation per time step and use this framework to obtain a provably admissibility preserving scheme with the flux limiter proposed in~\cite{babbar2024admissibility}.

Hyperbolic conservation laws are a first order system of partial differential equations (PDE) that represent the conservation of mass, momentum, and energy. They are used for various practical applications like Computational Fluid Dynamics (CFD), astrophysics and weather modeling. These applications, due to their complexity and scale, often necessitate substantial computational resources. Thus, developing efficient numerical methods for solving these PDE is an important area of research. Higher-order methods are particularly relevant in the current state of memory-bound High Performance Computing (HPC) hardware~\cite{attig2011,subcommittee2014} because their high arithmetic intensity makes them less likely to be constrained by memory bandwidth. However, lower order methods are more robust and thus still routinely used in practice. This work is in direction of taking the best of both methods by developing high order methods that adaptively switch to low order methods in regions where robustness in required.

Discontinuous Galerkin (DG) is a Spectral Element Method with a discontinuous finite element basis. It was first introduced by Reed and Hill~\cite{reed1973} for neutron transport equations and developed for fluid dynamics equations by Cockburn and Shu and others, see~\cite{Cockburn1991,Cockburn1989a,Cockburn1989} and other references in~\cite{cockburn2000}. The DG method preserves local conservation, allows simple $hp$-refinement~\cite{Kopriva2002} and can be fitted into complex geometries~\cite{Kopriva2006,kopriva2009}. It is also suitable for modern HPC as the neighbouring DG elements are coupled only through the numerical flux and thus bulk of computations are local to the element, minimizing data communication. Flux Reconstruction (FR) method introduced by Huynh~\cite{Huynh2007} is also a class of discontinuous Spectral Element Methods, having the key property of being quadrature-free. The central idea in this method is to use the numerical flux to construct a continuous flux approximation and then collocate the conservation law at solution points. The collocation nature of the scheme leads to an efficient implementation that can exploit optimized matrix-vector operations and vectorization capabilities of modern CPUs. A crucial ingredient of the FR scheme is the use of the correction functions to construct the continuous flux approximation. The choice of correction functions affects the accuracy and stability of the method~\cite{Huynh2007,Vincent2011a,Vincent2015,Vermeire2016,Trojak2021}. In~\cite{Cicchino2022}, a Flux Reconstruction scheme with nonlinear stability~\cite{Jameson2012} was constructed, one of whose key ingredients was application of correction functions to the volume terms. By properly choosing the correction function and solution points, FR method can be shown to be equivalent to some discontinuous Galerkin and spectral difference schemes~\cite{Huynh2007,Trojak2021}.

In the pioneering works of DG/FR methods~\cite{Cockburn1991,Cockburn1989a,Cockburn1989,Huynh2007}, these spectral element methods are used as the spatial discretization to obtain a semidiscrete system. Following the commonly used method of lines (MOL), the semidiscrete system is then solved using a multistage Runge-Kutta (RK) method. This requires applying the FR/DG method at every RK stage, necessitating interelement communication at each stage. The alternative to the multistage approach are single stage solvers which specify a complete space time discretization. We briefly mention some of the well known single stage methods, an in-depth review can be found in~\cite{babbar2022,babbar2024thesis}. ADER (Arbitrary high order DERivative) schemes are one class of single-stage evolution methods which originated as ADER Finite Volume (FV) schemes~\cite{Titarev2002,Titarev2005} and are also used as ADER-DG schemes~\cite{Dumbser2008,Dumbser2014}. Another class of single-stage evolution methods are the Lax-Wendroff schemes which were originally proposed in the finite difference framework in~\cite{Qiu2003} with the WENO approximation of spatial derivatives~\cite{Shu1989} and later extended to DG framework in~\cite{Qiu2005b,sun2017stability}. The cRKDG method of~\cite{chen2024} is a recent method that performs single stage evolution by doing the interelement communication only in the last stage of the Runge-Kutta method. For linear problems, the authors in~\cite{chen2024} showed that the cRKDG method is equivalent to Lax-Wendroff methods and to the ADER-DG method for a particular choice of RK method.

In this work, we propose a compact Runge-Kutta Flux Reconstruction (cRKFR) method with a focus on obtaining an admissibility preserving scheme. A crucial part of our method is regarding compact Runge-Kutta (cRK) as a procedure to approximate time average fluxes, similar to the Lax-Wendroff Flux Reconstruction methods~\cite{Qiu2005b,babbar2022,babbar2024admissibility}. In the inner Runge-Kutta stages, we perform local evolutions by computing derivatives of \textit{discontinuous flux} approximations within each element. The local evolutions are used to compute a \textit{local time averaged flux approximation}. In the final stage, the local flux approximation is made into the \textit{continuous time averaged flux approximation} using FR correction functions~\cite{Huynh2007}, and is used to perform the evolution to the next time level. Although the stencil of each element per time step of~\cite{chen2024} also consists only of its neighbouring element, the final stage in~\cite{chen2024} performed computation of a numerical flux corresponding to each inner stage requiring computation of $s$ numerical fluxes for a Runge-Kutta method with $s$ inner stages. In our approach with the time averaged flux view point, we only need to compute one numerical flux per time step, which saves computations and storage and can also reduce the volume of interelement communication.

Many practical problems involving hyperbolic conservation laws involve nonsmooth solutions arising from rarefactions, shocks and other discontinuities. The usage of high order methods for such problems requires the use of limiters to prevent spurious oscillations. A review of various limiters for DG/FR methods can be found in~\cite{babbar2024admissibility}. In the first work on cRKDG schemes~\cite{chen2024}, the TVB limiter of~\cite{Cockburn1989a} was used. The TVB limiter does not preserve any subcell information other than the element mean and trace values. In this work, we use the subcell based blending limiter of~\cite{babbar2024admissibility} inspired from~\cite{hennemann2021} to control spurious oscillations while using information from subcells. As in~\cite{babbar2024admissibility}, the subcell limiter is designed to minimize dissipation errors by using Gauss-Legendre solution points and performing MUSCL-Hancock reconstruction on subcells. These give us a scheme that produces better resolution than TVB limiters and the original blending scheme proposed in~\cite{hennemann2021}, which we demonstrate through numerical experiments.

The subcell based limiting is heuristically designed to balance accuracy and robustness, and it allows some spurious oscillations to remain. These oscillations may cause numerical solutions to be out of physically admissible bounds, e.g., leading to negative density and pressure values. Nonphysical solutions typically cause numerical simulations to crash, and it is thus essential to perform additional limiting to ensure that the physical bounds are maintained. A common approach, proposed in~\cite{Zhang2010b}, is to use a Strong Stability Preserving RK method~\cite{Gottlieb2001} that gives solutions whose element means are admissible in every stage, and combine them with the scaling limiter of~\cite{Zhang2010b} to obtain physically admissible solutions. In the recent work of~\cite{liu2024}, admissibility in means for each stage of several cRKDG schemes was proven. Thus, by applying the scaling limiter of~\cite{Zhang2010b} at each of the stages, admissibility preserving cRKDG schemes were obtained in~\cite{liu2024}. In~\cite{babbar2024admissibility}, a limiting procedure for the time averaged flux was proposed which gave admissibility of the solution. Interpreting the cRK procedure with   time averaging  fluxes, we develop a similar limiting procedure to obtain an admissibility preserving cRKFR scheme. Our limiting procedure works for any cRK scheme with the same limiting process irrespective of the order of the method, which includes one application of the flux limiter from~\cite{babbar2024admissibility} and one application of the scaling limiter from~\cite{Zhang2010b} for each time step. This is different from the admissibility preserving cRKDG schemes obtained in~\cite{liu2024}, where $s$ scaling limiters of~\cite{Zhang2010b} are required for an $s$ stage Runge-Kutta method. Additionally, the limiting procedure of~\cite{liu2024} requires computation of multiple numerical fluxes per time step as in~\cite{chen2024} while our proposed procedure is compatible with the single time averaged numerical flux computation proposed in this work. The limiting procedure also allows us to exploit the potential advantages of single stage methods that only require one application of limiters per time step.

The rest of the paper is organized as follows. In Section~\ref{sec:rkfr}, we review the standard Runge-Kutta Flux Reconstruction method by deriving it as the Discontinuous Galerkin method whose quadrature points are taken to be the nodal solution points. Section~\ref{sec:rkfr} also introduces the notations used in this work. In Section~\ref{sec:crkfr}, the compact Runge-Kutta Flux Reconstruction (cRKFR) scheme is introduced with its capability to perform evolution with only one \textit{time averaged numerical flux} computation. The choices for the computation of the central ({\extrapolate} and {\evaluate}) and dissipative part (D-CSX, D1, D2) of the time averaged numerical flux are discussed in Section~\ref{sec:numflux}.  In Section~\ref{sec:blending.scheme}, we describe the blending limiter for the cRKFR method. The flux limiter that ensures admissibility of means is described as a part of the blending scheme in Section~\ref{sec:flux.limiter}. In Section~\ref{sec:numerical.results}, we present numerical results to verify accuracy and robustness of our scheme using some scalar equations, compressible Euler equations and the ten moment problem. A summary of the proposed scheme is given in Section~\ref{sec:conclusion}.

In Appendix~\ref{app:time.averaged.flux}, we give a formal order of accuracy justification of viewing the cRK scheme as time averages. The viewpoint is then used to prove the linear equivalence of the proposed cRKFR scheme with the Lax-Wendroff schemes of~\cite{Qiu2005b,babbar2022} depending on the dissipation model (Section~\ref{sec:numflux}). In Appendix~\ref{app:source.term}, we describe the extension of the cRKFR scheme to handle source terms.

\section{Runge-Kutta Flux Reconstruction} \label{sec:rkfr}
\begin{figure}
\begin{center}
\begin{tabular}{cc}
\includegraphics[width=0.45\textwidth]{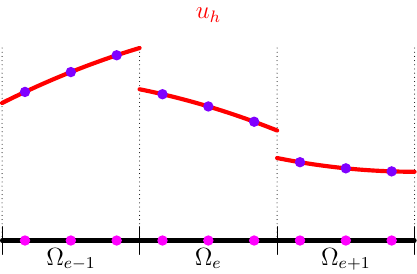} &
\includegraphics[width=0.45\textwidth]{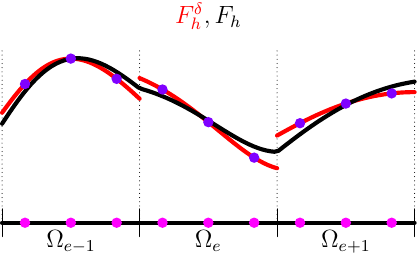} \\
(a) & (b)
\end{tabular}
\end{center}
\caption{(a) Piecewise polynomial solution at time $t_n$, and (b) discontinuous and continuous flux. The figure has been taken from~\cite{babbar2022}.}
\label{fig:solflux1}
\end{figure}

Consider the 1-D hyperbolic conservation laws
\begin{equation} \label{eq:con.law}
\uu_t + \pf(\uu)_x = \bzero, \qquad \uu(x,0) = \uu_0(x),\qquad x \in \Omega,
\end{equation}
where the solution $\uu \in \re^\nc$ is the vector of conserved quantities, $\Omega$ is the physical domain, $\pf(\uu)$ is the physical flux, $\uu_0$ is the initial condition, and some boundary conditions are additionally prescribed. In many of these applications, it is essential that the solution is physically correct, i.e. it belongs to an admissible set, denoted by $\Uad$. For example, in gas dynamics governed by the compressible Euler's equations~(\ref{eq:1deuler},~\ref{eq:2deuler}), the density and pressure are positive. In case of relativistic hydrodynamics equations, the velocity has to be less than unity. For the models that we study, we assume that the admissible set is a convex subset of $\re^\nc$, and can be written as
\begin{equation}
\label{eq:uad.form} \Uad = \{ \uu \in \re^\nc : \ad_k (\uu) > 0, 1 \le k \le K\},
\end{equation}
where each admissibility constraint $\ad_k$ is a concave function if $\ad_j > 0$ for all $j <
k$. For Euler's equations, $K = 2$ and $\ad_1, \ad_2$ are density, pressure
functions respectively; if the density is positive then pressure is a concave
function of the conserved variables. In this work, we cast the discretization of~\cite{chen2024} as the quadrature free Flux Reconstruction (FR) scheme of Huynh~\cite{Huynh2007}. This framework will be used to construct a cRKFR scheme requiring computation of only one numerical flux (Section~\ref{sec:crkfr}) and to apply the admissibility enforcing flux limiters of~\cite{babbar2024admissibility} (Section~\ref{sec:flux.limiter}). We now describe the notations for the finite element grid. Divide the physical domain $\Omega$ into disjoint elements $\{\Omega_e\}$ with
\begin{equation} \label{eq:reference.element}
\Omega_e = [x_{\emh}, x_{\eph}] \qquad \textrm{and} \qquad \Delta x_e =
x_{\eph} - x_{\emh}.
\end{equation}
Let us map each element to a reference element, $\Omega_e \to [0, 1]$, by
\[
x \mapsto \xi = \frac{x - x_{\emh}}{\Delta x_e}.
\]
Inside each element, we approximate the numerical solution to~\eqref{eq:con.law} by $\poly_N$ functions which are degree $N \geq 0$ polynomials so that the basis for numerical solution is
\begin{equation} \label{eq:fr.basis}
V_h = \{v_h : v_h|_{\Omega_e} \in \poly_N \}.
\end{equation}
Since the equation~\eqref{eq:con.law} admits solutions with discontinuities, the basis functions in~\eqref{eq:fr.basis} are allowed to be discontinuous at element interfaces, see Figure~\ref{fig:solflux1}a.
For each basis function $v_h \in V_h$~\eqref{eq:fr.basis} and a physical element $\Omega_e$~\eqref{eq:reference.element}, we define $\hat{v}_{h,e} = \hat{v}_{h,e}(\xi)$ through reference map~\eqref{eq:reference.element} as
\begin{equation} \label{eq:vh.reference}
\hat{v}_{h,e}(\xi) = v_h(x_{\emh} + \xi \Delta x_e) = v_h(x).
\end{equation}
The element index $e$ will often be suppressed for brevity. To construct $\poly_N$ polynomials in the basis, choose $N + 1$ distinct nodes
\begin{equation} \label{eq:soln.points}
0 \le \xi_0 < \xi_1 < \cdots < \xi_N \le 1,
\end{equation}
which will be taken to be Gauss-Legendre (GL) or Gauss-Legendre-Lobatto (GLL)
nodes, and will also be referred to as \textit{solution points}. There are
associated quadrature weights $w_j$ such that the quadrature rule is exact for
polynomials of degree up to $2 N + 1$ for GL points and up to degree $2 N - 1$
for GLL points. Note that the nodes and weights we use are with respect to the
interval $[0, 1]$ whereas they are usually defined for the interval $[- 1, +
1]$. The numerical solution $\uu_h \in V_h^\nc$, where $\nc$ is the number of conservative variables~\eqref{eq:con.law}, inside an element $\Omega_e$ is given in reference coordinates as
\begin{equation} \label{eq:soln.poly}
x \in \Omega_e : \qquad \hat{\uu}_{h,e} (\xi, t) = \sum_{p = 0}^N \uu_{e,p} (t) \ell_p (\xi),
\end{equation}
where each $\ell_p$ is a Lagrange polynomial of degree $N$ given by
\begin{equation} \label{eq:defn.lagrange}
\ell_q (\xi) = \prod_{p = 0, p \ne q}^N \frac{\xi - \xi_p}{\xi_q - \xi_p}
\in \poly_N, \qquad \ell_q (\xi_p) = \delta_{pq}.
\end{equation}
The numerical methods to solve~\eqref{eq:con.law} require computation of spatial derivatives. The spatial derivatives are computed on the reference interval $[0,1]$ using the differentiation matrix $\vD = [D_{ij}]$ whose entries are given by
\begin{equation} \label{eq:diff.matrix}
D_{ij} = \ell_j'(\xi_i),\qquad 0 \le i, j \le N.
\end{equation}

Figure~\ref{fig:solflux1}a illustrates a piecewise polynomial solution at
some time $t_n$ with discontinuities at the element boundaries. In order to clearly show our connection with the work of~\cite{chen2024}, we review how the standard Runge-Kutta Discontinuous Galerkin (RKDG) method in the notations of~\cite{chen2024} is equivalent to the Runge-Kutta Flux Reconstruction (RKFR) method of~\cite{Huynh2007}. Following~\cite{chen2024}, we define the discrete DG operator $\ddgx$ so that $\ddgx \pf (\uu_h) \in V_h^\nc$ is prescribed by its action on test functions $v_h \in V_h$ as
\begin{equation} \label{eq:defn.ddg}
\int_{\Omega_e} \ddgx  \pf(\uu_h) v_h \ud x = -\int_{\Omega_e} \pf(\uu_h) \pdx v_h
+ \pf_{\eph} \hat{v}_h(1) - \pf_{\emh} \hat{v}_h(0),
\end{equation}
where $\pf_{\eph}$ denotes the numerical flux at the interface $x_{\eph}$. Then $\pf_{\eph}$ is computed by
\begin{equation} \label{eq:num.flux.fr}
\pf_\eph(\uu_h) = \frac{1}{2}(\pf(\uu_\eph^-) + \pf(\uu_\eph^+)) - \frac {\lambda_{\eph}}{2} (\uu_\eph^+ - \uu_\eph^-), \qquad \uu_\eph^\pm = \uu_h(x_\eph^\pm),
\end{equation}
where $\lambda_{\eph} = \lambda(\uu_{\eph}^-,\uu_{\eph}^+)$ is an approximation of the wave speed at the interface, often computed using a Rusanov~\cite{Rusanov1962} or local Lax-Friedrichs approximation as
\begin{equation} \label{eq:rusanov.wave.speed}
\lambda(\uu_l, \uu_r) = \max_{\uu \in \{\uu_l, \uu_r\}} \sigma(\pf'(\uu)).
\end{equation}
Here $\sigma(A)$ denotes the spectral radius of a matrix $A$. Equation~\eqref{eq:defn.ddg} to define $\ddgx$ is obtained by a formal integration by parts in space. In order to cast the above scheme in a Flux Reconstruction framework, we define a degree $N$ \textit{discontinuous flux} approximation $\discf$ (Figure~\ref{fig:solflux1}b) in reference coordinates for element $e$ as
\begin{equation} \label{eq:discts.flux}
\discfref(\xi) = \sum_{p=0}^N \pf(\uu_{e,p}) \ell_p(\xi).
\end{equation}
Then, assuming that the quadrature points to approximate the integral in~\eqref{eq:defn.ddg} are same as the solution points~\eqref{eq:soln.points} used to define the solution polynomial~\eqref{eq:soln.poly}, the action of the DG operator $\ddgx$ can be written in terms of the discontinuous flux polynomial~\eqref{eq:discts.flux} as
\begin{equation*}
\int_{\Omega_e} \ddgx  \pf(\uu_h) v_h \ud x = -\int_{\Omega_e} \discf  \pdx v_h \ud x + \pf_{\eph} \hat{v}_h(1) - \pf_{\emh} \hat{v}_h(0),\qquad v_h \in V_h.
\end{equation*}
When the solution points are chosen to be the commonly used Gauss-Legendre (GL) or Gauss-Legendre-Lobatto (GLL) points, the flux integral is exact and we can perform another integration by parts to obtain the \textit{strong form} of the operator $\ddgx$ as
\[
\int_{\Omega_e} \ddgx \pf(\uu_h) v_h \ud x = \int_{\Omega_e} \pdx \discf v_h \ud x + (\pf_\eph - \discfrefe(1)) \hat{v}_h(1) + (\pf_\emh - \discfrefe(0)) \hat{v}_h(0),\qquad v_h \in V_h.
\]
Substituting $v_h$ with the Lagrange polynomial $\ell_p$~\eqref{eq:defn.lagrange} for each element $\Omega_e$, this is further equivalent to
\begin{align}
\int_{\Omega_e} \ddgx  \pf(\uu_h) \ell_p \ud x &= \int_{\Omega_e} \pdx \discf  \ell_p \ud x + (\pf_\eph - \discfrefe(1)) \ell_p(1) + (\pf_\emh - \discfrefe(0)) \ell_p(0),\quad 0 \le p \le N \nonumber.
\end{align}
Since we have assumed that the quadrature points are same as the solution points, performing quadrature and dividing by $\Delta x_e w_p$ gives
\begin{align}
\ddgx  \pf(\uu_h) (\xi_p) &= \pdx \discf(\xi_p) + (\pf_\eph - \discfrefe(1)) \frac{\ell_p(1)}{\Delta x_e \wg_p} + (\pf_\emh - \discfrefe(0)) \frac{\ell_p(0)}{\Delta x_e \wg_p}, \quad 0 \le p \le N \nonumber \nonumber \\
& = \pdx \discf(\xi_p) + (\pf_\eph - \discfrefe(1)) \frac{g_R'(\xi_p)}{\Delta x_e} + (\pf_\emh - \discfrefe(0)) \frac{g_L'(\xi_p)}{\Delta x_e}. \label{eq:fr.equiv}
\end{align}
For simplicity, this can be rewritten as
\begin{equation}
  \ddgx \pf(\uu_h)  = \pdx \pf_h=: \dfrx \pf(\uu_h), \quad \pf_h  = \discf + (\pf_\eph - \discfrefe(1)) g_R + (\pf_\emh - \discfrefe(0)) g_L, \label{eq:dfrx.defn}
\end{equation}
where $g_L, g_R \in \poly_{N+1}$ are Flux Reconstruction (FR) correction functions satisfying
\[
g_L(0) = g_R(1) = 1, \qquad g_L(1) = g_R(0) = 0.
\]
The $\pf_h$ defined in~\eqref{eq:dfrx.defn} is called the \textit{continuous flux approximation} in the FR literature~\cite{Huynh2007} as it is globally continuous, taking the numerical flux value $\pf_{e+1/2}$ at the each element interface $e+1/2$. We choose correction functions $g_L, g_R \in \mathbb{P}_{N + 1}$ to be
$g_{\text{Radau}}, g_2$~{\cite{Huynh2007}} if the solution points are taken to be GL, GLL respectively. For both of the cases, equation~\eqref{eq:fr.equiv} is obtained by the following identities from the proof
of equivalence of FR and DG in~\cite{Huynh2007}\footnote{For
Radau correction functions, the identities follow from taking $\phi = \ell_p$
in (7.8) of~{\cite{Huynh2007}}. For $g_2$, see Section 2.4
of~{\cite{Grazia2014}}. A detailed proof is also available in Appendix B of~\cite{babbar2024thesis}.}
\[
g_R'(\xi_p) = \ell_p (1) / w_p, \qquad g_L' (\xi_p) = - \ell_p (0) / w_p.
\]
The $\ddgx$ operator is used to obtain a semidiscretization of~\eqref{eq:con.law} in~\cite{chen2024}. By~\eqref{eq:dfrx.defn} the semidiscretization of~\cite{chen2024}, when quadrature points are taken to be the same as the solution points, is equivalent to the following
\begin{equation} \label{eq:semidiscretization.dg}
\partial_t \uu_h + \dfrx \pf(\uu_h) = 0,
\end{equation}
which is the original FR method of~\cite{Huynh2007}. A fully discrete scheme to solve~\eqref{eq:con.law} is obtained by applying an explicit Runge-Kutta (RK) method to discretize~\eqref{eq:semidiscretization.dg} in time. An explicit RK method is specified by its Butcher tableau
\begin{equation}\label{eq:butcher}
\begin{array}{c|c}
c&A\\\hline
&b\\
\end{array},\quad A = (a_{ij})_{s\times s}, \quad b = (b_1,\dots, b_s), \quad c_i = \sum_j a_{ij},
\end{equation}
where $A$ is a lower triangular matrix in \eqref{eq:butcher}, namely, $a_{ij} = 0$ if $i \ge j$. The corresponding RKFR scheme is given by
\begin{subequations}\label{eq:rkfr}
\begin{align}
\uu_h^{(i)} &= \uu_h^n -  \Delta t\sum_{j = 1}^{i-1}  a_{ij} \dfrx \pf(\uu_h^{(j)}), \quad  i = 1, 2, \dots, s,\label{eq:rkfr1}\\
\uu_h^{n+1} &= \uu_h^n - \Delta t \sum_{i = 1}^s b_i \dfrx \pf(\uu_h^{(i)} ).\label{eq:rkfr2}
\end{align}
\end{subequations}

\section{Compact Runge-Kutta Flux Reconstruction} \label{sec:crkfr}
Following~\cite{chen2024}, we define a local discrete spatial operator $\dlocx$ where the inter-element terms from~\eqref{eq:dfrx.defn} are dropped to be
\begin{equation*}
\dlocx\pf(\uu_h)(\xi_p) = \pdx \discf(\xi_p).
\end{equation*}
The operator $\dlocx$ computes derivatives of the degree $N$ discontinuous flux approximation~\eqref{eq:discts.flux}. In practice, this operation is performed with a differentiation matrix~\eqref{eq:diff.matrix}. With this notation, a compact Runge-Kutta scheme in Flux Reconstruction framework with the Butcher tableau~\eqref{eq:butcher} is obtained as follows
\begin{subequations}\label{eq:crkfr}
\begin{align}
\uu_h^{(i)} &= \uu_h^n -  \Delta t\sum_{j = 1}^{i-1}  a_{ij} \dlocx \pf(\uu_h^{(j)}), \quad  i = 1, 2, \dots, s,\label{eq:crkfr1}\\
\uu_h^{n+1} &= \uu_h^n - \Delta t \sum_{i = 1}^s b_i \dfrx \pf(\uu_h^{(i)} ).\label{eq:crkfr2}
\end{align}
\end{subequations}
The difference between the above scheme~\eqref{eq:crkfr} and the standard Flux Reconstruction scheme lies in the fact that inner stages to evaluate $\{\uu_h^{(i)}\}_{i=1}^s$ in~\eqref{eq:crkfr1} use the local operator $\dlocx$ instead of $\dfrx$ in~\eqref{eq:rkfr1}. The cRK Discontinuous Galerkin (DG) scheme of~\cite{chen2024} uses the operator $\ddgx$ in place of~\eqref{eq:crkfr2}. Thus, the scheme~\eqref{eq:crkfr} is equivalent to the cRKDG scheme of~\cite{chen2024}, under the assumptions previously made assumptions~\eqref{eq:fr.equiv} that the quadrature points of the DG scheme are the same as the solution points. In order to perform evolution with a single numerical flux computation and to use the flux limiters of~\cite{babbar2024admissibility}, we will replace~\eqref{eq:crkfr} with an evolution in terms of approximations of the time average flux $\F = \int_{t^n}^{t^{n+1}} \pf(\uu) \ud t$, similar to the Lax-Wendroff schemes~\cite{Qiu2005b, babbar2022}. As discussed in Section~\ref{sec:numflux}, the proposed scheme can recover the cRKFR scheme~\eqref{eq:crkfr} by a choice of the dissipation term~\eqref{eq:numflux.general} in the numerical flux. A \textit{discontinuous time averaged flux approximation} similar to the approximation in~\eqref{eq:discts.flux} is obtained by using the inner stages in~\eqref{eq:crkfr} as
\begin{equation} \label{eq:disc.avg.flux}
\discFref = \sum_{p=0}^N \F_{e,p} \ell_p(\xi),\qquad \F_{e,p} = \sum_{i=1}^s b_i \pf({\uu}_{e,p}^{(i)}).
\end{equation}
A formal justification for why~\eqref{eq:disc.avg.flux} is an approximation to the time averaged flux can be found in Appendix~\ref{app:time.averaged.flux}. Thus, the proposed compact Runge-Kutta Flux Reconstruction (cRKFR) scheme is given by replacing the final stage of~\eqref{eq:crkfr} by the following application of the $\dfrx$ operator~\eqref{eq:dfrx.defn} on the time average flux
\begin{equation}\label{eq:crkfr.avg}
\uu_h^{n+1} = \uu_h^n - \Delta t \dfrx \F,\qquad \dfrx \F = \pdx \F_h  = \pdx \discF + (\F_\eph - \discFrefe(1)) \pdx g_R + (\F_\emh - \discFrefe(0)) \pdx g_L.
\end{equation}
Following~\cite{babbar2022}, $\F_\eph$ denotes the \textit{time averaged numerical flux}. The possible numerical flux approximations are discussed in Section~\ref{sec:numflux} and have an influence on stability and accuracy. As is further discussed in Section~\ref{sec:diss.models}, since the numerical flux computation is the only non-linear operation in terms of the flux in~\eqref{eq:crkfr}, it is also the choice of numerical flux that leads to equivalence of~\eqref{eq:crkfr.avg} with~\eqref{eq:crkfr} and thus with~\cite{chen2024}. Due to its generality over~\eqref{eq:crkfr}, cRKFR will always refer to~\eqref{eq:crkfr.avg}. We list the cRKFR schemes of orders 2,3, and 4 that are used in this work. To obtain a cRKFR scheme of order $N+1$, the RK scheme~\eqref{eq:butcher} of the same order is used and the solution polynomial is taken to be of degree $N$. Since the cRK scheme is cast in the FR framework as in~\eqref{eq:crkfr.avg}, we simply need to explain how the time averaged flux~\eqref{eq:disc.avg.flux} is computed in the cRK scheme using the inner stages. Since this step computes the local time average flux $\F$, we call this the \textit{local compact Runge-Kutta (cRK)} procedure.
\paragraph{cRK22.}
\begin{equation} \label{eq:crk22}
\F = \pf (\uu_h^{(2)}),
\end{equation}
where
\begin{align*}
\uu^{(2)}_h & = \uu_h^n - \frac{\mathLaplace t}{2} \dlocx
 \pf(\uu_h^n).
\end{align*}
\paragraph{cRK33.}

\begin{equation} \label{eq:crk33}
\F = \frac{1}{4}  \pf ( \uu_h^n ) + \frac{3}{4}  \pf( \uu_h^{(3)}),
\end{equation}
where
\begin{align*}
\uu^{(2)}_h & = \uu_h^n - \frac{\mathLaplace t}{3} \dlocx
 \pf ( \uu_h^n ), \\
\uu_h^{(3)} & =  \uu_h^n - \frac{2}{3} \mathLaplace t
\dlocx  \pf (\uu_h^{(2)}).
\end{align*}

\paragraph{cRK44.}

\begin{equation} \label{eq:crk44}
\F = \frac{1}{6}  \pf (\uu_h^n) + \frac{1}{3}  \pf (\uu_h^{(2)}) + \frac{1}{3}  \pf (\uu_h^{(3)}) + \frac{1}{6}  \pf (\uu_h^{(4)}),\\
\end{equation}
where
\begin{align*}
\uu_h^{(2)} & = \uu_h^n - \frac{\mathLaplace t}{2}
\dlocx  \pf (\uu_h^n),\\
\uu_h^{(3)} & = \uu_h^n - \frac{\mathLaplace t}{2}
\dlocx  \pf (\uu_h^{(2)}),\\
\uu_h^{(4)} & = \uu_h^n - \mathLaplace t
\dlocx  \pf (\uu_h^{(3)}).
\end{align*}

\subsection{Time average numerical flux computation} \label{sec:numflux}

As is seen from~\eqref{eq:crkfr.avg}, the cRKFR scheme requires computation of only one numerical flux $\F_\eph$, which is an approximation of the time average numerical flux at the interfaces (Appendix~\ref{app:time.averaged.flux}). The general form of time averaged numerical flux approximation discussed in this work is given by
\begin{equation} \label{eq:numflux.general}
\F_\eph = \frac{1}{2} (\F_\eph^- + \F_\eph^+) - \bD_\eph,
\end{equation}
where $\F_{\eph}^{\pm}$ is an approximation to the time average flux at the element interfaces. The $\F_{\eph}^\pm$ terms are called the central part of the numerical flux, and the $\bD_\eph$ is called the dissipative part. In the standard RKFR schemes, both the central and dissipative parts are computed using the solution at the current time level~\eqref{eq:num.flux.fr}.

\subsubsection{Dissipative part of the numerical flux} \label{sec:diss.models}

We now discuss the computation of the dissipative part $\bD_{\eph}$ of the numerical flux~\eqref{eq:numflux.general}.

\paragraph{D-CSX (Chen, Sun, Xing~\cite{chen2024}) dissipation.} As is mentioned in (3.14) of~\cite{chen2024}, the numerical flux in the cRKDG scheme of~\cite{chen2024} is computed as
\begin{equation} \label{eq:Feph.crk}
\F_\eph = \sum_{i=1}^s b_i \pf_\eph(\uu_h^{(i)}),
\end{equation}
where $\{b_i\}$ are the coefficients of the final stage of RK scheme~(\ref{eq:butcher}), $\{\uu_h^{i}\}$ are the inner stage evolutions of cRK scheme that are locally computed~\eqref{eq:crkfr1}. By the definition of $\dfrx$~\eqref{eq:dfrx.defn}, the approximation~\eqref{eq:Feph.crk} is also the numerical flux used in~\eqref{eq:crkfr2}. If the numerical flux in the cRKFR scheme~\eqref{eq:crkfr.avg} is computed as in~\eqref{eq:Feph.crk}, the obtained cRKFR scheme will be equivalent to~\eqref{eq:crkfr}. We can see this by observing that~\eqref{eq:crkfr.avg} and~\eqref{eq:crkfr} compute the same inner evolutions $\{\uu_h^{(i)}\}$ using the local operator $\dlocx$. The final stage in~\eqref{eq:crkfr2} applies $\dfrx$~\eqref{eq:dfrx.defn} on the $\{\pf(\uu_h^{(i)})\}$ and then  takes their linear combination, while the final stage in~\eqref{eq:crkfr.avg} applies $\dfrx$ on a linear combination of $\{\pf(\uu_h^{(i)})\}$. Since the only nonlinear part in an application of $\dfrx$ on $\{\pf(\uu_h^{(i)})\}$ is the numerical flux, the scheme~\eqref{eq:crkfr.avg} will be equivalent to~\eqref{eq:crkfr} if it uses the same numerical flux~\eqref{eq:Feph.crk}. Thus, in this case, the scheme~\eqref{eq:crkfr.avg} will also be equivalent to the cRKDG scheme of~\cite{chen2024} when the quadrature points of the DG scheme are taken to be the same as solution points (Section~\ref{sec:rkfr}).

Expanding equation~\eqref{eq:Feph.crk} using the numerical flux expression~\eqref{eq:num.flux.fr}, we obtain that the dissipative part $\bD_{\eph}$~\eqref{eq:numflux.general} of the numerical flux~\eqref{eq:Feph.crk} is
\begin{equation} \label{eq:dcsx}
\bD_\eph = \sum_{i=1}^s b_i \frac{\lambda^i_{\eph}}{2} (\uu_{\eph}^{(i)+} - \uu_{\eph}^{(i)-}),\qquad \uu_{\eph}^{(i)\pm} = \uu_h^{(i)}(x_{\eph}^\pm).
\end{equation}
Equation~\eqref{eq:dcsx} will be referred to as the D-CSX dissipation model. This dissipation model requires computation of wave speed estimates $\lambda^i_{\eph} = \lambda(\uu_{\eph}^{(i)-}, \uu_{\eph}^{(i)+})$~\eqref{eq:rusanov.wave.speed} and interelement communication of $\uu_{\eph}^{(i)\pm}$ for each inner stage $i$ of the cRK scheme.

\paragraph{D1 dissipation.}
The D1 dissipation model (terminology of~\cite{babbar2022}), used in the earlier works like~\cite{Qiu2005b}, computes the dissipative part $\bD_{\eph}$~\eqref{eq:numflux.general} using the solution at the current time level as
\begin{equation}  \label{eq:D1}
\bD_\eph = \frac{\lambda_{\eph}}{2} (\uu_\eph^+ - \uu_\eph^-),\qquad \uu_\eph^\pm = \uu_h(x_\eph^\pm).
\end{equation}

The D1 dissipation model uses only one wave speed estimate $\lambda_{\eph} = \lambda(\uu_{\eph}^-, \uu_{\eph}^+)$~\eqref{eq:rusanov.wave.speed}, requiring less interelement communication than the D-CSX dissipation. However, it leads to lower Courant–Friedrichs–Lewy (CFL) numbers compared to D-CSX dissipation, see Table 1.1 of~\cite{chen2024}.

\paragraph{D2 dissipation.}

The D2 dissipation model, proposed for the LWFR scheme in~\cite{babbar2022}, is given by
\begin{equation} \label{eq:D2}
\bD_\eph = \frac{\lambda_{\eph}}{2} (\uU_\eph^+ - \uU_\eph^-), \qquad \uU_\eph^\pm = \uU_h(x_{\eph}^\pm),
\end{equation}
\begin{equation}\label{eq:time.avg.sol}
\uU_h = \sum_{p=0}^N \uU_{e,p} \ell_p(\xi),\qquad \uU_{e,p} = \sum_{i=1}^s b_i {\uu}_{e,p}^{(i)}.
\end{equation}
The D2 dissipation model~\eqref{eq:D2} also requires only one wave speed estimate~$\lambda_{\eph} = \lambda(\uu_{\eph}^-, \uu_{\eph}^+)$~\eqref{eq:rusanov.wave.speed}, and interelement communication of $\uu_{\eph}^\pm, \uU_{\eph}^\pm$. The interelement communication of the D2 dissipation model is more than the D1 model, but still less than the D-CSX model especially as higher order RK methods are used. For linear problems, the D2 dissipation model and D-CSX model are equivalent and thus both have the same CFL numbers~\cite{babbar2022,chen2024}. For nonlinear problems, the two models are numerically observed to give very similar performance. Thus, in all numerical results, unless specified otherwise, D2 dissipation is used because it has higher CFL number than D1 dissipation and lower interelement communication than D-CSX dissipation.
\subsubsection{Central part of the numerical flux}

The $\F_\eph^\pm$ in~\eqref{eq:numflux.general} are approximations of the time averaged flux obtained from the right and left element of the interface $\eph$ respectively. A natural idea to obtain the fluxes at the interfaces is by extrapolating the discontinuous time averaged flux~\eqref{eq:disc.avg.flux} to obtain
\begin{equation} \label{eq:extrapolate}
\F_\eph^\pm = \F_h^\delta(x_\eph^\pm).
\end{equation}
In the context of Lax-Wendroff Flux Reconstruction schemes, the above flux approximation is called the \textit{Average and Extrapolate (\extrapolate)} procedure~\cite{babbar2022}. It was observed in~\cite{babbar2022} that the above procedure leads to suboptimal convergence rates and thus \textit{Extrapolate and Average (\evaluate)} procedure was introduced in~\cite{babbar2022} for Lax-Wendroff schemes. In the context of cRK, we propose the following as the {\evaluate} scheme
\begin{equation} \label{eq:evaluate}
\F_\eph^\pm = \sum_{i=1}^s b_i \pf(\uu_\eph^{(i)\pm}).
\end{equation}
The two schemes~(\ref{eq:extrapolate},~\ref{eq:evaluate}) are equivalent for linear problems. Similar to~\cite{babbar2022}, we see that the cRK scheme with {\evaluate} flux computation~\eqref{eq:evaluate} gives optimal order of accuracy for all problems, while {\extrapolate} are suboptimal for some nonlinear problems (Section~\ref{sec:burg}). The schemes are equivalent if the solution points~\eqref{eq:soln.points} are taken to include the boundary points, e.g., for Gauss-Legendre-Lobatto (GLL) solution points. We also remark that suboptimal convergence rates observed in Section~\ref{sec:burg} are also observed when GLL solution points are used.


\subsection{Conservation property and admissibility preservation}
To reliably compute weak solutions with discontinuities for non-linear conservation laws, it is essential to use conservative numerical schemes. According to the Lax-Wendroff theorem, any consistent and conservative numerical method that converges will give a weak solution in the limit. Although the conservative property of the cRKFR scheme~\eqref{eq:crkfr.avg} may not be immediately obvious, it can be obtained by multiplying~\eqref{eq:crkfr.avg} with the quadrature weights associated with the solution points and summing the results across all points in the $e^\text{th}$ element as
\begin{equation} \label{eq:upmean}
\au_e^{n+1} = \au_e^n - \frac{\Delta t}{\Delta x_e} ( \F_{\eph} - \F_{\emh}),\qquad \au_e^n = \sum_{p=0}^N \uu_{e,p} w_p.
\end{equation}
As will be discussed in Section~\ref{sec:flux.limiter}, the conservation property is also crucial for the procedure for enforcing admissibility of a solution~\eqref{eq:uad.form}. For the solution of a conservation law, admissibility preservation is the property that
\[
\uu(\cdot, t_0) \in \Uad \qquad  \implies \qquad \uu(\cdot, t) \in \Uad,\qquad t>t_0,
\]
and thus an admissibility preserving numerical scheme is defined as follows.
\begin{definition}\label{defn:admissibility.preserving}
A Flux Reconstruction scheme is said to be admissibility preserving if
\[
\uu_{e,p}^n \in \Uad \quad \forall e,p \qquad \implies \qquad \uu_{e,p}^{n+1} \in \Uad \quad \forall e,p,
\]
where $\Uad$ is the admissible set~\eqref{eq:uad.form} of the conservation law~\eqref{eq:con.law}.
\end{definition}
The relevance of~\eqref{eq:upmean} can already be motivated as our idea is to first develop a scheme that preserves the following weaker property of \textit{admissibility preservation in means}.
\begin{definition}\label{defn:admissibility.preserving.means}
A Flux Reconstruction scheme is said to be admissibility preserving in means if
\begin{equation}
\uu_{e,p}^n \in \Uad \quad \forall e,p \qquad \implies \qquad \au_{e}^{n+1} \in \Uad \quad \forall e,
\end{equation}
where $\Uad$ is the admissible set~\eqref{eq:uad.form} of the conservation law~\eqref{eq:con.law}.
\end{definition}
Once the scheme preserves admissibility in means, an admissibility preserving scheme (Definition~\ref{defn:admissibility.preserving}) is obtained by using the scaling limiter of Zhang and Shu~\cite{Zhang2010b}.

\subsection{High order boundary treatment} \label{sec:boundary}
Similar to the boundary condition treatment for LWFR schemes in~\cite{babbar2022} (also detailed in Chapter 4 of~\cite{babbar2024thesis}), we impose boundary conditions for the cRKFR scheme weakly by specifying the numerical flux~\eqref{eq:numflux.general} at the boundary faces. In this section, we describe how the time averaged numerical flux~\eqref{eq:numflux.general} for the boundaries, called the \textit{boundary numerical flux}, is computed for the cRKFR scheme. In all numerical experiments, we computed the boundary numerical flux by performing computations only at the boundary face, which is described for various types of boundary conditions in this section. We also propose an alternate approach that uses the local cRK procedure (Section~\ref{sec:crkfr}) in a \textit{ghost element}, which also gives the optimal order of accuracy. Since this work is limited to Cartesian grids, the former approach is more computationally efficient as it only requires ghost values at the interfaces. However, we also describe the latter approach as it has the potential to be more general, e.g., extendable to high order curvilinear meshes for some types of boundary conditions like solid walls. Curvilinear meshes for cRKFR schemes will be considered in future works.
Consider the grid elements to be $\Omega_e = [x_{\emh}, x_{\eph}]$ for $e=1,\dots, \nel$. Thus, the boundary elements correspond to $e = 1, \nel$. The numerical flux at the left most boundary is denoted by $\F_\half$ and at the right boundary by $\F_{\nel + \half}$. In the following, we describe the procedure to compute the numerical flux $\F_\half$ at the left boundary $x_\half$. The procedure for the right boundary $x_{\nel + \half}$ is similar. In order to compute~\eqref{eq:numflux.general}, we require $\uu_{\half}^-, \F_\half^-, \uU_\half^-$, which we will refer to as the \textit{interface ghost values}.

\paragraph{Periodic boundaries.} In case of periodic boundary conditions, the interface ghost values are computed as
\[
\uu_{\half}^- = \uu_{\nel + \half}^-, \qquad \F_{\half}^- = \F_{\nel + \half}^-, \qquad \uU_{\half}^- = \uU_{\nel + \half}^-.
\]
\paragraph{Inflow/outflow boundaries.} The boundary conditions will be a mix of inflow and outflow boundary conditions when eigenvalues of the flux Jacobian $\pf'$ of~\eqref{eq:con.law} at interface $x_\half$ are both positive and negative. Such a case is handled through an upwind numerical flux like HLLC/Roe's (see Appendix A of~\cite{babbar2022} for computation of these numerical fluxes from time averages fluxes). We assume that the inflow boundary value is specified by $\pg(x_\half, t)$ for $t \ge 0$. Then, the interface ghost values are given by
\[
\uu_{\half}^- = \pg(x_\half, t^n),\qquad \F_{\half}^- = \sum_{i=1}^s b_i \pf(\pg(x_\half, t^n + \Delta t c_i)), \qquad \uU_\half^- = \sum_{i=1}^s b_i \pg(x_\half, t^n + \Delta t c_i),
\]
where $b_i, c_i$ are as defined in~\eqref{eq:butcher}. With these interface ghost values, an upwind numerical flux is used to obtain the boundary numerical flux. An HLLC flux was used at the left boundary in the jet flow simulation of Section~\ref{sec:m2000} because the left boundary is a mix of inflow and outflow. In case the boundary is fully inflow (all eigenvalues are positive\footnote{For the right boundary $x_{\nel+\half}$, inflow boundaries are those where $\pf'$~\eqref{eq:con.law} has only negative eigenvalues. In general, for higher dimensions, inflow boundaries are those where eigenvalues of the Jacobian of the normal flux are positive.}), we simply compute the numerical flux as $\F_\half = \F_\half^-$, while it is computed as $\F_\half = \F_\half^+$ in the fully outflow case (all eigenvalues are negative). This approach gives optimal convergence rates as shown in Section~\ref{sec:var.adv}.
\paragraph{Solid wall boundaries.} In case of solid wall boundary conditions for Euler's equations~\eqref{eq:1deuler}, where the conservative variables are $\uu = (\rho, \rho v_1, E)$, the solution at the ghost interface is prescribed by reflecting the velocity. Since the flux is given by $\pf = (\rho v_1, p + \rho v_1^2,(E+p)v_1)$, the solid wall boundary conditions are enforced by specifying the ghost values as
\begin{equation} \label{eq:reflect.time.avg}
\uu_{\half}^- = ((\uu_{\half}^+)_1, -(\uu_{\half}^+)_2, (\uu_{\half}^+)_3), \ \ \F_{\half}^- = (-(\F_{\half}^+)_1, (\F_{\half}^+)_2, -(\F_{\half}^+)_3),  \ \  \uU_{\half}^- = ((\uU_{\half}^+)_1, -(\uU_{\half}^+)_2, (\uU_{\half}^+)_3).
\end{equation}
\subsubsection{Local cRK procedure in ghost elements (Alternate boundary treatment)}
An alternative approach for boundary condition treatment is to perform the local cRK procedure (Section~\ref{sec:crkfr}) in the ghost elements. The idea is to consider $e=0,\nel+1$ as the ghost elements with the relevant interfaces being $x_{\half}, x_{\nel+\half}$ respectively. We explain it by considering the left ghost cell $\Omega_0 = [x_{-\half}, x_\half]$, where $x_{-\half} = x_\half - \Delta x_1$. Similar to~\eqref{eq:soln.points}, we can place $N+1$ solution points in the ghost cell and define the numerical solution $\uu_{0,p}^n$ in the ghost element $\Omega_e$. The particular choices for the solution values will depend on the boundary condition. In case of periodic boundary conditions, the solution values are specified as $\uu_{0, p} = \uu_{\nel, p}$, for $0 \le p \le N$. In case of inflow/outflow boundaries with inflow part $\pg(x, t)$ for $x \in \Omega_{0}$, we prescribe the boundary values as $\uu_{0,p} = \pg(x_{0,p}, t^n)$. In case of solid wall boundary conditions, the solution values are specified as
\begin{equation} \label{eq:reflect.bc}
\uu_{0, p} = \uu^\text{reflect}_{1, N-p+1},\qquad \uu^\text{reflect}_{1, p} = (\rho, -\rho v_1, E)_{1,p}.
\end{equation}
Once the solution values are available in the ghost elements, the local cRK procedure of Section~\ref{sec:crkfr} is performed to compute the time averaged flux $\F_{\half}^-$ and solution $\uU_{\half}^-$ (\ref{eq:disc.avg.flux},~\ref{eq:time.avg.sol}). This gives us all the quantities needed to compute the numerical flux at the left boundary $\F_{\half}$ using~\eqref{eq:numflux.general}. The reason why this boundary treatment is more general is that equation~\eqref{eq:reflect.bc} can be extended to curvilinear boundaries by reflecting the solution along the normal direction. However, it is not possible to do the same for the time average flux and time average solution in~\eqref{eq:reflect.time.avg}. Since this work is restricted to Cartesian grids, we use the approach of~\eqref{eq:reflect.time.avg} in this work. The alternate will be considered in future works when the cRKFR scheme is extended to curvilinear meshes.

\section{Blending scheme} \label{sec:blending.scheme}
By the formal high order accuracy of the RK scheme (Appendix~\ref{app:time.averaged.flux}), the cRKFR scheme~\eqref{eq:crkfr.avg} is a fully discrete high order accurate scheme. Thus, as demonstrated in the numerical results (like Section~\ref{sec:burg}), in problems with smooth solutions, the scheme gives high order accuracy. In this section, we discuss how to apply this high order method to non-smooth solutions of hyperbolic conservation laws~\eqref{eq:con.law} where structures like rarefactions, shocks and other discontinuities appear. Such solutions occur in a wide variety of practical applications. A direct usage of high order methods for such problems is known to produce Gibbs oscillations, as is also formalized in Godunov's order barrier theorem~\cite{godunov1959}. The approach that we take in this work is to adaptively switch to a robust lower order method in spatial regions where the solution is less regular, while using the high order scheme in regions of smoothness. This is consistent with Godunov's order barrier theorem because the obtained scheme will be non-linear even for linear equations. The switching process is called the blending scheme because it takes a combination of high and low order methods. In~\cite{babbar2024admissibility}, a subcell based blending for the high order
single stage Lax-Wendroff Flux Reconstruction was introduced, inspired by~\cite{hennemann2021}. A crucial part of the subcell scheme of~\cite{babbar2024admissibility} was the limiting of the time average numerical flux. Since we have cast the cRK scheme in terms of time average fluxes~\eqref{eq:crkfr.avg}, we will now be able to apply the ideas of~\cite{babbar2024admissibility} to obtain a subcell based cRK scheme with blending. Similar to~\cite{babbar2024admissibility}, the limiter is applied only once per time step in the cRKFR scheme.

The description of the blending scheme begins by writing a single evolution step of the cRKFR scheme~\eqref{eq:crkfr.avg} as
\begin{equation*}
\vu^{H, n + 1}_e = \vu^n_e - \frac{\Delta t}{\Delta x_e}  \vR^H_e,
\end{equation*}
where $\vu_e$ is the vector of nodal values in the element $e$.  We use a lower order scheme written in the same form as
\begin{equation*}
\vu^{L, n + 1}_e = \vu^n_e - \frac{\Delta t}{\Delta x_e}  \vR^L_e.
\end{equation*}
Then the blended scheme is obtained by taking a convex combination of the two schemes as
\begin{equation}
\begin{split}
\vu^{n + 1}_e & = (1 - \alpha_e)  \vu^{H, n + 1}_e + \alpha_e  \vu^{L, n
+ 1}_e = \vu^n_e - \frac{\Delta t}{\Delta x_e}  [(1 - \alpha_e) \vR^H_e + \alpha_e  \vR^L_e],
\end{split} \label{eq:blended.scheme}
\end{equation}
where $\alpha_e \in [0, 1]$ is chosen based on a local smoothness
indicator taken from~\cite{hennemann2021} (Also see Section 3.3 of~\cite{babbar2024admissibility} for a discussion on the indicator of~\cite{hennemann2021}). If $\alpha_e = 0$, the blending scheme is the high order cRKFR scheme, while
if $\alpha_e = 1$, the scheme becomes the low order scheme that does not produce spurious oscillations and is admissibility preserving (Definition~\ref{defn:admissibility.preserving}). In subsequent sections, we explain the details of the lower order
scheme. The lower order scheme will
either be a first order finite volume scheme or a high resolution scheme based
on MUSCL-Hancock idea.
\begin{figure}
\centering
\includegraphics[width=0.7\textwidth]{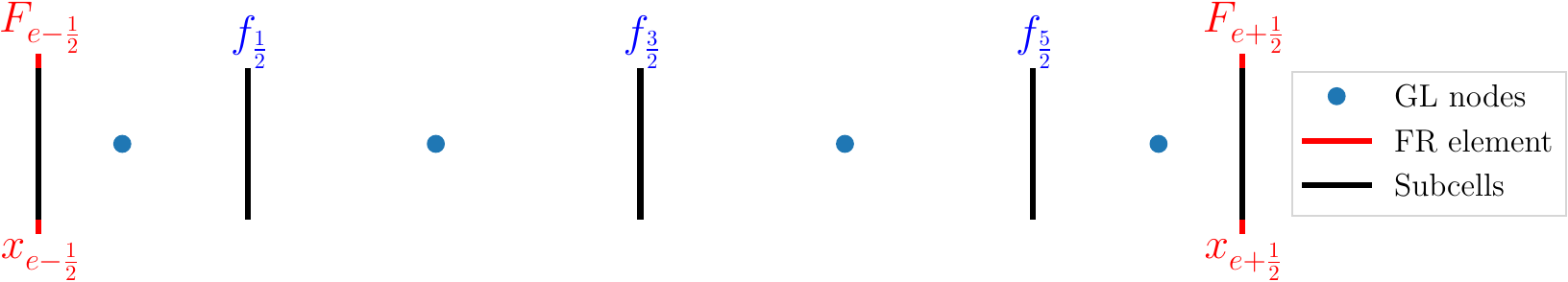} \\
(a) \\
\includegraphics[width=0.7\textwidth]{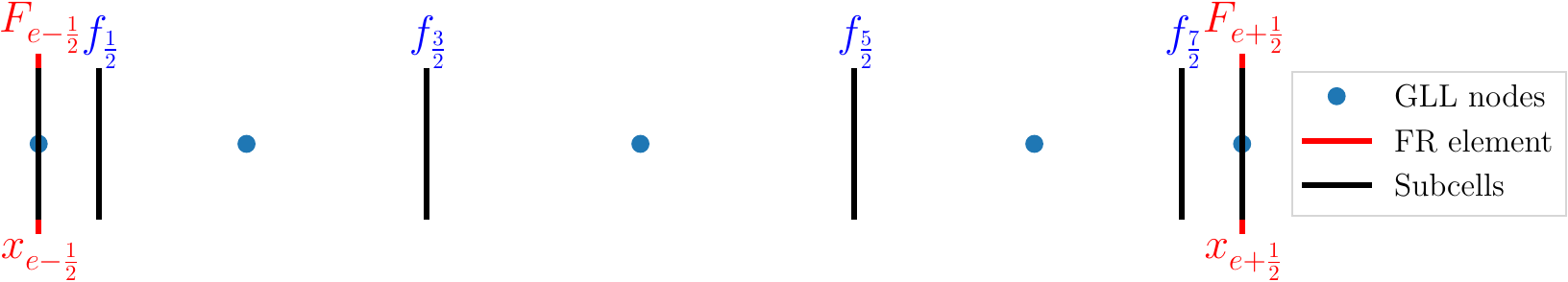} \\
(b)
\caption{Subcells used by lower order scheme for degree $N = 3$ using (a) Gauss-Legendre (GL) solution points (GL), (b) Gauss-Legendre-Lobatto (GLL) solution points. \label{fig:subcells}}
\end{figure}
Let us subdivide each element $\Omega_e$ into $N + 1$ subcells associated
to the solution points $\{ \xep, p = 0, 1, \ldots, N\}$ of the cRKFR scheme.
Thus, we will have $N + 2$ subfaces denoted by $\{x^e_{\pph}, p = - 1, 0,
\ldots, N\}$ with $x^e_{- \half} = x_{\emh}$ and $x^e_{\Nph} = x_{\eph}$. For
maintaining a conservative scheme, the $p^{\tmop{th}}$ subcell is chosen so
that
\begin{equation}
\label{eq:subcell.defn} x_{\pph}^e - x_{\pmh}^e = w_p \Delta x_e, \qquad 0
\le p \le N,
\end{equation}
where $w_p$ is the $p^{\text{th}}$ quadrature weight associated with the
solution points. Figure~\ref{fig:subcells} gives an illustration of the
subcells for the degree $N = 3$ case for Gauss-Legendre and Gauss-Legendre-Lobatto solution points. The low order scheme is obtained by updating the solution in each
of the subcells by a finite volume scheme,
\begin{subequations}
\begin{align}
\uez^{L, n + 1} & = \uez^n - \frac{\Delta t}{w_0 \Delta x_e}
[\pf_{\half}^e - \F_{\emh}], \label{eq:low.order.update.a}\\
\uep^{L, n + 1} & = \uep^n - \frac{\Delta t}{w_p \Delta x_e}
[\pf_{\pph}^e - \pf_{\pmh}^e], \qquad 1 \le p \le N - 1, \\
\ueN^{L, n + 1} & = \ueN^n - \frac{\Delta t}{w_N \Delta x_e}  [\F_{\eph} -
\pf_{\Nmh}^e] \label{eq:low.order.update.b}.
\end{align} \label{eq:low.order.update}
\end{subequations}
The inter-element fluxes $\F_{\eph}$ used in the low order scheme are same as
those used in the high order cRK scheme in equation~\eqref{eq:crkfr.avg}.
Usually, Rusanov's flux~\eqref{eq:rusanov.wave.speed} will be used for the
inter-element fluxes and in the lower order scheme. The element mean value
obtained by the low order scheme satisfies
\begin{equation}
\label{eq:low.order.cell.avg.update} \overline{\uu}_e^{L, n + 1} = \sum_{p =
0}^N \uep^{L, n + 1} w_p = \overline{\uu}_e^n - \frac{\Delta t}{\Delta x_e}
(\F_{\eph} - \F_{\emh}),
\end{equation}
which is identical to the element mean evolution of the cRKFR scheme~\eqref{eq:upmean} in the absence of source terms. Thus, the element mean in the blended scheme evolves according to
\begin{equation}
\label{eq:blended.cell.avg.update}
\begin{split}
\overline{\uu}_e^{n + 1} & = (1 - \alpha_e)  (\overline{\uu}_e)^{H, n + 1}
+ \alpha_e  (\overline{\uu}_e)^{L, n + 1}\\
& = (1 - \alpha_e)  \left[ \overline{\uu}_e^n - \frac{\Delta t}{\Delta
x_e} (\F_{\eph} - \F_{\emh}) \right] + \alpha_e  \left[ \overline{\uu}_e^n -
\frac{\Delta t}{\Delta x_e} (\F_{\eph} - \F_{\emh}) \right]\\
& = \overline{\uu}_e^n - \frac{\Delta t}{\Delta x_e}  (\F_{\eph} -
\F_{\emh}),
\end{split}
\end{equation}
and hence the blended scheme is also conservative. Overall, all three schemes, i.e., lower order, cRKFR and the blended scheme, predict the same mean value. The same inter-element flux $\F_{\eph}$ is used both in the low and high order schemes. To achieve high order accuracy, this flux needs to be high order accurate in smooth regions. However, usage of a high order flux may produce spurious oscillations near discontinuities when used in the low order scheme. A natural choice to balance accuracy and oscillations is to take
\begin{equation}
\F_{\eph} = (1 - \alpha_{\eph})  \F_{\eph}^\text{HO} + \alpha_{\eph}  \pf_{\eph},
\label{eq:Fblend}
\end{equation}
where $\F^\text{HO}_{\eph}$ is the high
order inter-element time-averaged numerical fluxes used in the cRKFR
scheme~\eqref{eq:numflux.general} and $\pf_{\eph}$ is a lower order flux at the face
$x_{\eph}$ shared between FR elements and
subcells~(\ref{eq:fo.at.face},~\ref{eq:mh.at.face}). The blending coefficient
$\alpha_{\eph}$ will be based on a local smoothness indicator which will bias
the flux towards the lower order flux $\pf_{\eph}$ near regions of lower
solution smoothness. However, to enforce admissibility in means
(Definition~\ref{defn:admissibility.preserving.means}), the flux has to be further limited, as
explained in Section~\ref{sec:flux.limiter}.

\subsection{First order blending}\label{sec:fo}

The lower order scheme is taken to be a first order finite volume scheme, for
which the subcell fluxes in~\eqref{eq:low.order.update} are given by
\[ \pf_{\pph}^e = \pf (\uep, \uu_{e, p + 1}). \]
At the interfaces that are shared with FR elements, we define the lower order
flux used in computing inter-element flux as
\begin{equation}
\pf_{\eph} = \pf (\ueN, \uepoz) \label{eq:fo.at.face}.
\end{equation}
In this work, the numerical flux $\pf (\cdot, \cdot)$ is taken to be
Rusanov's flux~\eqref{eq:rusanov.wave.speed}.
\subsection{Higher order blending}\label{sec:mh}
The MUSCL-Hancock scheme is a second order finite volume method that uses a flux whose stencil at FR element $e+1/2$ is given by
\begin{equation}
\pf_{\eph} = \pf (\uu_{e, N - 1}, \ueN, \uepoz, \uu_{e + 1, 1}).
\label{eq:mh.at.face}
\end{equation}
The standard MUSCL-Hancock scheme is described for finite volume grids where the solution points are the cell centres. In order to maintain the conservation property, the solution points are chosen to be as in~\eqref{eq:subcell.defn} and are thus not at cell centres. In order for the admissibility preservation for the blended scheme, it is also essential that the MUSCL-Hancock scheme used on subcells is admissibility preserving. Thus, we use the admissibility preserving MUSCL-Hancock scheme on non-uniform, non-centred grids described in Appendix A of~\cite{babbar2024admissibility}. The idea for admissibility preservation of the MUSCL-Hancock scheme in~\cite{babbar2024admissibility} is an extension of the cell-centred case of~\cite{Berthon2015}.
\subsection{Flux limiter for admissibility in means}\label{sec:flux.limiter}
The blending scheme~\eqref{eq:blended.scheme} with the interelement flux $\F_\eph$ chosen as in~\eqref{eq:Fblend} provides a reasonable control on spurious oscillations while maintaining accuracy.
However, since it allows some spurious oscillations, it does not guarantee that the numerical solution will be in the admissible set $\Uad$~\eqref{eq:uad.form}.
A standard approach to construct admissibility preserving schemes (Definition~\ref{defn:admissibility.preserving}) is to use Strong Stability Preserving (SSP) Runge-Kutta (RK) schemes~\cite{Gottlieb2001} for temporal discretization of FR schemes as they give admissibility preserving in means schemes (Definition~\ref{defn:admissibility.preserving.means}).
The admissibility preservation in means of SSP-RKFR schemes is obtained by decomposing the element mean evolutions into a convex combination of first order finite volume evolutions~\cite{Zhang2010b,babbar2024generalized}, and then exploiting the fact that each RK stage of an SSPRK method is a convex combination of forward Euler evolutions of the previous stages. For such admissibility preserving in means schemes, we can use the scaling limiter of Zhang and Shu~\cite{Zhang2010b} to obtain an admissibility preserving scheme (Definition~\ref{defn:admissibility.preserving}).
However, unlike the SSP-RKFR schemes, the general cRKFR schemes of the form~\eqref{eq:crkfr.avg} are designed in the Butcher form with no obvious SSP structures. In the recent work of~\cite{liu2024}, for the cRKDG schemes, which can be equivalent cRKFR schemes of the form~\eqref{eq:crkfr} (or of the form~\eqref{eq:crkfr.avg} with D-CSX dissipation~\eqref{eq:dcsx}), admissibility in means of each stage of the cRK method was proven with conditions on CFL numbers. Thus, by applying the scaling limiter of~\cite{Zhang2010b} at every stage in the cRK scheme, an admissibility preserving cRKDG scheme was obtained. This approach requires $s$ applications of scaling limiters of~\cite{Zhang2010b} for an $s$ stage Runge-Kutta method and cannot directly be applied to our D1/D2 dissipation models, which only use one wave speed estimate for each time level.

Therefore, to ensure admissibility preserving in means for the cRKFR scheme, also as an alternative to~\cite{liu2024}, we propose to use the flux limiter of~\cite{babbar2024admissibility} to compute the time average numerical flux used at the element interfaces. This approach can be applied since we have cast our scheme in terms of time averaged flux~\eqref{eq:crkfr.avg}. For any number of stages of the cRK method, the process requires only one flux limiter to obtain admissibility preservation in means of the cRKFR method. In addition, the scheme helps to obtain admissibility in means for the general cRKFR schemes of the form~\eqref{eq:crkfr.avg}. The theoretical basis for the flux limiter is the same as for the LWFR scheme (Theorem 2 of~\cite{babbar2024admissibility}), and is given in the following theorem.
\begin{theorem}\label{thm:lwfr.admissibility}
Consider the cRKFR blending scheme~\eqref{eq:blended.scheme} where low and high order schemes use the same numerical flux $\F_\eph$ at every element interface. Then we obtain the following results on  admissibility preserving in means (Definition~\ref{defn:admissibility.preserving.means}) of the scheme:
\begin{enumerate}
\item Suppose  element means of both low and high order schemes are equal. The blended scheme~\eqref{eq:blended.scheme} is admissibility preserving in means if and only if the lower order scheme is admissibility preserving in means;
\item If the finite volume method using the lower order flux $\pf_\eph$ as the interface flux is admissibility preserving, and the blended numerical flux $\F_\eph$ is chosen to preserve the admissibility of lower-order updates at solution points adjacent to the interfaces, then the blending scheme~\eqref{eq:blended.scheme} will preserve admissibility in means.
\end{enumerate}
\end{theorem}
By Theorem~\ref{thm:lwfr.admissibility}, an admissibility preserving in means cRKFR scheme (Definition~\ref{defn:admissibility.preserving.means}) can be obtained if we choose the flux $\F_\eph$ so that ~(\ref{eq:low.order.update.a},~\ref{eq:low.order.update.b}) are admissible. The idea is to update the coefficient $\alpha_{\eph}$ in~\eqref{eq:Fblend} to ensure that all the lower order evolutions in~\eqref{eq:low.order.update} are admissible. The procedure is the same as for LWFR scheme, and can be found in Section 5 of~\cite{babbar2024admissibility}, and also as Algorithm 4 in~\cite{babbar2025}. We also remark that, by generalizing the work of~\cite{Zhang2010b}, the ideas of flux limiter apply even when a blending scheme is not used~\cite{babbar2024generalized}. The limiting process has been extended to curvilinear grids in~\cite{babbar2025} and to multiderivative Runge-Kutta schemes in~\cite{babbar2024multiderivative}.

\section{Numerical results} \label{sec:numerical.results}
We show various numerical results to demonstrate that the cRKFR scheme is high order accurate and that the cRKFR scheme with blending is robust and admissibility preserving. We first show a few results for comparison between D1/D2 dissipation~(\ref{eq:D1},~\ref{eq:D2}), and between \extrapolate/{\evaluate} fluxes~(\ref{eq:extrapolate},~\ref{eq:evaluate}). The performance of D-CSX dissipation model~\eqref{eq:dcsx} is not shown as it is virtually the same as the D2 model although a script verifying that is provided for Burgers' equation in the file~\href{https://github.com/Arpit-Babbar/paper-crkfr/blob/main/generate/dcsx.jl}{generate/dcsx.jl} at~\cite{crkrepo}. After the comparison, we use the D2 dissipation for its higher CFL stability (Table 1 of~\cite{babbar2022}) and the {\evaluate} flux for its higher accuracy in all results. When a degree $N$ polynomial discretization is used, the cRK scheme of order $N+1$ is used (e.g., cRK33~\eqref{eq:crk33} for $N=2$). The most general equations that we work with are the 2-D hyperbolic conservation laws
\begin{equation} \label{eq:2d.con.law}
\uu_t + \pf_x + \pg_y = 0,
\end{equation}
where $\uu \in \re^\nc$ is the vector of conserved variables and $\pf, \pg$ are $x,y$ fluxes. The description of the numerical scheme in previous sections has been for 1-D conservation laws~\eqref{eq:con.law}, but the extension to 2-D conservation laws~\eqref{eq:2d.con.law} can be done by applying the similar procedure to both $x,y$ fluxes in~\eqref{eq:2d.con.law} (see~\cite{basak2024} for a detailed description of the LWFR scheme with blending limiter in 2-D). Recall $\sigma(A)$ denotes the spectral radius of a matrix $A$. We compute the time step size $\Delta t$ for the scheme with polynomial degree $N$ as
\begin{align} \label{eq:time.step.2d}
\Delta t = C_s \min_e \left( \frac{\sigma(\ff'(\au_e))}{\Delta x_e} + \frac{\sigma(\fg'(\au_e))}{\Delta y_e} \right)^{-1} \text{CFL}(N),
\end{align}
where $e$ is the element index, $\Delta x_e, \Delta y_e$ are the element lengths in $x,y$ directions respectively,  $C_s \le 1$ is a safety factor, $\text{CFL}(N)$ is the optimal CFL number obtained by Fourier stability analysis. The safety factor is taken to be $C_s = 0.98$ unless specified otherwise. Since the cRKFR scheme is linearly equivalent to the LWFR scheme (Appendix~\ref{app:time.averaged.flux}), we refer to Table 1 of~\cite{babbar2022} or Table 1.1 of~\cite{chen2024} for the CFL numbers. For the 1-D equations,~\eqref{eq:time.step.2d} is used by setting $\fg = \bzero$. The developments have been contributed to the Julia~\cite{Bezanson2017} package \texttt{Tenkai.jl}~\cite{tenkai} and the run files for generating the results in this paper are available at~\cite{crkrepo}.
\subsection{Scalar equations}
In this section, we look at numerical results involving scalar equations. The $L^2$ error is plotted  versus degrees of freedom with different polynomial degrees in different colours. In all plots, degrees $N=1,2,3$ are shown with orange, blue and green colours respectively.
\subsubsection{Linear advection equation}

\begin{figure}
\centering
\begin{tabular}{cc}
\includegraphics[width=0.45\textwidth]{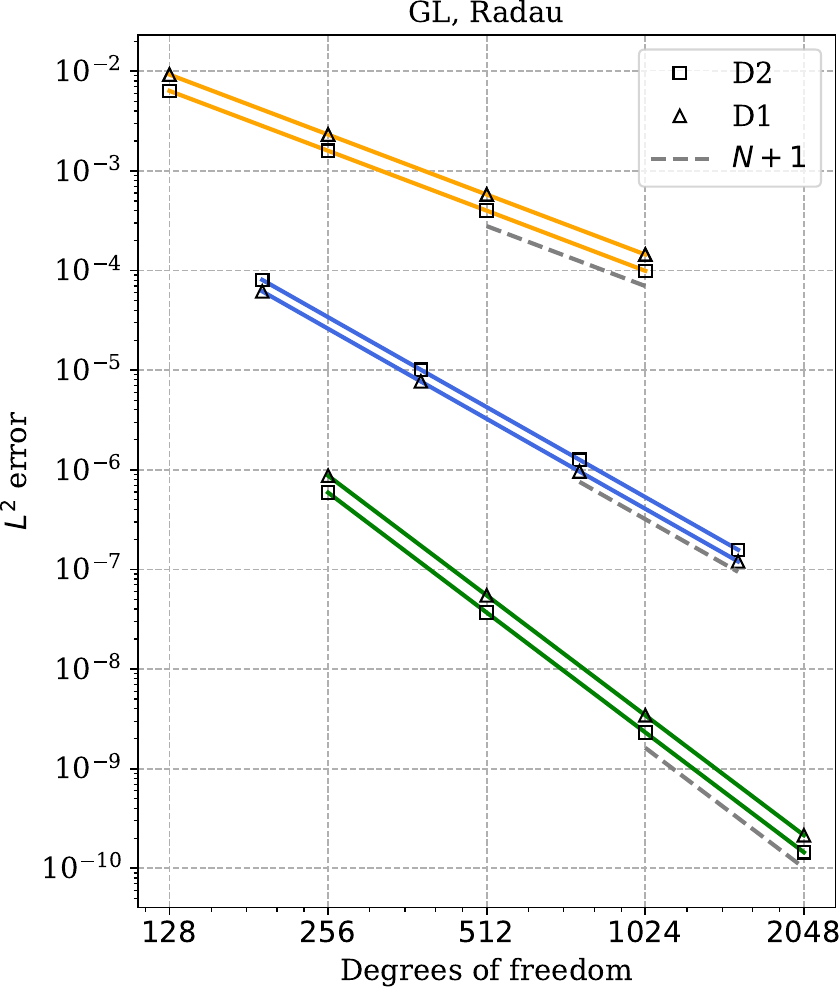} &
\includegraphics[width=0.45\textwidth]{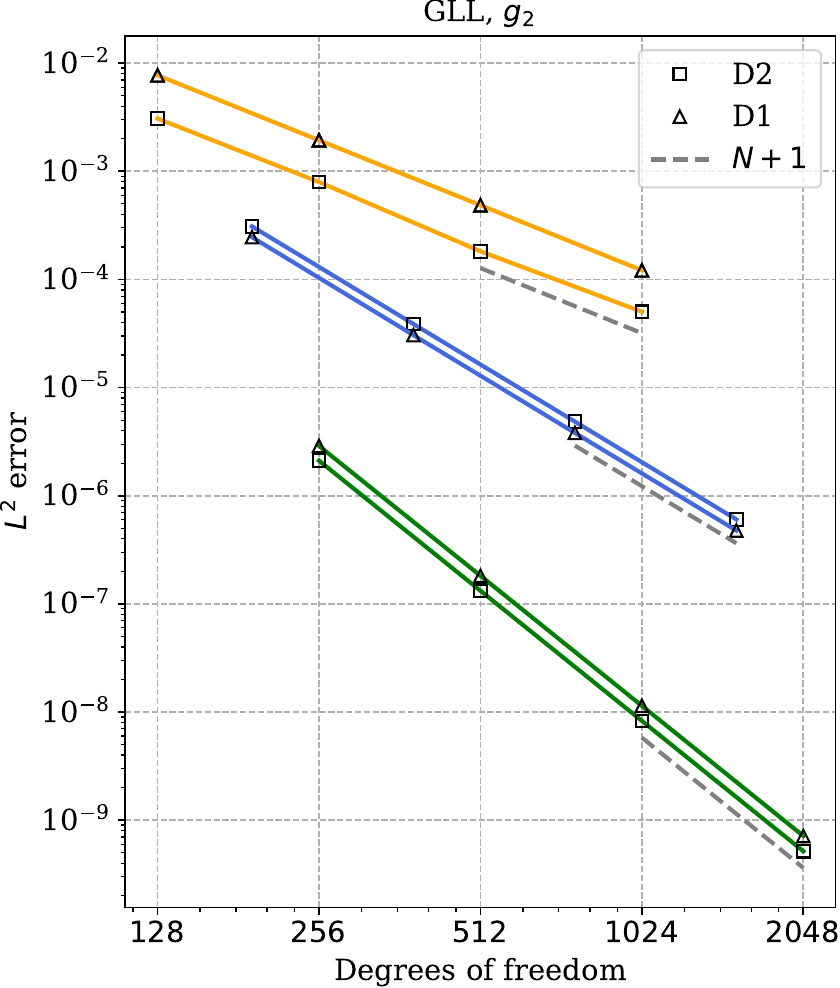} \\
(a) & (b)
\end{tabular}
\caption{$L^2$ error versus number of degrees of freedom for degrees $N=1,2,3$ for 1-D linear advection equation with periodic boundary conditions at $t = 2$ comparing D1 and D2 dissipation schemes using (a) Gauss-Legendre (GL) solution points with Radau correction functions, (b) Gauss-Legendre-Lobatto (GLL) solution points with $g_2$ correction functions.}
\label{fig:linadv1d}
\end{figure}

We consider the scalar 1-D conservation law~\eqref{eq:con.law} with the flux $f(u) = u$ and the initial condition $u(x,0) = \sin(2\pi x)$ with periodic boundary conditions on the domain $[0,1]$. The $L^2$ error is computed at $t=2$ for polynomial degrees $N=1,2,3$ with the respective order $N+1$ cRK method~(\ref{eq:crk22},~\ref{eq:crk33},~\ref{eq:crk44}). The $L^2$ error versus total degrees of freedom plots are shown in Figure~\ref{fig:linadv1d} comparing D1 and D2 dissipation for Radau and $g_2$ correction functions. The D1 and D2 dissipation schemes are used at their optimal CFL numbers (Table 1 of~\cite{babbar2022}). Optimal order of accuracy is seen for both correction functions. Radau correction function give lower errors, as is consistent with observations in the earlier works~\cite{Huynh2007,babbar2022}. For odd degrees $N=1,3$, the D2 dissipation scheme is showing better accuracy which is promising since the scheme is running at larger time steps. For the even degree $N=2$, the D1 dissipation has smaller errors although we also observed that D2 dissipation could give similar accuracy if it is using the same time step sizes as that of the D1 dissipation.

\subsubsection{Variable advection equation} \label{sec:var.adv}
\begin{figure}
\centering
\begin{tabular}{ccc}
\includegraphics[width=0.305\textwidth]{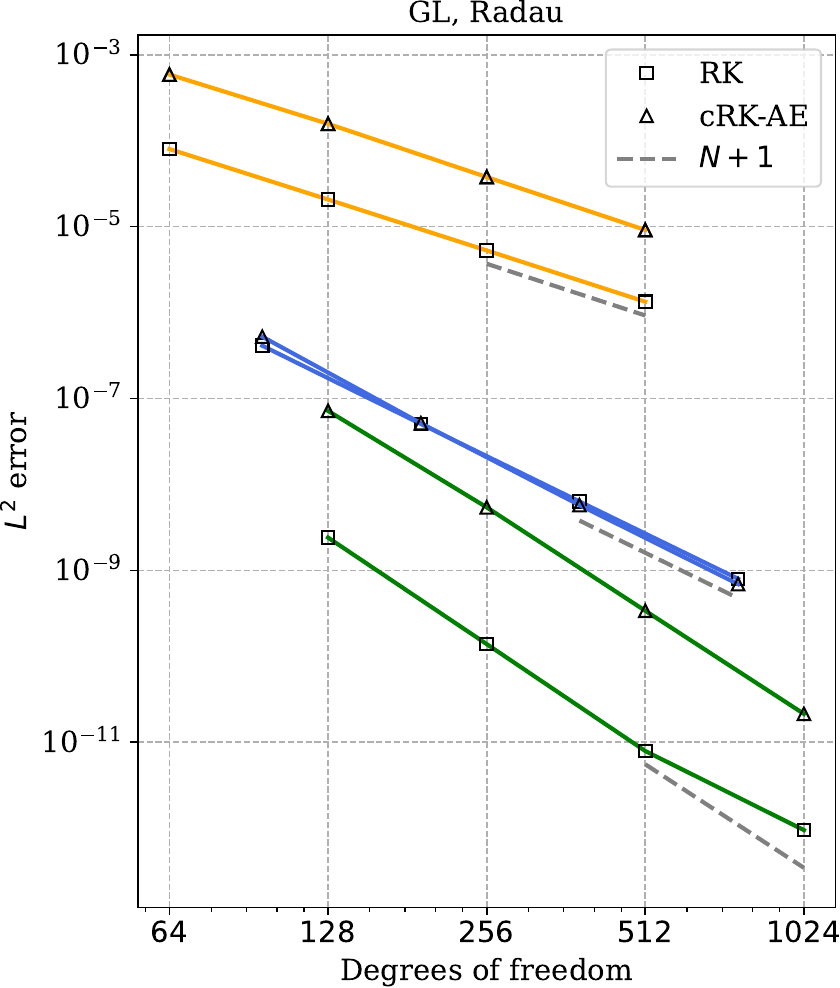} &
\includegraphics[width=0.305\textwidth]{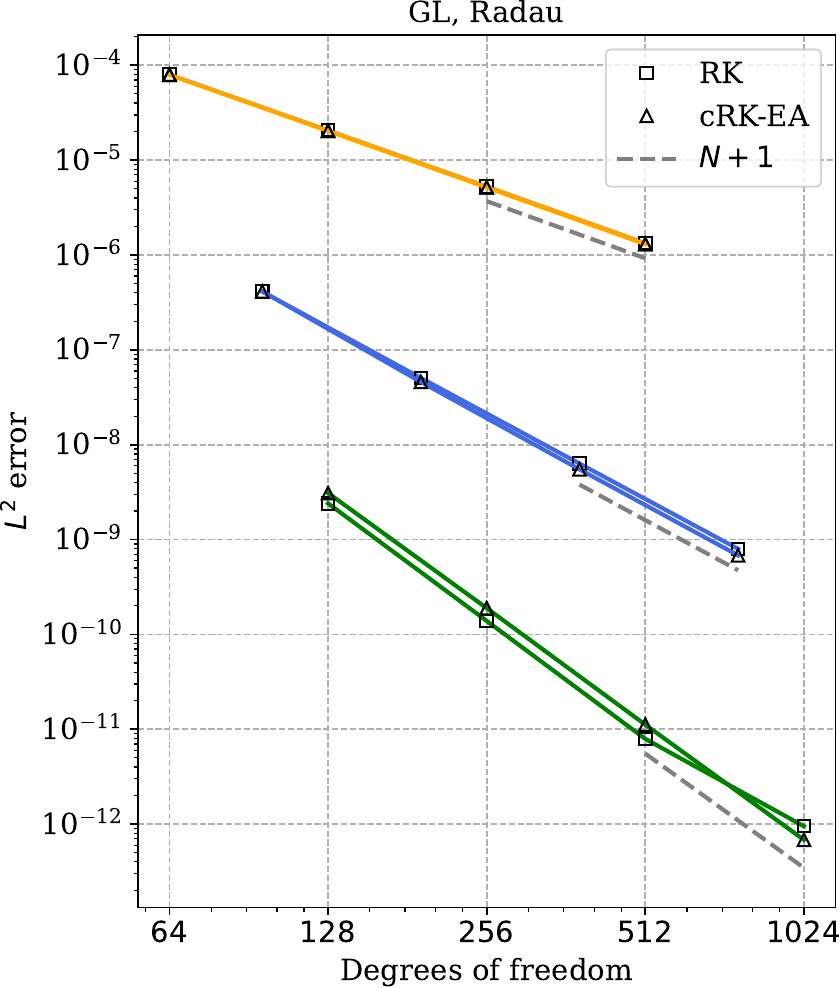} &
\includegraphics[width=0.305\textwidth]{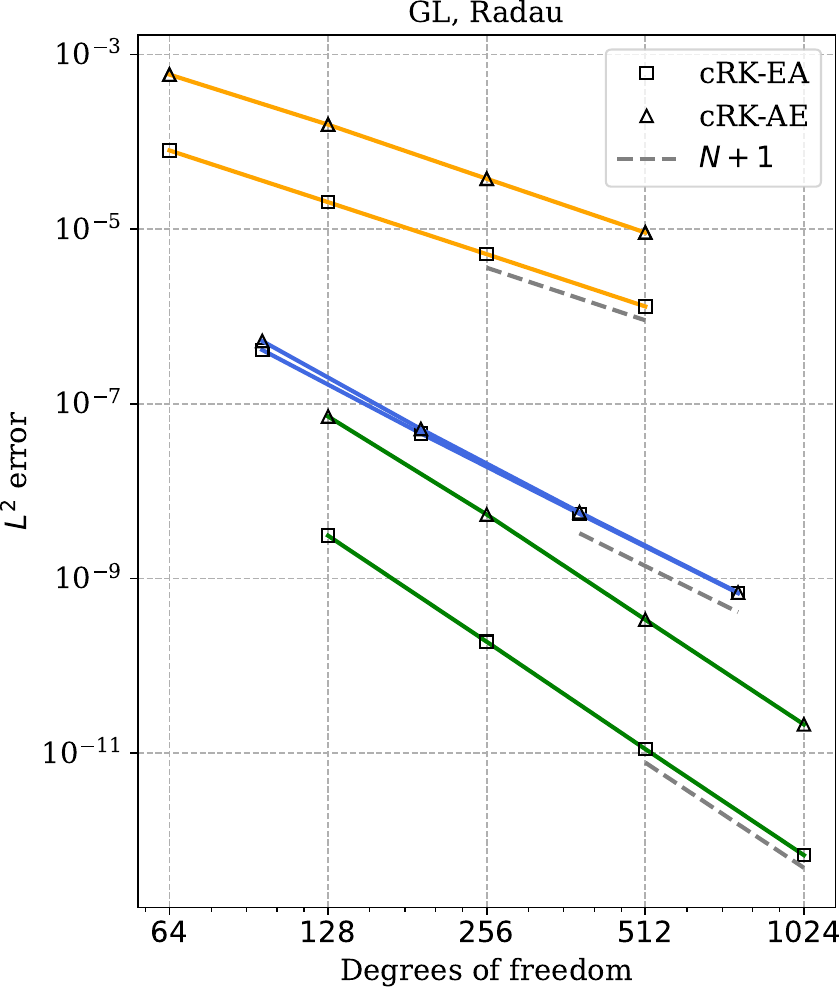} \\
(a) & (b) & (c)
\end{tabular}
\caption{$L^2$ error versus number of degrees of freedom for degrees $N=1,2,3$ for variable advection equation with flux $f(x,u) = x^2 u$ and non-periodic boundary conditions using Gauss-Legendre (GL) solution points with Radau correction functions. (a) Comparison of the cRKFR scheme using the {\extrapolate} flux~\eqref{eq:extrapolate} with the standard RKFR scheme, (b) comparison of the cRKFR scheme using the {\evaluate} flux~\eqref{eq:evaluate} with the standard RKFR scheme, (c) comparison of the cRKFR schemes using {\extrapolate}~\eqref{eq:extrapolate} and {\evaluate}~\eqref{eq:evaluate} fluxes.}
\label{fig:varadv1d}
\end{figure}

We now work with conservation laws where the flux has an explicit spatial dependence so that the PDE is given by
\[
u_t + f(x,u)_x = 0, \qquad f(x,u) = a(x) u.
\]
We tried various such test cases from~\cite{Offner2019} and found similar conclusions. To save the space, we only present the case $a(x) = x^2$ on domain $[0.1,1]$ with initial condition $u_0(x) = \cos (\pi x /2)$. The exact solution is given by $u(x,t) = u_0(x/(1+tx))/(1+tx)^2$ and it is used to impose inflow boundary conditions (Section~\ref{sec:boundary}) on the left, and outflow boundary conditions are imposed on the right. While the flux $f(x,u)$ is linear in the $u$ variable, the problem is non-linear in the spatial variable in the sense that $I_h(a(x) u_h) \neq I_h(a) I_h(u_h)$ where $I_h$ is the interpolation operator. Thus, we use this problem to compare the \extrapolate~\eqref{eq:extrapolate} and \evaluate~\eqref{eq:evaluate} schemes for computation of trace values of the time averaged flux. The $L^2$ error versus total degrees of freedom plots are shown in Figure~\ref{fig:varadv1d} where a comparison with the standard Runge-Kutta method of the same order is made for \extrapolate, {\evaluate} fluxes. In Figure~\ref{fig:varadv1d}a, the cRKFR scheme with {\extrapolate} flux~\eqref{eq:extrapolate} is compared with the standard RKFR scheme. For degree $N=2$, the two schemes give similar performance. For the odd degrees $N=1,3$, it is seen that the cRKFR scheme with {\extrapolate} flux gives much larger errors than the RKFR scheme. However, optimal order of accuracy is shown by all schemes. In Figure~\ref{fig:varadv1d}b, the cRKFR scheme with {\evaluate} flux~\eqref{eq:evaluate} is compared with the standard RKFR scheme. In this case, the two schemes give similar performance for all degrees. In Figure~\ref{fig:varadv1d}c, the cRKFR schemes with {\extrapolate}~\eqref{eq:extrapolate} and {\evaluate}~\eqref{eq:evaluate} fluxes are compared. The {\evaluate} flux gives much better accuracy for degrees $N=1,3$ while for $N=2$, the two schemes give similar performance.

\subsubsection{Burgers' equation} \label{sec:burg}
\begin{figure}
\centering
\begin{tabular}{cc}
\includegraphics[width=0.45\textwidth]{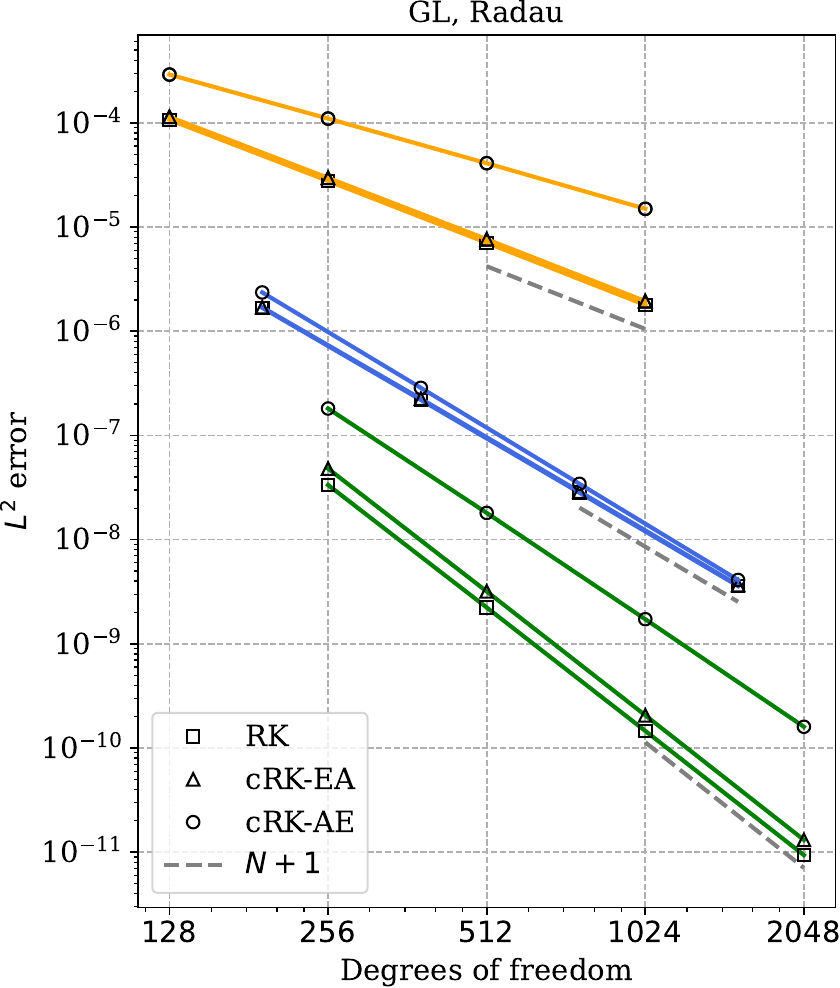} &
\includegraphics[width=0.45\textwidth]{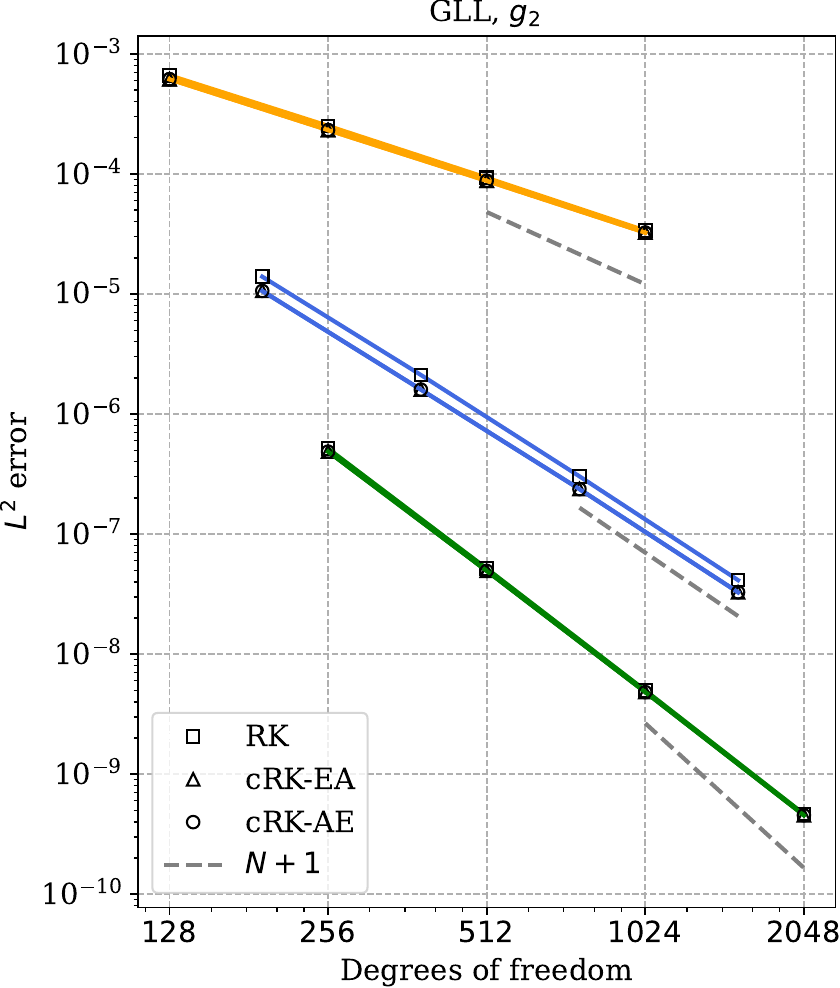}\\
(a) & (b)
\end{tabular}
\caption{Comparing {\extrapolate} and {\evaluate} schemes for 1-D Burgers' equation at $t = 2$ using Gauss-Legendre (GL) solution points with (a) Radau correction function with Gauss-Legendre (GL) solution points, (b) $g_2$ correction function with Gauss-Legendre-Lobatto (GLL) solution points.}
\label{fig:burg1d}
\end{figure}

We now test the simplest non-linear hyperbolic conservation law, the 1-D Burgers' equation with flux~\eqref{eq:con.law} given by $f(u) = u^2/2$. The initial condition is taken to be $u_0(x) = 0.2 \sin(x)$ on the domain $[0,2 \pi]$ with periodic boundary conditions. The solution develops a shock at time $t=5$ and thus we compute the error against the exact solution at time $t=2$ when the solution is smooth and optimal order of accuracy is expected. The exact solution is obtained by solving a nonlinear equation obtained through the method of characteristics. The $L^2$ error versus total degrees of freedom plots are shown in Figure~\ref{fig:burg1d} comparing the {\extrapolate} and {\evaluate} schemes for the cRKFR scheme, along with the standard RKFR scheme of the same order. In Figure~\ref{fig:burg1d}a, Radau correction functions with Gauss-Legendre solution points are used. The {\evaluate} scheme gives optimal order of accuracy for all degrees, and has comparable performance to the standard RKFR scheme. For odd degrees $N=1,3$, it is seen that the {\extrapolate} scheme gives noticeably larger errors than the standard RKFR scheme and suboptimal convergence of order close to $N+1/2$. In Figure~\ref{fig:burg1d}b, $g_2$ correction functions with Gauss-Legendre-Lobatto (GLL) solution points are used. In this case with GLL points, the {\extrapolate} and {\evaluate} schemes are equivalent as is verified through the error performance. The same performance is seen for cRKFR and RKFR for all degrees. Optimal order of accuracy is shown for degree $N=2$, but for the odd degrees $N=1,3$, the schemes give suboptimal convergence of order close to $N+1/2$. As discussed in Section~\ref{sec:numflux}, the D2 and D-CSX dissipation models give very similar performance. We provide a simple script which does the comparison of the two dissipation models for Burgers' equations (Section~\ref{sec:burg}) in the file~\href{https://github.com/Arpit-Babbar/paper-crkfr/blob/main/generate/dcsx.jl}{generate/dcsx.jl} at~\cite{crkrepo}.

\subsubsection{Composite signal}
\begin{figure}
\centering
\begin{tabular}{cccc}
\newcommand{\spacingcom}{\kern-2.2em}
\kern-2.2em \includegraphics[width=0.27\textwidth] {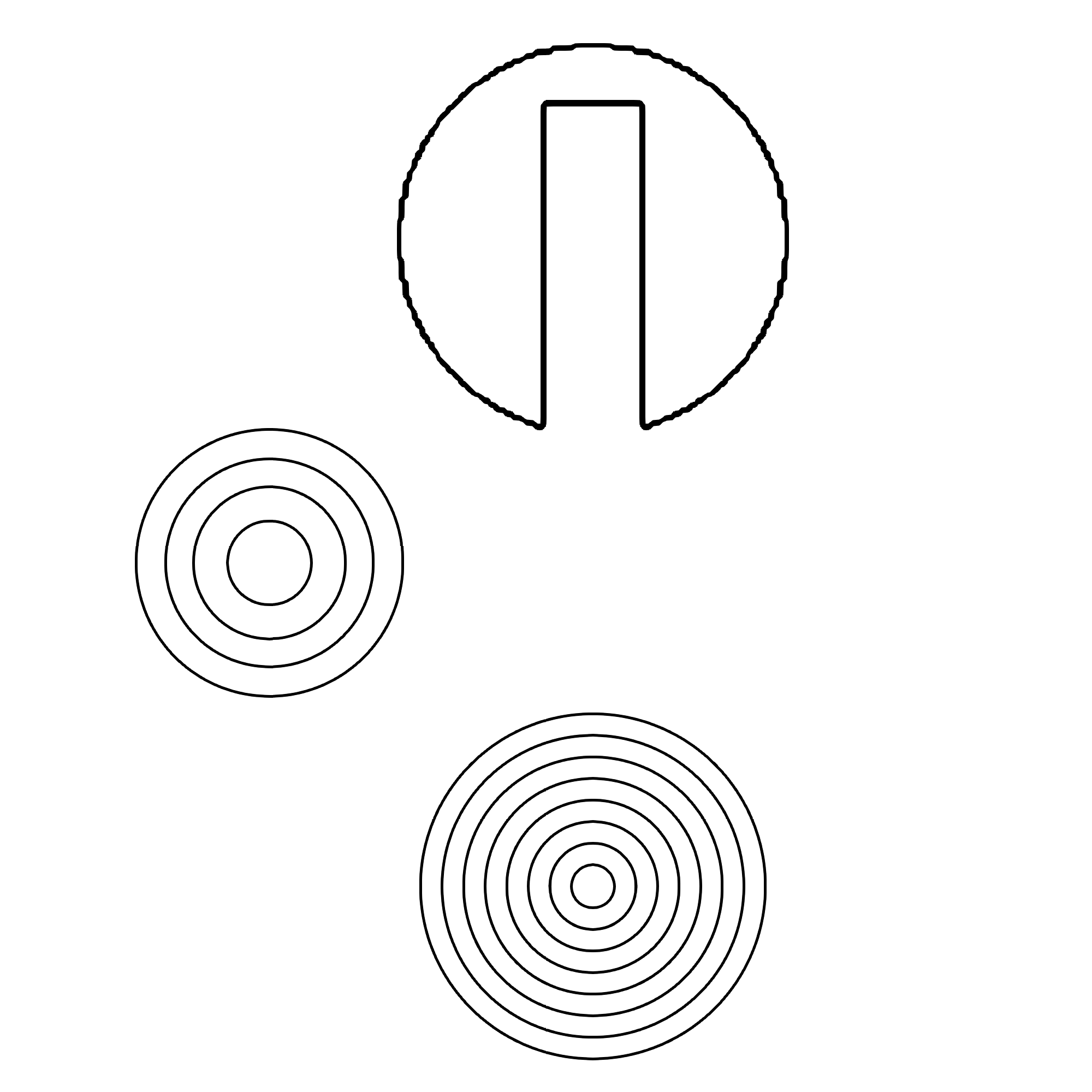} \kern-2.2em &
\includegraphics[width=0.27\textwidth]{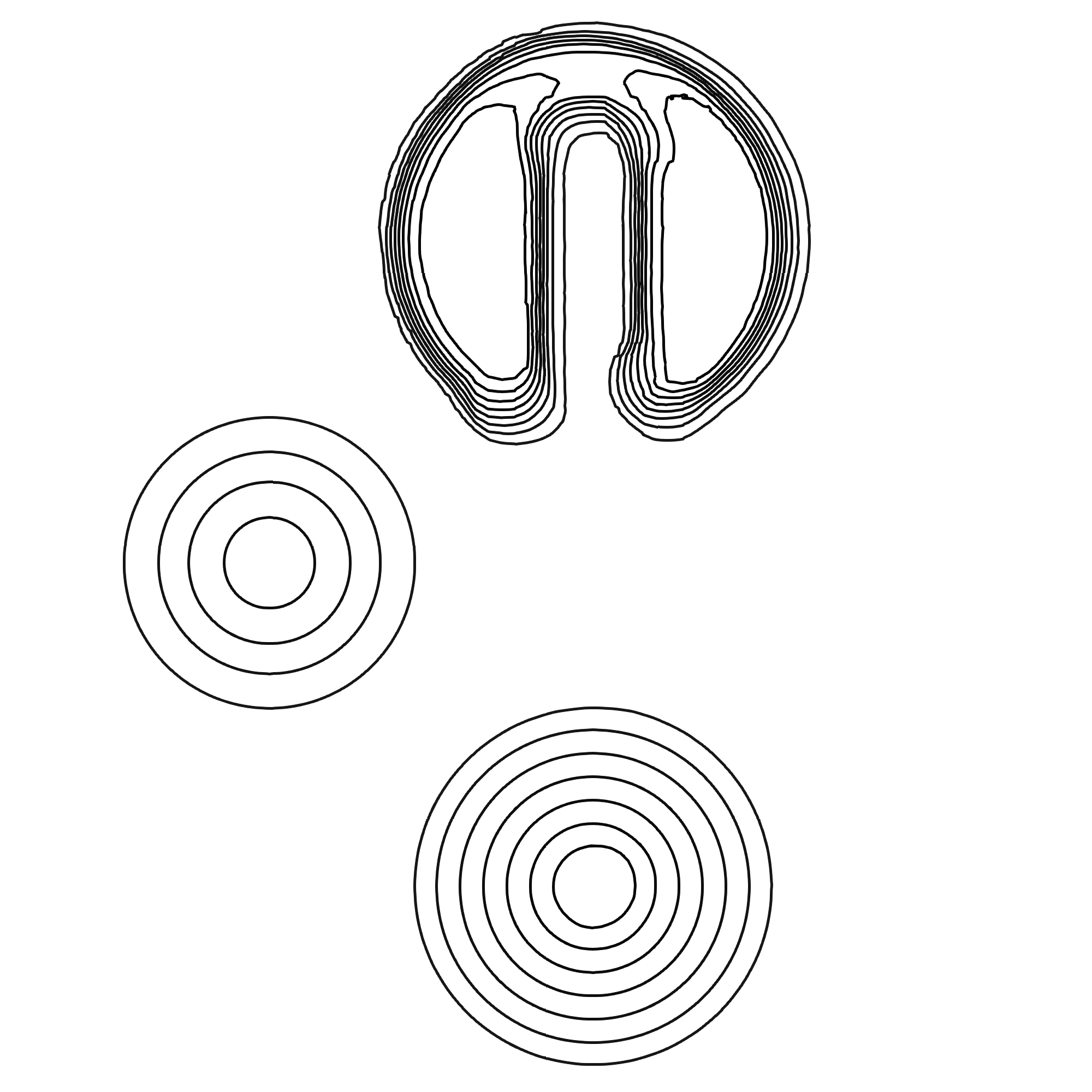} \kern-2.2em &
\includegraphics[width=0.27\textwidth]{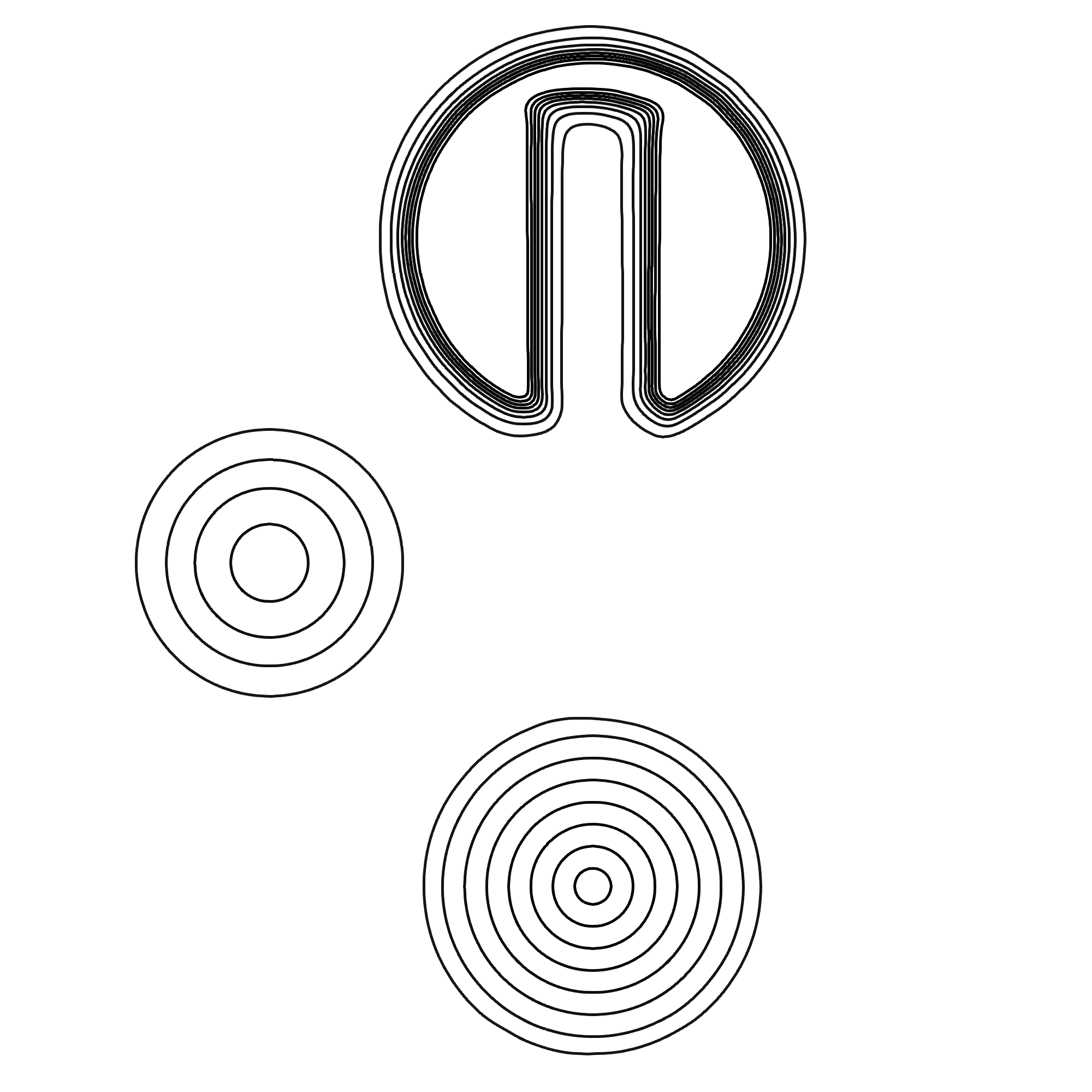} \kern-2.2em &
\includegraphics[width=0.27\textwidth]{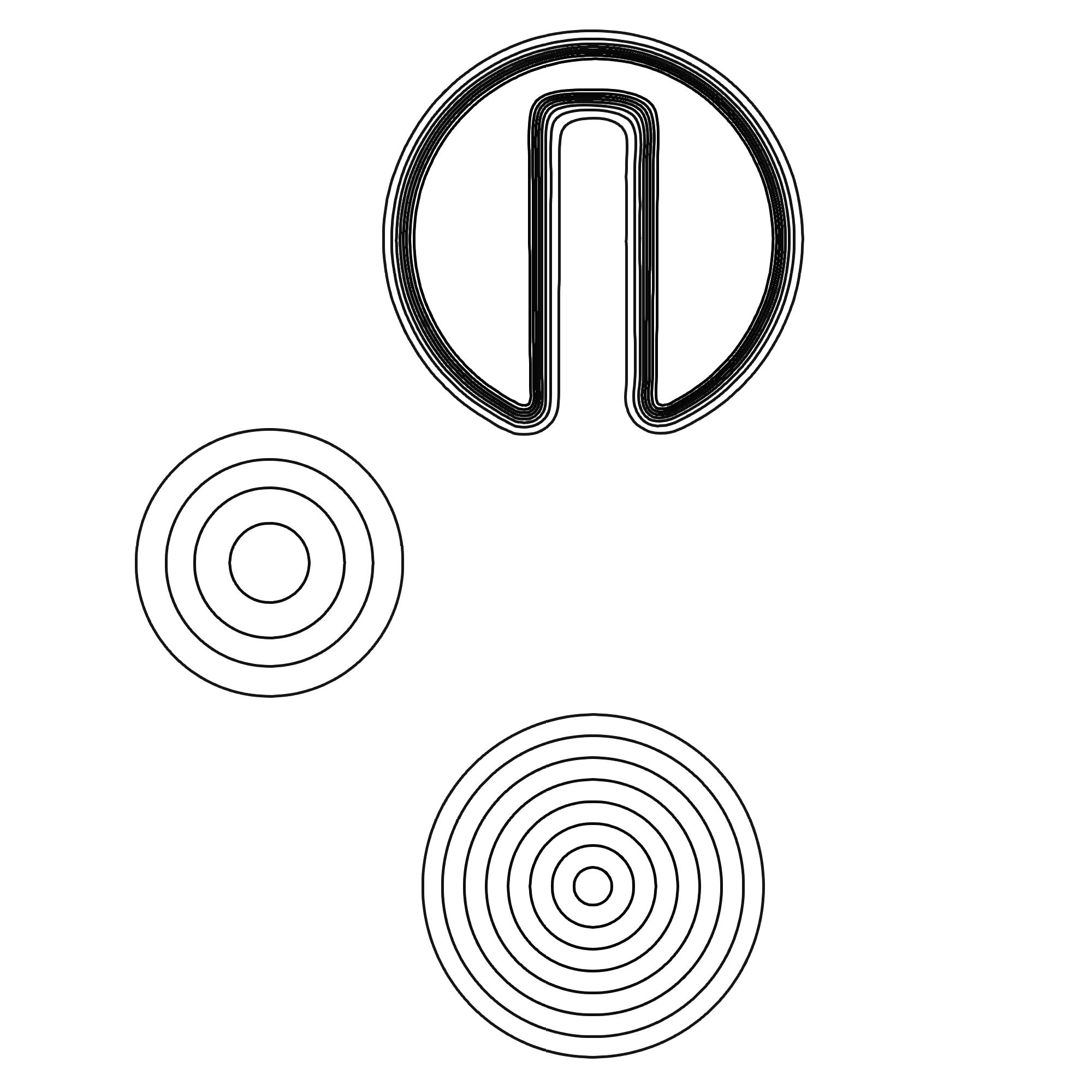} \kern-2.2em \\
\kern-2.2em (a) Exact \kern-2.2em & (b) TVB with $M=100$ \kern-2.2em & (c) FO blending \kern-2.2em & (d) MH blending \kern-2.2em
\end{tabular}
\caption{Rotation of a composite signal with velocity $\mathbf{a} = (\frac{1}{2} - y, x - \frac{1}{2})$,
numerical solution with polynomial degree $N=3$ on a mesh of $100^2$ elements.}
\label{fig:composite.signal.2d}
\end{figure}

The following is a test for the dissipation errors produced by the blending scheme (Section~\ref{sec:blending.scheme}). We compare the performance of the blending scheme with the TVB limiter of~\cite{Cockburn1989a} which was also used for the cRKDG scheme in~\cite{chen2024}. The following 2D advection equation with a divergence free velocity field is considered
\begin{equation}\label{eq:2dvaradv}
u_t+\nabla \cdot [a(x,y)u] =0, \qquad a=(1/2-y,x-1/2).
\end{equation}
The domain is taken to be $[0,1]^2$. The characteristics of the PDE~\eqref{eq:2dvaradv} are circles around the centre $(1/2,1/2)$ and thus the solution rotates around the centre. The initial condition, taken from~\cite{LeVeque1996}, consists of the following three profiles - a slotted disc, a cone and a smooth hump (see Figure~\ref{fig:composite.signal.2d}a). They are given by
\begin{align*}
u(x,y,0) &= u_1(x,y)+u_2(x,y)+u_3( x,y),\quad  (x,y)\in [0,1]\times[0,1]\\
u_1(x,y) &= \frac{1}{4}(1+\cos (\pi q( x,y))),\quad q(x,y) = \min ( \sqrt{(x-\bar x)^2+(y-\bar y)^2},r_0 )/r_0,{(\bar x,\bar y)} = (1/4,1/2), r_0=0.15\\
u_2(x,y) &= \begin{cases}
1-\dfrac{1}{r_0} \sqrt{(x-\bar x)^2+(y-\bar y)^2} & \mbox{ if } (x-\bar x)^2+(y-\bar y)^2\le r_0^2\\
0 & \mbox{otherwise}
\end{cases}, \quad {(\bar x,\bar y)} =(1/2,1/4), r_0=0.15\\
u_3(x,y) &= \begin{cases}
1 & \mbox{ if } (x,y) \in \mathrm{C}\\
0 & \mbox{otherwise}
\end{cases}
\end{align*}
where $\mathrm{C}$ is a slotted disc with center at $(0.5,0.75)$ and radius of $r_0 = 0.15.$ That is, $C$ is obtained by removing the set $\{(x,y): |x - 1/2| < r_0/4, y < 3/4 + 0.7 r_0\}$ from the disc.

We look at the numerical solution after it completes one time period at $t= 2 \pi$ in Figure~\ref{fig:composite.signal.2d}. We compare contour plots of polynomial solutions obtained using the cRKFR method of degree $N=3$ with the TVB limiter using a fine-tuned parameter $M=100$ (Figure~\ref{fig:composite.signal.2d}b), and with blending limiter using the first order (FO) (Figure~\ref{fig:composite.signal.2d}c) and MUSCL-Hancock  (MH) reconstructions (Figure~\ref{fig:composite.signal.2d}d). The smooth Gaussian profile is equally resolved by all three methods. For the cone profile, both the blending schemes show equally well-resolved profiles, which are slightly less dissipated than the one obtained by the TVB limiter. The most benefit of the blending limiter is seen in resolution of the sharp features of the slotted disc. The solution by the TVB limiter is a lot more smeared than both the blending schemes. The MUSCL-Hancock blending scheme is slightly better than the first order blending scheme in this case.

\subsection{1-D Euler's equations}
We now consider the 1-D Euler equations, which are a system of non-linear hyperbolic conservation laws given by
\begin{equation}
\label{eq:1deuler} \pd{}{t}  \left(\begin{array}{c}
\rho\\
\rho v_1\\
E
\end{array}\right) + \pd{}{x}  \left(\begin{array}{c}
\rho v_1\\
p + \rho v_1^2\\
(E + p) v_1
\end{array}\right) = \bzero
\end{equation}
where $\rho, v_1, p$ and $E$ denote the density, velocity, pressure and total
energy per unit volume of the gas, respectively. For a polytropic gas, an equation of state $E
= E (\rho, v_1, p)$ which leads to a closed system is given by
\begin{equation}
\label{eq:state} E = E (\rho, v_1, p) = \frac{p}{\gamma - 1} + \frac{1}{2}
\rho v_1^2
\end{equation}
where $\gamma > 1$ is the adiabatic constant, that will be taken as $1.4$
which is the value for air. The time step size is computed using~\eqref{eq:time.step.2d} with the spectral radii for 1-D Euler's equations~\eqref{eq:1deuler} approximated as
\begin{equation} \label{eq:spectral.radius.1d}
\sigma(\pf'(\uu)) = |v_1| + c, \qquad \sigma(\pg'(\uu)) = 0, \qquad c = \sqrt{\gamma p/\rho}.
\end{equation}
In the following results, wherever it is not mentioned, we use the Rusanov flux~\eqref{eq:rusanov.wave.speed}. We use Radau correction function with Gauss-Legendre solution points because of their superiority seen in the scalar equations, as is consistent with the LWFR results in~\cite{babbar2022,babbar2024admissibility}. We use {\evaluate} flux~\eqref{eq:evaluate} because of its accuracy benefit and D2 dissipation~\eqref{eq:D2} because of its higher CFL numbers (Table 1 of~\cite{babbar2022}). The tests consist of problems with strong shocks with high pressure ratios in order to demonstrate the scheme's capability to preserve admissibility of the solution. We also show the performance comparison of the TVB limiter and the blending scheme to show that the blending scheme is able to resolve extremum and contact discontinuities well.

\subsubsection{Blast wave}
\begin{figure}
\centering
\begin{tabular}{cc}
\includegraphics[width=0.45\textwidth]{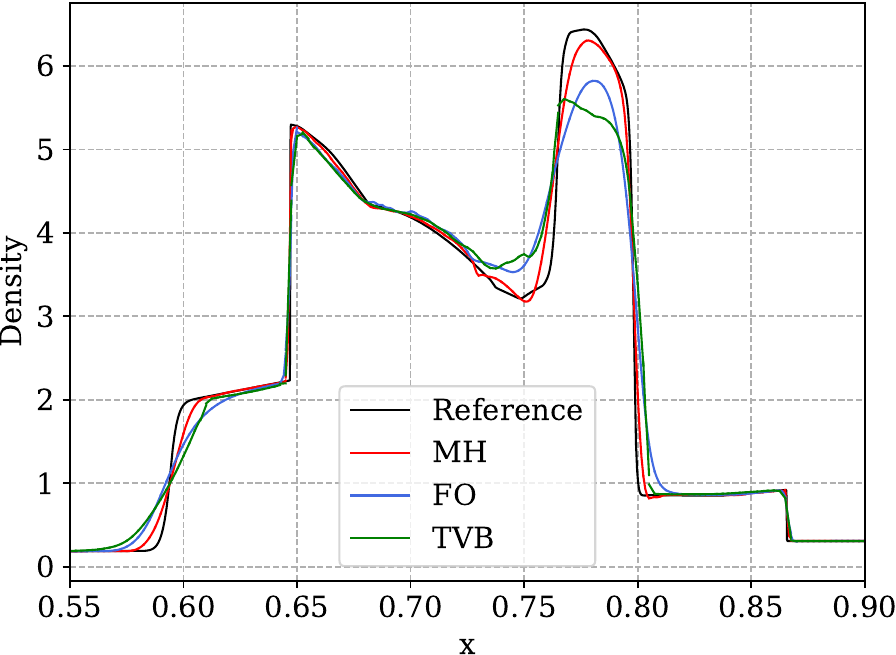} &
\includegraphics[width=0.45\textwidth]{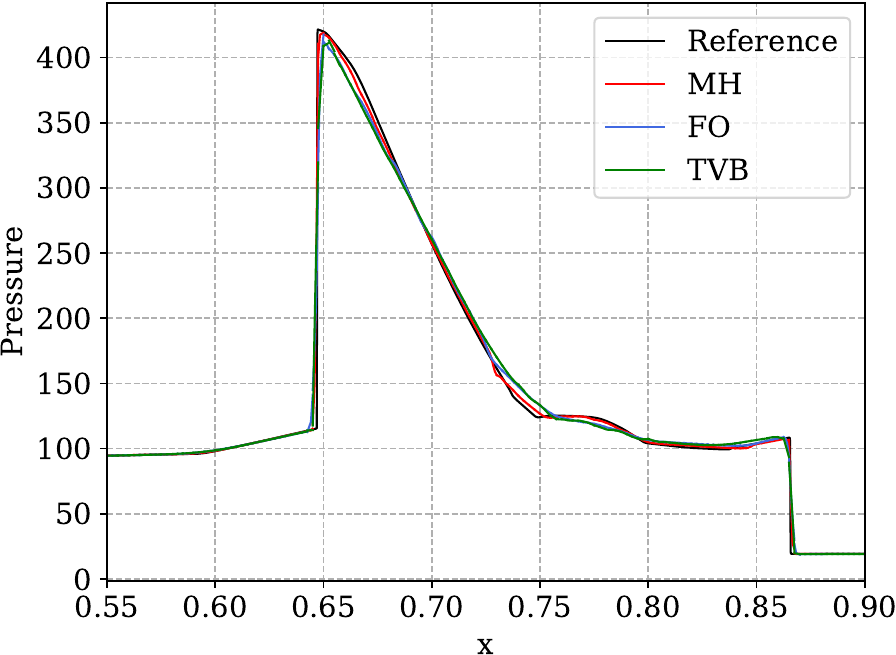} \\
(a) & (b)
\end{tabular}
\caption{Blast wave problem, numerical solution with degree $N=3$ using first order (FO) and MUSCL-Hancock (MH) blending schemes, and TVB limited scheme (TVB-300) with parameter $M=300$. (a) Density, (b) pressure profiles are shown at time $t=0.038$ on a mesh of 400 elements.}
\label{fig:blast}
\end{figure}

This blast wave test from~\cite{Woodward1984} consists of an initial condition with two discontinuities given as
\begin{equation*}
(\rho,v_1,p)=\begin{cases}
(1,0, 1000), & \mbox{ if } x<0.1,\\
(1,0,0.01), & \mbox{ if }  0.1 < x <0.9,\\
(1,0, 100), & \mbox{ if } x> 0.9,
\end{cases}
\end{equation*}
in the domain $[0,1]$. The boundaries are set as solid walls by imposing the reflecting boundary conditions at $x=0$ and $x=1$~\eqref{eq:reflect.bc}. The two discontinuities lead to two Riemann problems, each of which lead to a shock, rarefaction and a contact discontinuity. The boundary conditions cause reflection of shocks and expansion waves off the solid wall and several wave interactions inside the domain. As time goes on, the two reflected Riemann problems interact, which is the main point of interest for this test~\cite{Woodward1984}. Without positivity correction, negative pressures are obtained by the cRKFR scheme during the shock interaction. The interaction of the two Riemann problems also cause the formation of an extremum in the density profile. Dissipation errors of limiters are often seen in the presence of extrema. Thus, we use this test to demonstrate the performance of the First Order (FO) and MUSCL-Hancock (MH) blending schemes by comparing them with the TVB limiter. The numerical results are shown in Figure~\ref{fig:blast} with a grid of 400 elements using polynomial degree $N=3$ at time $t=0.038$. The TVB limiter is used with the fine tuned parameter $M=300$. Figure~\ref{fig:blast} shows the density and pressure profiles of the numerical solution, compared with a reference solution computed using a very fine mesh. The superior of blending scheme is evident, especially when compared with the TVB limiter, by looking at peak amplitude (extremum) and contact discontinuity of the density profile (Figure~\ref{fig:blast}a).

\subsubsection{Titarev Toro}
\begin{figure}
\centering
\begin{tabular}{cc}
\includegraphics[width=0.45\textwidth]{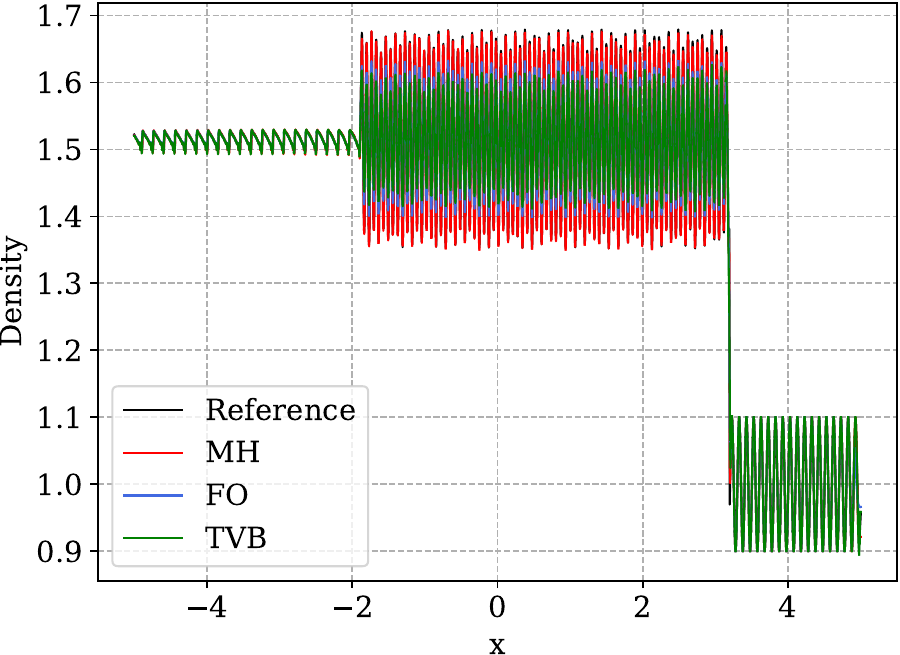} &
\includegraphics[width=0.45\textwidth]{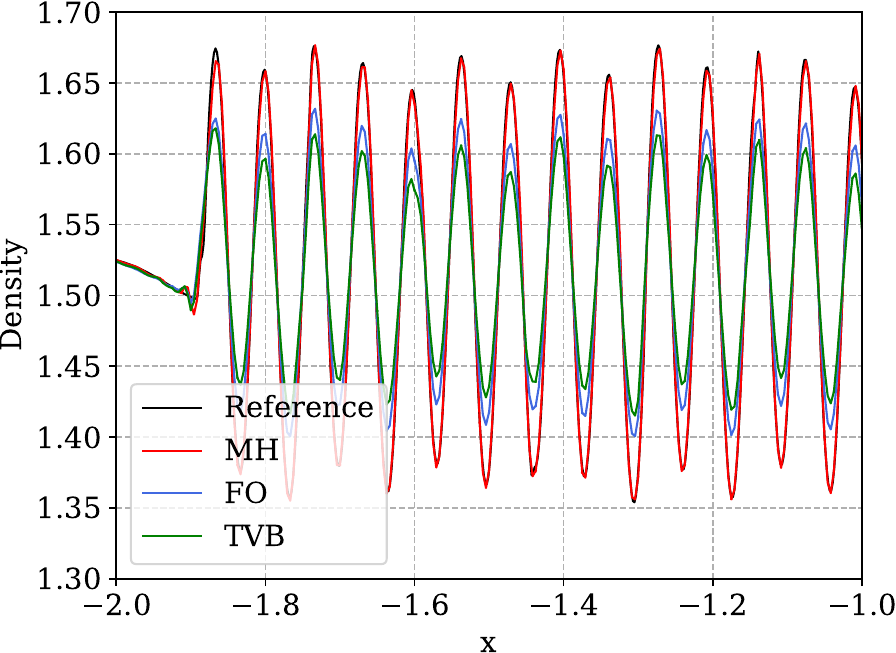} \\
(a) & (b)
\end{tabular}
\caption{Titarev-Toro problem, comparing First Order (FO) and MUSCL-Hancock
(MH) blending (a) Complete, (b) Zoomed density profile near smooth extrema on a
mesh of 800 elements.}
\label{fig:titarev.toro}
\end{figure}

This test case by Titarev and Toro~{\cite{Titarev2004}} consists of severely
oscillatory waves (extrema) interacting with a shock. The initial condition is given by
\[ (\rho, v_1, p) = \begin{cases}
(1.515695, 0.523346, 1.805), \qquad & - 5 \leq x \leq - 4.5\\
(1 + 0.1 \sin (20 \pi x), 0, 1),  & - 4.5 < x \leq 5
\end{cases} \]
The physical domain is $[- 5, 5]$ with outflow boundary conditions (Section~\ref{sec:boundary}). The challenge in this test case is to have a robust scheme near the discontinuity without adding too much dissipation to the smooth extrema as they pass through the shock. Thus, we compare the performance of the First Order (FO) and MUSCL-Hancock (MH) blending schemes with the TVB limiter by looking at how they resolve the smooth extrema. The
density profile at $t = 5$ is shown in Figure~\ref{fig:titarev.toro}. The TVB limiter is using the fine tuned parameter $M=300$. The zoomed Figure~\ref{fig:titarev.toro}b clearly shows that the blending schemes give better resolution of the smooth extrema than the TVB limiter. The MUSCL-Hancock has nearly resolved the extrema, and is noticeably better than the first order scheme.

\subsubsection{Sedov's blast}\label{sec:sedov.blast.1d}
\begin{figure}
\centering
\begin{tabular}{cc}
\includegraphics[width=0.45\textwidth]{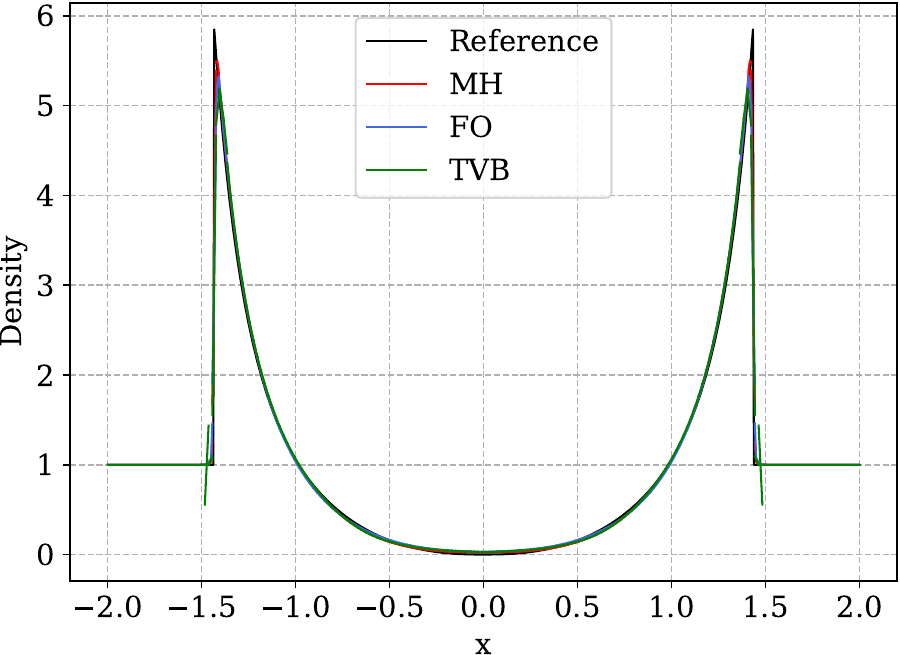} &
\includegraphics[width=0.45\textwidth]{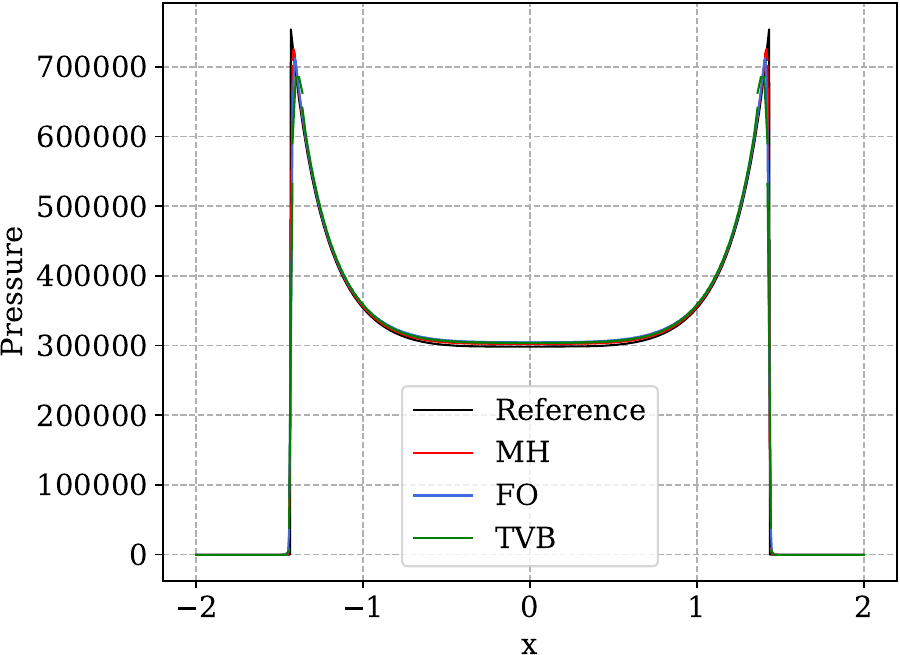} \\
(a) & (b)
\end{tabular}
\caption{Sedov's blast wave problem, numerical solution with degree $N=3$ using first order (FO) and MUSCL-Hancock blending schemes. (a) Density and (b) pressure profiles of numerical solutions are plotted at time $t=0.001$ on a mesh of $201$ elements.}
\label{fig:sedov.blast}
\end{figure}

Sedov's blast wave test~\cite{sedov1959} describes a point injection in a gas, e.g., from an explosion. The injection generates a radial shock wave that compresses the gas and causes extreme temperatures and pressures. The problem can be formulated in one dimension as a special case, where the explosion occurs at $x=0$ and the gas is confined to the interval $[-1,1]$ by solid walls~\eqref{eq:reflect.bc}. For the simulation, on a grid of 201 elements with solid wall boundary conditions, we use the following initial data from~\cite{Zhang2012},
\[
\rho = 1, \qquad
v_1 = 0, \qquad
E(x) = \begin{cases}
\frac{3.2 \times 10^6}{\Delta x}, \qquad & |x| \le \frac{\Delta x}{2}, \\
10^{-12}, \qquad & \text{otherwise},
\end{cases}
\]
where $\Delta x$ is the element width. The extremely high pressure ratios at the shock make it challenging to preserve admissibility of the solution. This test is known to only work for schemes which are provably admissibility preserving. Thus, nonphysical solutions are obtained for the cRKFR scheme if the proposed admissibility preservation corrections from Section~\ref{sec:flux.limiter} are not applied. The density and pressure profiles at $t=0.001$ obtained using blending schemes are shown in Figure~\ref{fig:sedov.blast}. The TVB limiter is used by generalizing the positivity preserving framework in Section~\ref{sec:flux.limiter}, following~\cite{babbar2024generalized}. The parameter $M=0$ was found to be the best for this test. Compared to the TVB limiter, the blending scheme gives a sharper resolution of the shock while also giving better control of spurious oscillations.

\subsection{2-D Euler's equations} \label{sec:2d.euler}
We consider the two-dimensional Euler equations of gas dynamics given by
\begin{equation}
\label{eq:2deuler} \pd{}{t}  \left(\begin{array}{c}
\rho\\
\rho v_1\\
\rho v_2\\
E
\end{array}\right) + \pd{}{x}  \left( \begin{array}{c}
\rho v_1\\
p + \rho v_1^2\\
\rho v_1 v_2\\
(E + p) v_1
\end{array} \right) + \pd{}{y}  \left( \begin{array}{c}
\rho v_2\\
\rho v_1 v_2\\
p + \rho v_2^2\\
(E + p) v_2
\end{array} \right) = \bzero,
\end{equation}
where $\rho, p$ and $E$ denote the density, pressure and total energy of the
gas, respectively and $(v_1, v_2)$ are Cartesian components of the fluid velocity.
For a polytropic gas, an equation of state $E = E (\rho, v_1, v_2, p)$ which leads
to a closed system is given by
\begin{equation}
\label{eq:2dstate} E = E (\rho, v_1, v_2, p) = \frac{p}{\gamma - 1} +
\frac{1}{2} \rho (v_1^2 + v_2^2),
\end{equation}
where $\gamma > 1$ is the adiabatic constant, that will be taken as $1.4$ unless specified otherwise, which is the typical value for air. The time step size is computed using~\eqref{eq:time.step.2d} with the spectral radii for Euler's equations~\eqref{eq:2deuler} approximated as
\[
\sigma(\pf'(\uu)) = |v_1| + c, \qquad \sigma(\pg'(\uu)) = |v_2| + c, \qquad c = \sqrt{\gamma p/\rho}.
\]
Unless specified otherwise, we use the Rusanov numerical flux~\eqref{eq:rusanov.wave.speed} at the element interfaces.

Since we have already established the superiority of the blending scheme over the TVB limiter, for further verification of the cRKFR scheme with blending, we will compare our results with {\tt Trixi.jl}~\cite{Ranocha2022}. {\tt Trixi.jl} is a Julia package for solving hyperbolic conservation laws using the Runge-Kutta Discontinuous Galerkin methods and contains the blending limiter of~\cite{hennemann2021} with a first order finite volume scheme on the subcells using Gauss-Legendre-Lobatto (GLL) solutions points. As discussed in Section~\ref{sec:blending.scheme}, in order to minimize dissipation errors, our blending scheme instead uses the MUSCL-Hancock reconstruction on subcells, and Gauss-Legendre (GL) solution points. Thus, by making the comparison with {\tt Trixi.jl}, we numerically verify the accuracy gain of the proposed cRKFR blending scheme. In all numerical experiments, both solvers use the same time step sizes in ~\eqref{eq:time.step.2d}. We have also performed experiments using cRKFR scheme with first order blending scheme and Gauss-Legendre (GL) points, and observed lower accuracy than the MUSCL-Hancock blending scheme, but higher accuracy than the first order blending scheme implementation of {\tt Trixi.jl} using GLL points. However, to save space, we have not presented the results of cRKFR scheme with first order blending.
\subsubsection{Isentropic vortex problem}
\begin{figure}
\centering
\includegraphics[width=0.6\textwidth]{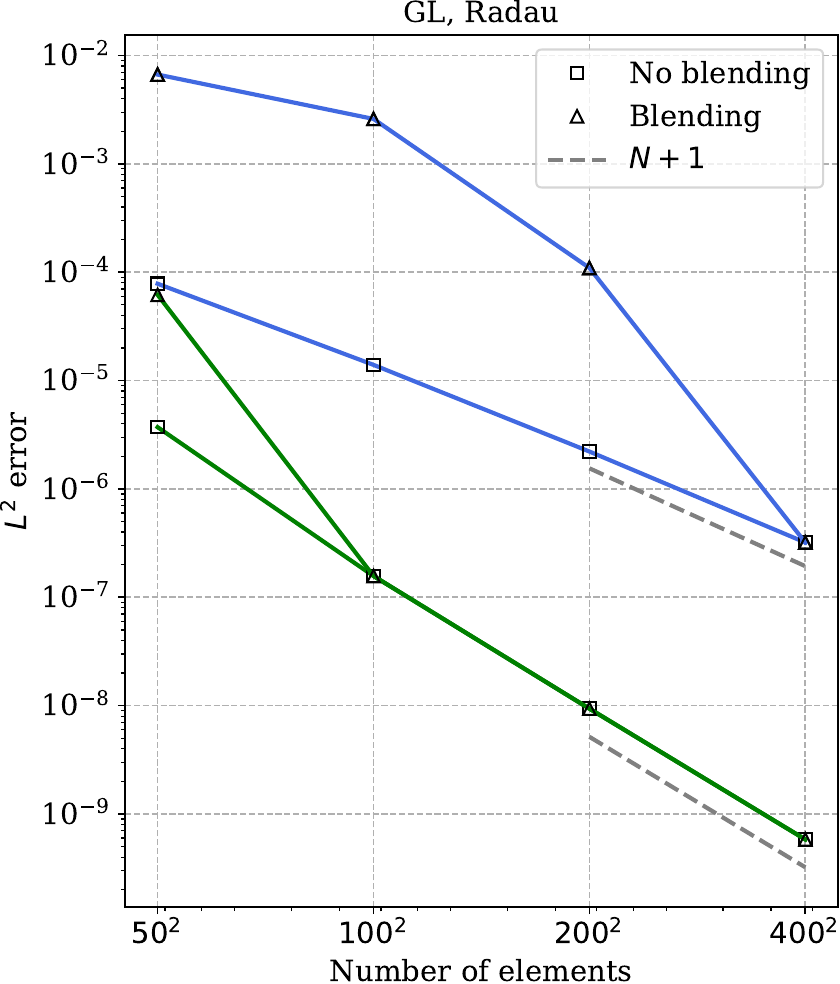}
\caption{$L^2$ error versus number of elements for polynomial degrees $N=2,3$ for the isentropic vortex test using Radau correction function with Gauss-Legendre (GL) solution points and Radau correction function.}
\label{fig:isentropic}
\end{figure}

We first test the cRKFR scheme for a setup where the exact solution is smooth and is analytically known. The setup, taken from~\cite{Yee1999,Spiegel2016}, consists of an isentropic vortex that advects at a constant velocity. The initial conditions with periodic boundary conditions are taken to be
\[
\rho = \left[ 1 - \frac{\beta^2 (\gamma - 1)}{8 \gamma \pi^2} \exp (1 -
r^2) \right]^{\frac{1}{\gamma - 1}}, \qquad v_1 = M_\infty \cos \alpha -
\frac{\beta (y - y_c)}{2 \pi} \exp \left( \frac{1 - r^2}{2} \right)
\]
\[
v_2 = M_\infty \sin \alpha + \frac{\beta (x - x_c)}{2 \pi} \exp \left( \frac{1 -
r^2}{2} \right), \qquad r^2 = (x - x_c)^2 + (y - y_c)^2,
\]
and the pressure is set as $p = \rho^{\gamma}$. The parameters of the vortex setup are chosen to be $\beta = 5$, $M_\infty = 0.5$, $\alpha = 45^o$, $(x_c, y_c) = (0, 0)$. The physical domain is taken to be $[- 10, 10] \times [- 10, 10]$.  For the initial condition $\uu_0 = \uu_0(x,y)$, the exact solution at time $t$ is given by $\uu_0(x-u_0t, y-v_0t)$ where $u_0 = M_\infty \cos \alpha, v_0 = M_\infty \sin \alpha$.
The simulation is run till $t=1$ units\footnote{Optimal order of accuracy and same WCT conclusions as in Table~\ref{tab:wall.clock} are also seen for the full time period $T = 20 \sqrt{2} M_\infty$. Results are shown for a shorter time for quick reproducibility of benchmarking file~\href{https://github.com/Arpit-Babbar/paper-crkfr/blob/main/generate/isentropic.jl}{generate/isentropic.jl} at~\cite{crkrepo}} and Figure~\ref{fig:isentropic} shows the $L^2$ error of the density variable sampled at $N+3$ equispaced points against grid resolution with and without the blending limiter. It can be seen that the limiter does not activate for adequately high resolution, yielding the same optimal convergence rates as those achieved without the limiter, as shown in Figure~\ref{fig:isentropic}. A performance comparison between Lax-Wendroff (LW) and standard multistage Runge-Kutta (RK) schemes has been made in previous works~\cite{Qiu2003,Qiu2005b,babbar2022} revealing the benefit of LW discretization. Thus, we only make a performance comparison of the cRKFR scheme with the LWFR scheme in Table~\ref{tab:wall.clock} implemented in {\tt Tenkai.jl} used for~\cite{babbar2022,babbar2024admissibility,babbar2024generalized}. It is seen that the cRKFR scheme is faster than the LWFR scheme up to a factor of $1.38$. The $L^2$ errors of the LWFR scheme are very close to cRKFR scheme and are thus not shown.

\begin{table}
\centering
\begin{tabular}{|c|c|c|c|c|c|c|c|c|c|}
\hline
$\nel$ & \multicolumn{3}{|c|}{$N=1$} & \multicolumn{3}{|c|}{$N=2$} & \multicolumn{3}{|c|}{$N=3$} \\
\hline
\textbf{ } & \textbf{LW} & \textbf{cRK} & \textbf{Ratio} & \textbf{LW} & \textbf{cRK} & \textbf{Ratio} & \textbf{LW} & \textbf{cRK} & \textbf{Ratio} \\\hline
$50^2$ & 0.14 & 0.12 & 1.14 & 0.5 & 0.52 & 0.96 & 1.7 & 1.2 & 1.37 \\
$100^2$ & 1.1 & 0.9 & 1.16 & 3.9 & 3.0 & 1.31 & 14.0 & 10.2 & 1.36 \\
$200^2$ & 9.2 & 7.9 & 1.15 & 32.1 & 25.9 & 1.23 & 115.2 & 84.3 & 1.36 \\
$400^2$ & 78.2 & 67.9 & 1.15 & 262.5 & 202.3 & 1.29 & 937.4 & 677.7 & 1.38 \\\hline
\end{tabular}
\caption{Wall clock time performance (seconds) comparison of LWFR and cRKFR schemes for different polynomial degrees $N$ and mesh with number of elements $\nel$.}\label{tab:wall.clock}
\end{table}

\subsubsection{2-D Riemann problem}
\begin{figure}
\centering
\begin{tabular}{cc}
\includegraphics[width=0.45\textwidth]{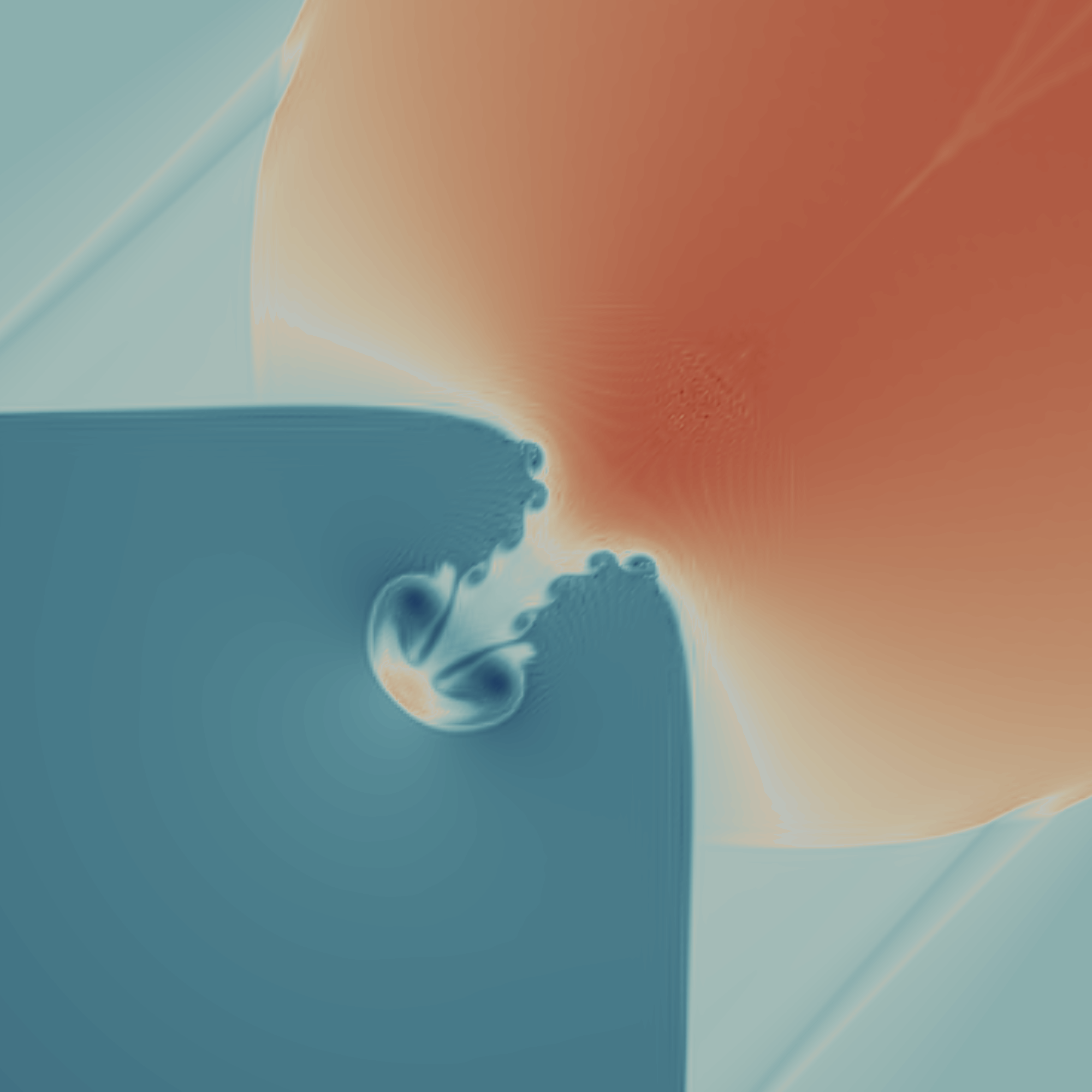} &
\includegraphics[width=0.45\textwidth]{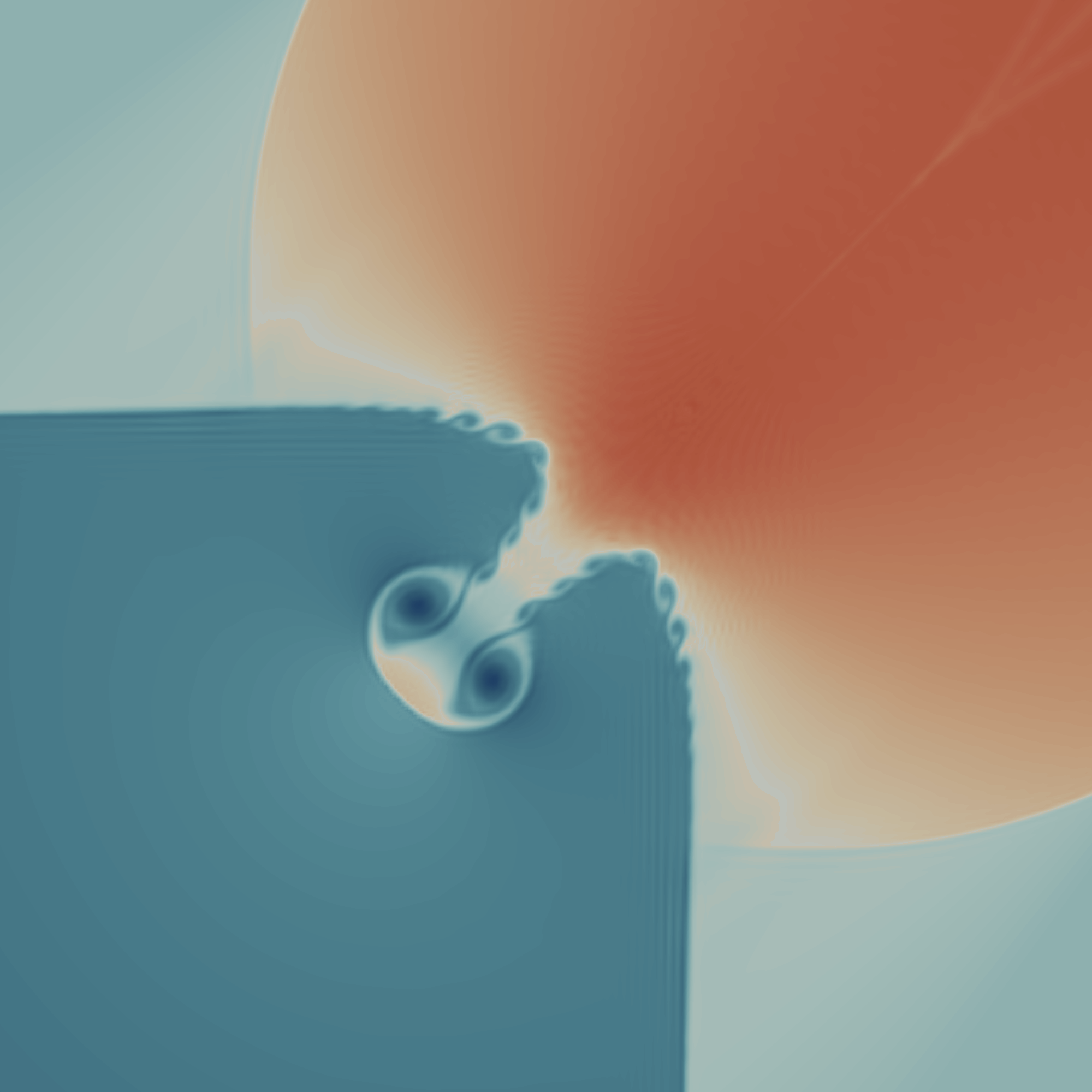} \\
(a) {\tt Trixi.jl} & (b) cRK-MH
\end{tabular}
\caption{2-D Riemann problem, density plots of numerical solution at $t=0.25$ for degree $N=3$ on a $256 \times 256$ mesh.}
\label{fig:rp2d}
\end{figure}

In this test, we show the accuracy of cRKFR scheme with the blending limiter by showing its ability to capture small scale shear structures that appear in a 2-D Riemann problem. The initial condition of a 2-D Riemann problem consists of four constant states. A theoretical and numerical study for the 2-D Riemann problem particular to gas dynamics can be found in~\cite{Glimm1985}. The physical domain is taken to be $[0,1]^2$, and the four initial states are chosen so that the jump across any two has an elementary planar wave solution, i.e., a shock, rarefaction or contact discontinuity. There are only 19 such genuinely different configurations possible~\cite{Zhang1990,Lax1998}. We consider configuration 12 of~\cite{Lax1998} consisting of 2 positive slip lines (contact discontinuities) and two forward shocks. The initial condition is given by
\begin{align*}
(\rho, v_1, v_2, p) = \begin{cases}
(0.5313, 0, 0, 0.4)  \qquad & \text{if } x \ge 0.5, y \ge 0.5, \\
(1, 0.7276, 0, 1)  & \text{if } x < 0.5, y \ge 0.5, \\
(0.8, 0, 0, 1) & \text{if } x < 0.5, y < 0.5, \\
(1, 0, 0.7276, 1) & \text{if } x \ge 0.5, y < 0.5.
\end{cases}
\end{align*}
The simulations are performed with outflow boundary conditions (Section~\ref{sec:boundary}) on an enlarged domain up to time $t=0.25$. The density profiles obtained from the MUSCL-Hancock blending scheme and {\tt Trixi.jl} are shown in Figure~\ref{fig:rp2d}. We see that both schemes give similar resolution in most regions. The MUSCL-Hancock blending scheme gives better resolution of the small scale structures arising across the slip lines, demonstrating the accuracy of the blending scheme.

\subsubsection{Double Mach reflection}
Similar to the 2-D Riemann problem, the double Mach reflection problem demonstrates the capability of the blending scheme to capture small scale shear structures that occur at a contact discontinuity. Additionally, it consists of a Mach 10 shock, making it a test for admissibility preservation as well. It was originally proposed by Woodward and Colella~\cite{Woodward1984} and consists of a shock impinging on a wedge/ramp which is inclined by 30 degrees. The solution consists of a self similar shock structure with two triple points. By a change of coordinates, the situation is simulated in the rectangular domain $\Omega = [0,4] \times [0,1],$ where the wedge/ramp is positioned at $x=1/6, y=0.$ Defining $\mathbf{u}_b = \mathbf{u}_b(x,y,t)$ with primitive variables given by
\begin{equation*}
(\rho,v_1,v_2,p)=\begin{cases}
(8, 8.25 \cos\left( \frac{\pi}{6} \right), -8.25 \sin\left( \frac{\pi}{6} \right), 116.5), & \mbox{ if } x < \frac{1}{6} + \frac{y + 20 t}{\sqrt{3}}\\
(1.4, 0, 0, 1), & \mbox{ if } x > \frac{1}{6} + \frac{y + 20 t}{\sqrt{3}}
\end{cases}
\end{equation*}
we define the initial condition to be $\mathbf{u}_0(x,y) = \mathbf{u}_b(x,y,0)$. With $\mathbf{u}_b$, we impose inflow boundary conditions at the left side $\{0\} \times [0,1]$, outflow boundary conditions both at $[0,1/6] \times \{0\}$ and $\{4\} \times [0,1]$, reflecting boundary conditions at $[1/6, 4] \times \{0\}$ and inflow boundary conditions at the upper side $[0,4] \times \{1\}$.

\begin{figure}
\centering
\includegraphics[width=0.6\textwidth]{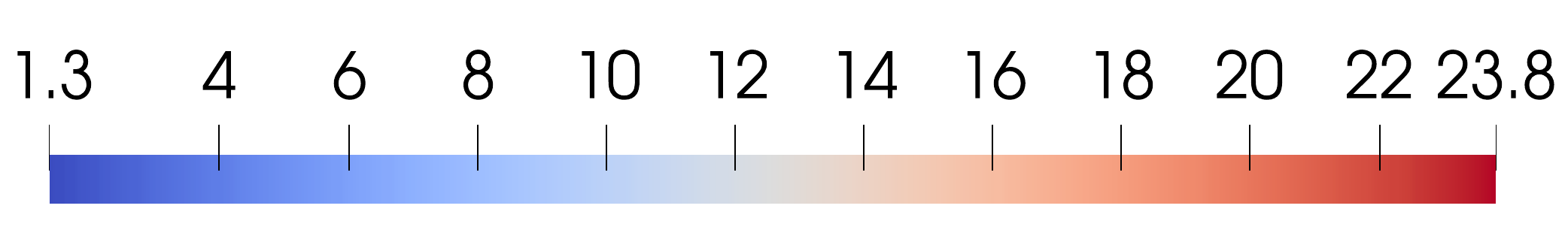}\\
{\noindent}\begin{tabular}{ccc}
\includegraphics[width=0.29\textwidth]{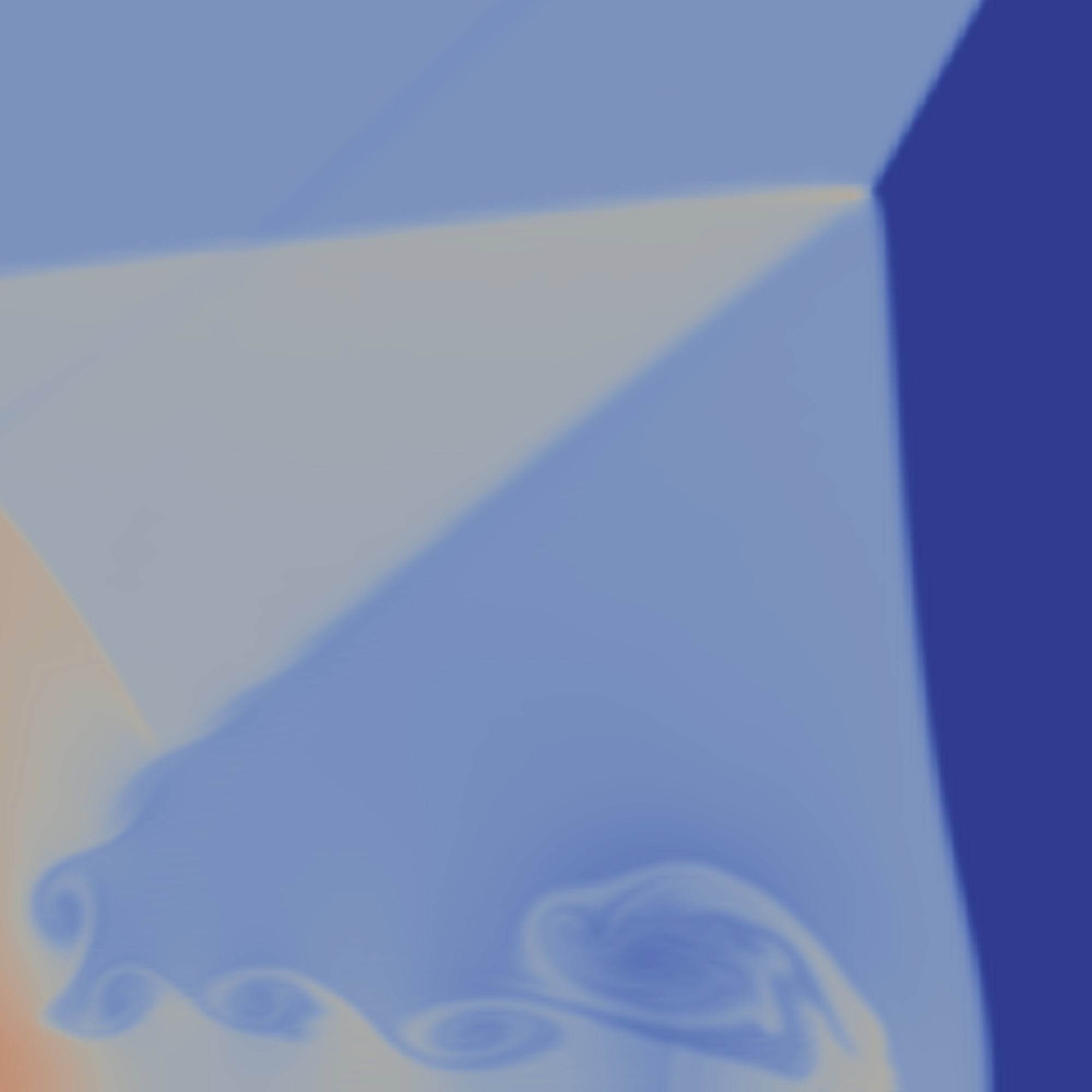}
&
\includegraphics[width=0.29\textwidth]{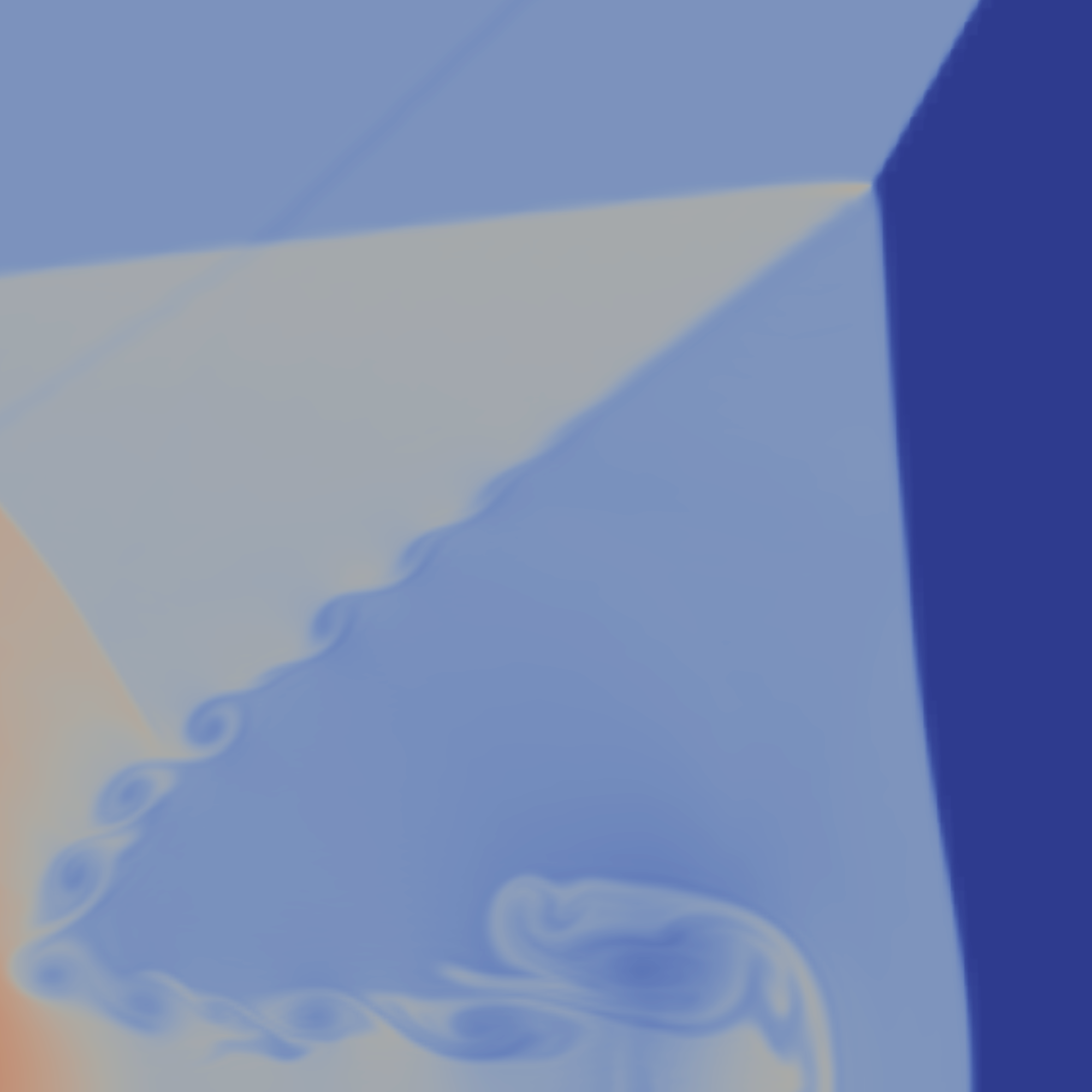}
&
\includegraphics[width=0.34\textwidth]{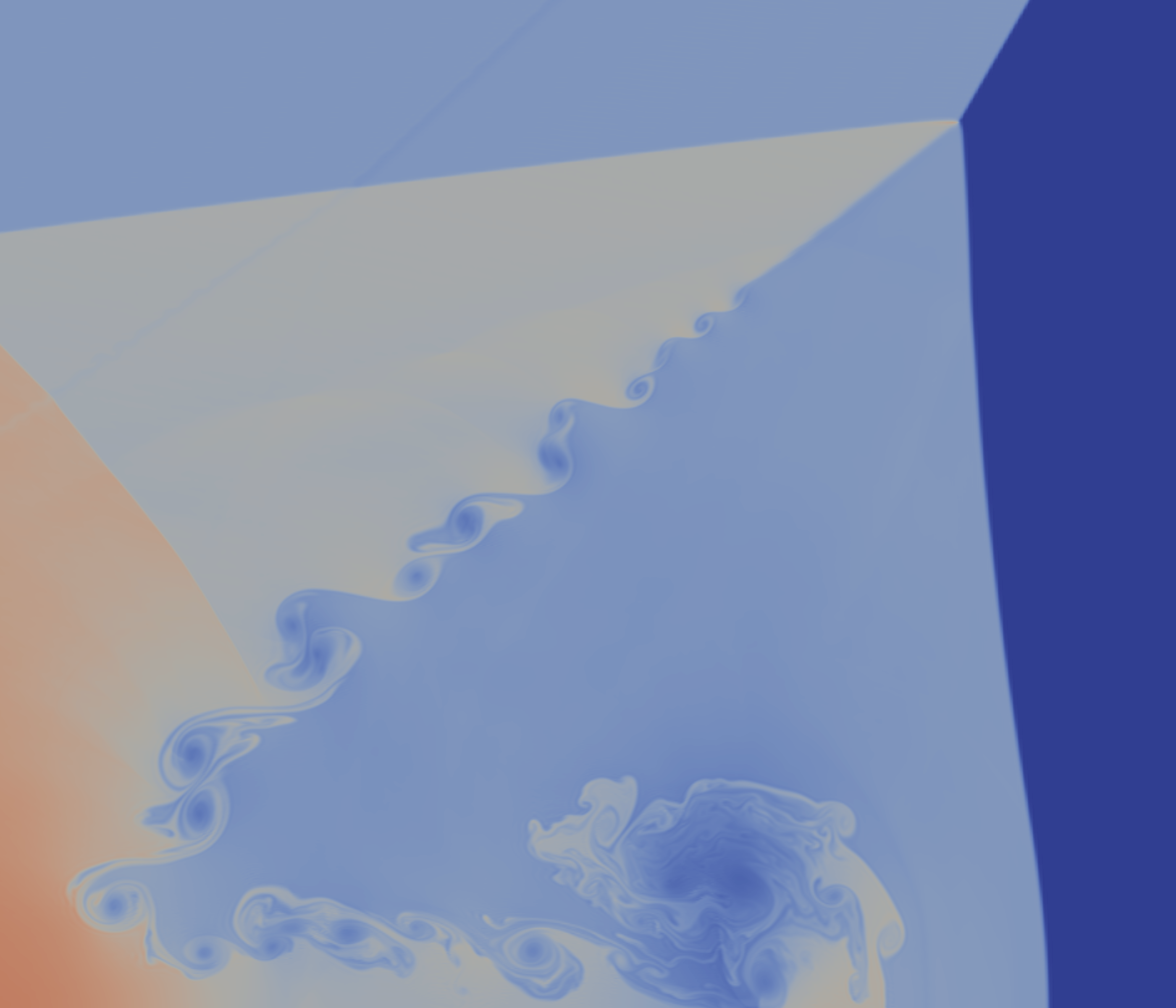}
\\
(a) {\tt Trixi.jl} & (b) cRK-MH & (c) Reference
\end{tabular}
\caption{Double Mach reflection problem, density plot of numerical solution at $t = 0.2$ on a $600 \times 150$ mesh zoomed near the primary triple point.\label{fig:dmr}}
\end{figure}

The simulation is run on a mesh of $600 \times 150$ elements using degree $N=3$ polynomials upto time $t=0.2$. In Figure~\ref{fig:dmr}, we compare the results of {\tt Trixi.jl} (Figure~\ref{fig:dmr}a) with the MUSCL-Hancock blended scheme (Figure~\ref{fig:dmr}b) zoomed near the primary triple point. In Figure~\ref{fig:dmr}c, we have a reference solution that is computed by a {\tt Trixi.jl} on a much finer mesh to demonstrate that the number of small scale structures increase with higher resolution. As expected, the small scale structures are captured better by the MUSCL-Hancock blended scheme making it look closer to the reference solution than {\tt Trixi.jl} on the same mesh.

\subsubsection{Sedov's blast test}
\begin{figure}
\centering
\begin{tabular}{cc}
\includegraphics[width=0.45\textwidth]{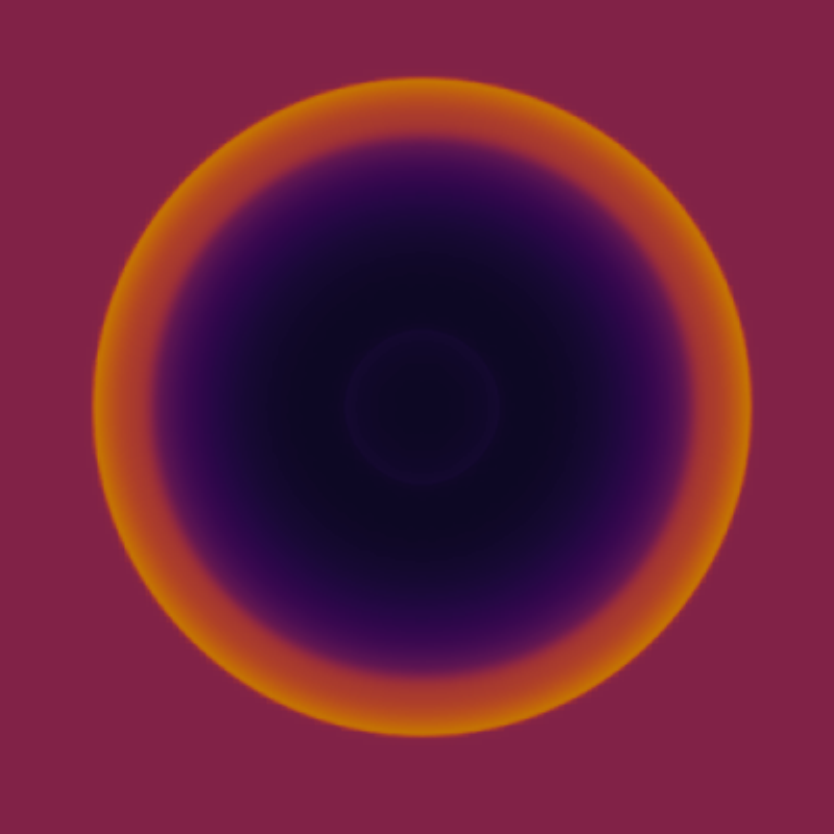} & \qquad \includegraphics[width=0.45\textwidth]{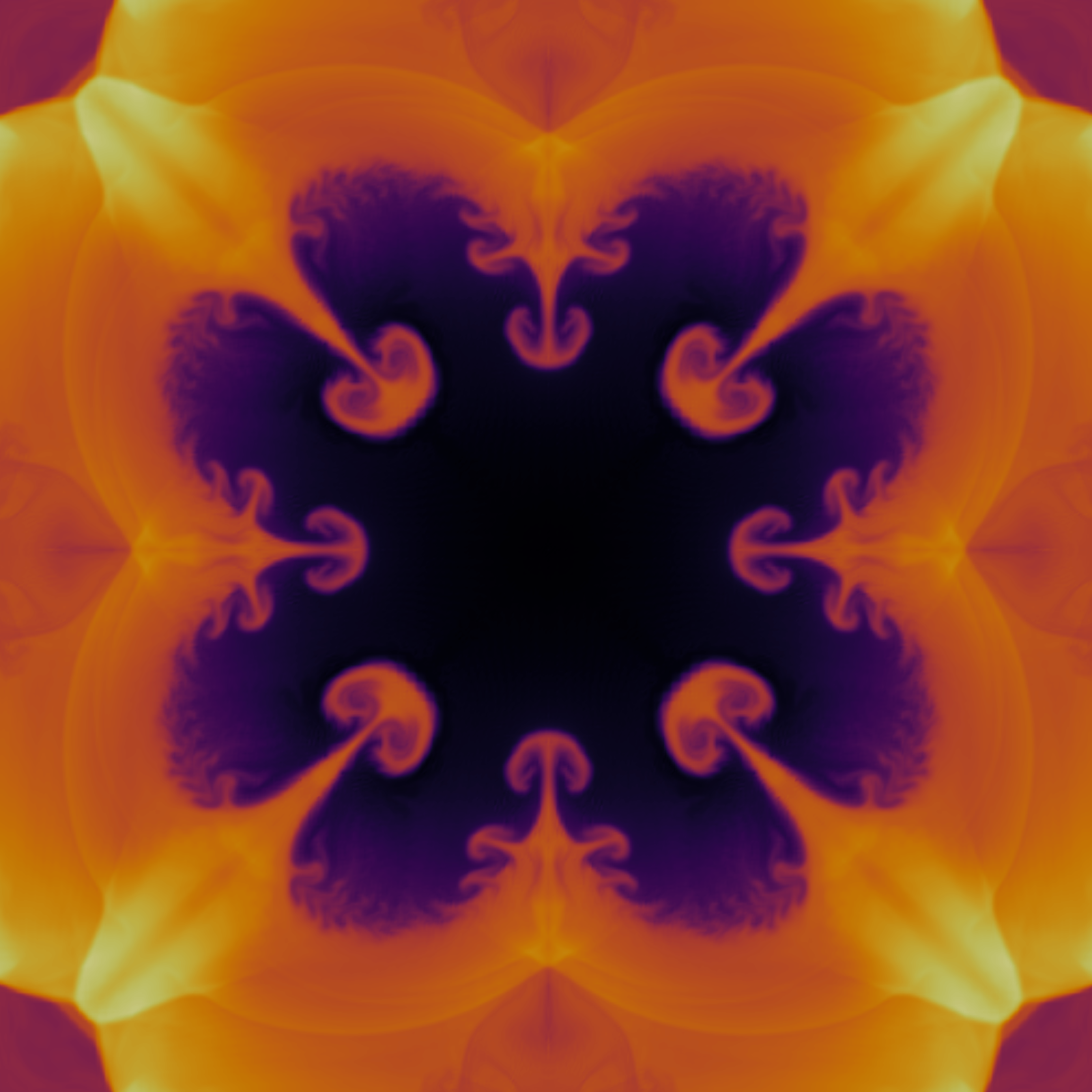} \\
(a) $t=2$ & (b) $t=20$
\end{tabular}
\caption{Sedov's blast test with periodic domain, density plot of numerical solution on $128 \times 128$ mesh in log scales with degree $N=3$ at (a) $t=2$ and (b) $t=20$ with polynomial degree $N=3$ computed using {\tt Trixi.jl}.}
\label{fig:blast.periodic.reference}
\end{figure}

\begin{figure}
\centering
\begin{tabular}{cc}
\includegraphics[width=0.45\textwidth]{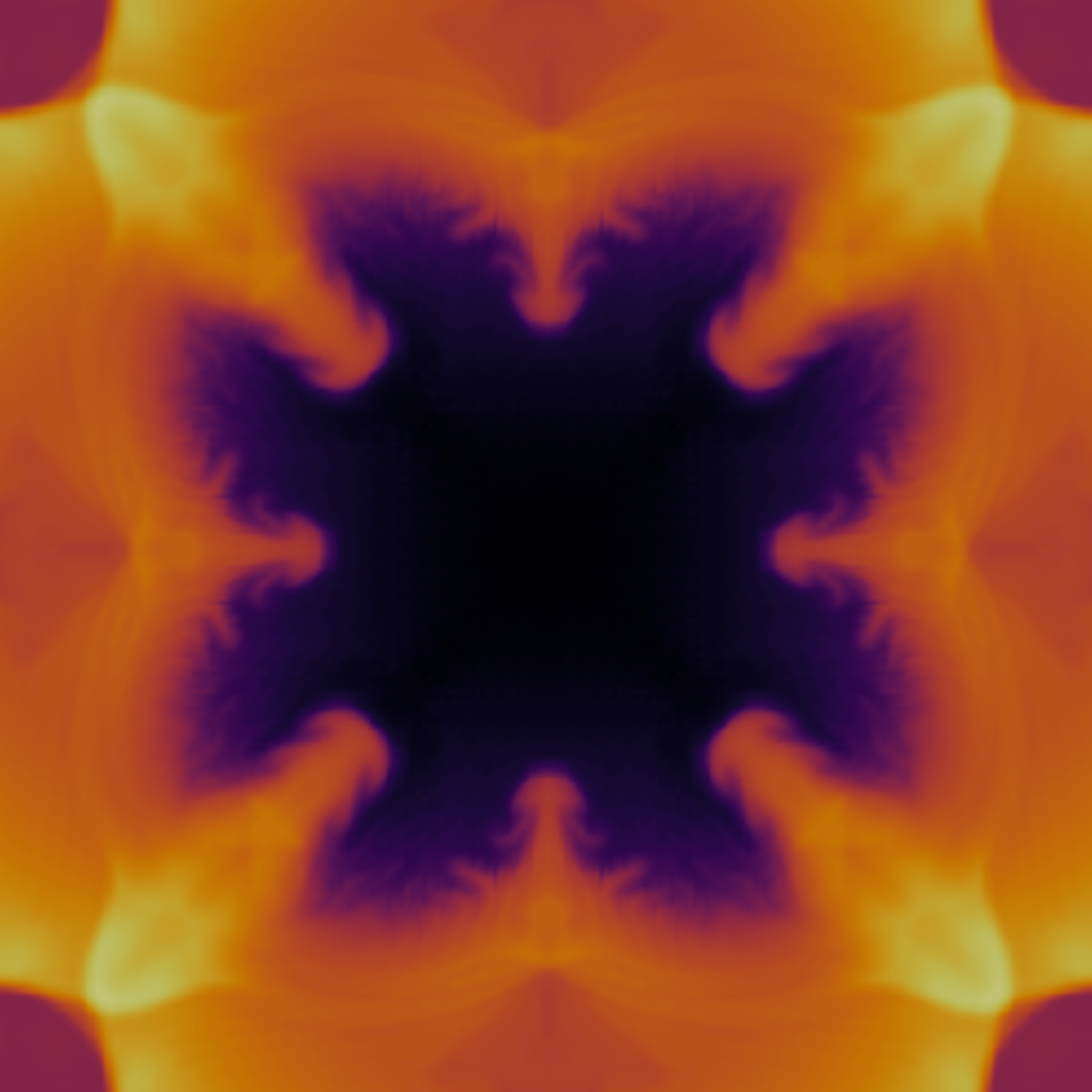} & \qquad \includegraphics[width=0.45\textwidth]{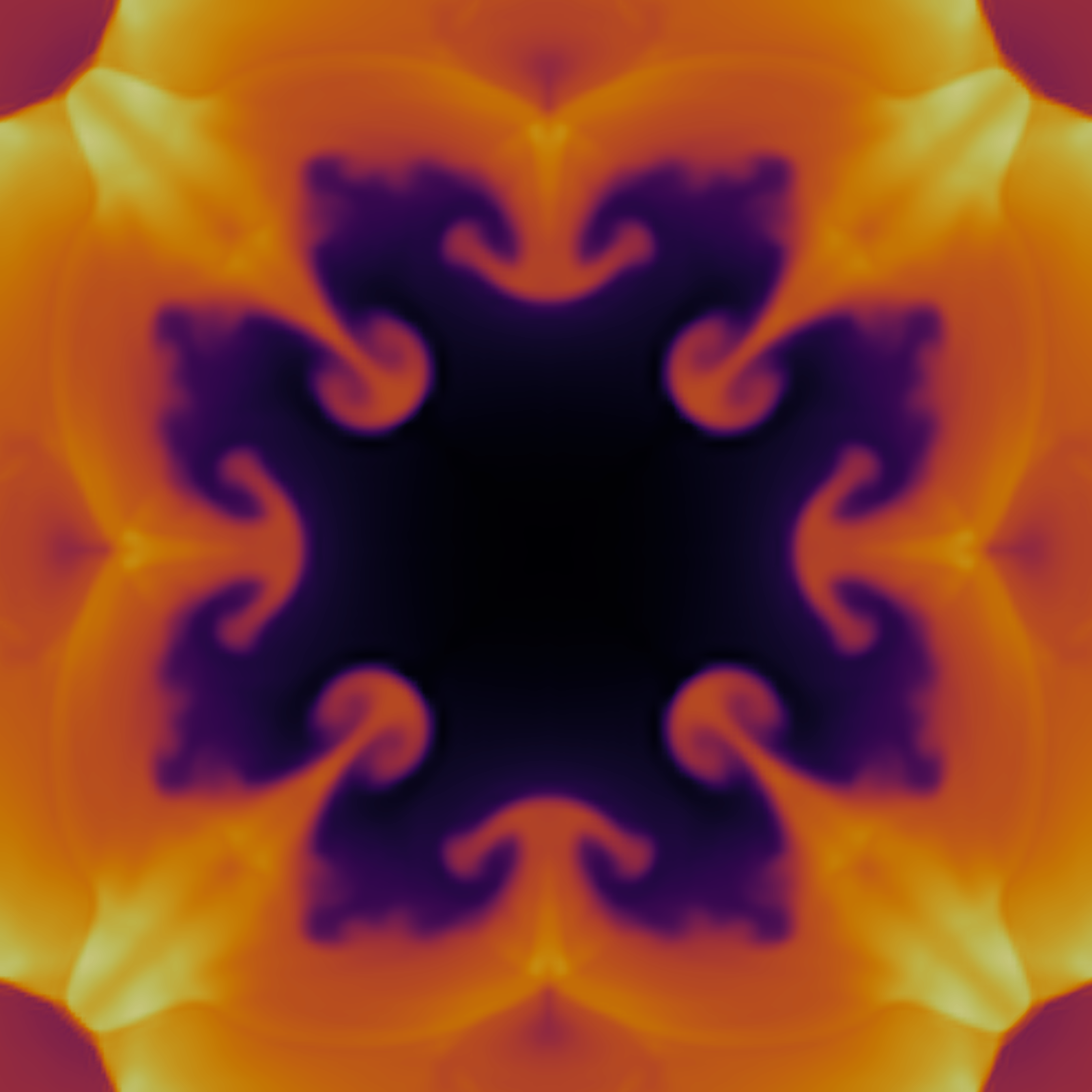} \\
(a) {\tt Trixi.jl} & (b) cRK-MH
\end{tabular}
\caption{Sedov's blast test with periodic domain, density plot of numerical solution on $64 \times 64$ mesh in log scales at $t=20$ with degree $N=3$.}
\label{fig:blast.periodic}
\end{figure}

This Sedov's blast test was introduced in~\cite{ramirez2021} to demonstrate the positivity preserving property of their scheme. Similar to Sedov's blast test in Section~\ref{sec:sedov.blast.1d},
this test on domain $[-1.5,1.5]^2$ has energy concentrated at the origin. More precisely, for the initial condition, we assume that the gas is at rest ($v_1 = v_2 = 0$) with Gaussian distribution of density and pressure
\begin{equation}
\rho(x,y) = \rho_0 + \frac{1}{4\pi\sigma_\rho^2} \exp \left( -\frac{r^2}{2\sigma_\rho^2} \right), \qquad p(x,y) = p_0  + \frac{\gamma - 1}{4 \pi \sigma_p^2} \exp\left( -\frac{r^2}{2\sigma_p^2} \right), \qquad r^2 = x^2 + y^2,
\end{equation}
where $\sigma_\rho = 0.25$ and $\sigma_p = 0.15$. The ambient density and ambient pressure are set to $\rho_0 = 1$, $p_0 = 10^{-5}$. The boundaries are taken to be periodic, leading to interaction between shocks as they re-enter the domain through periodic boundaries and form small scale fractal structures. In Figure~\ref{fig:blast.periodic}, we compare the density profiles of the numerical solutions of polynomial degree $N=3$ on a mesh of $64^2$ elements using {\tt Trixi.jl} and the proposed MUSCL-Hancock blending scheme in log scales. Looking at the reference solution on a finer $128^2$ element mesh (Figure~\correction{\ref{fig:blast.periodic.reference}}), we see that the MUSCL-Hancock scheme resolves the small scale structures better. This is particularly seen in the resolution of \textit{mushroom-like} structures that emanate from the corners of the black region in the pseudocolor plot.

\subsubsection{Astrophysical jet} \label{sec:m2000}
\begin{figure}
\centering
\begin{tabular}{cc}
\includegraphics[width=0.45\textwidth]{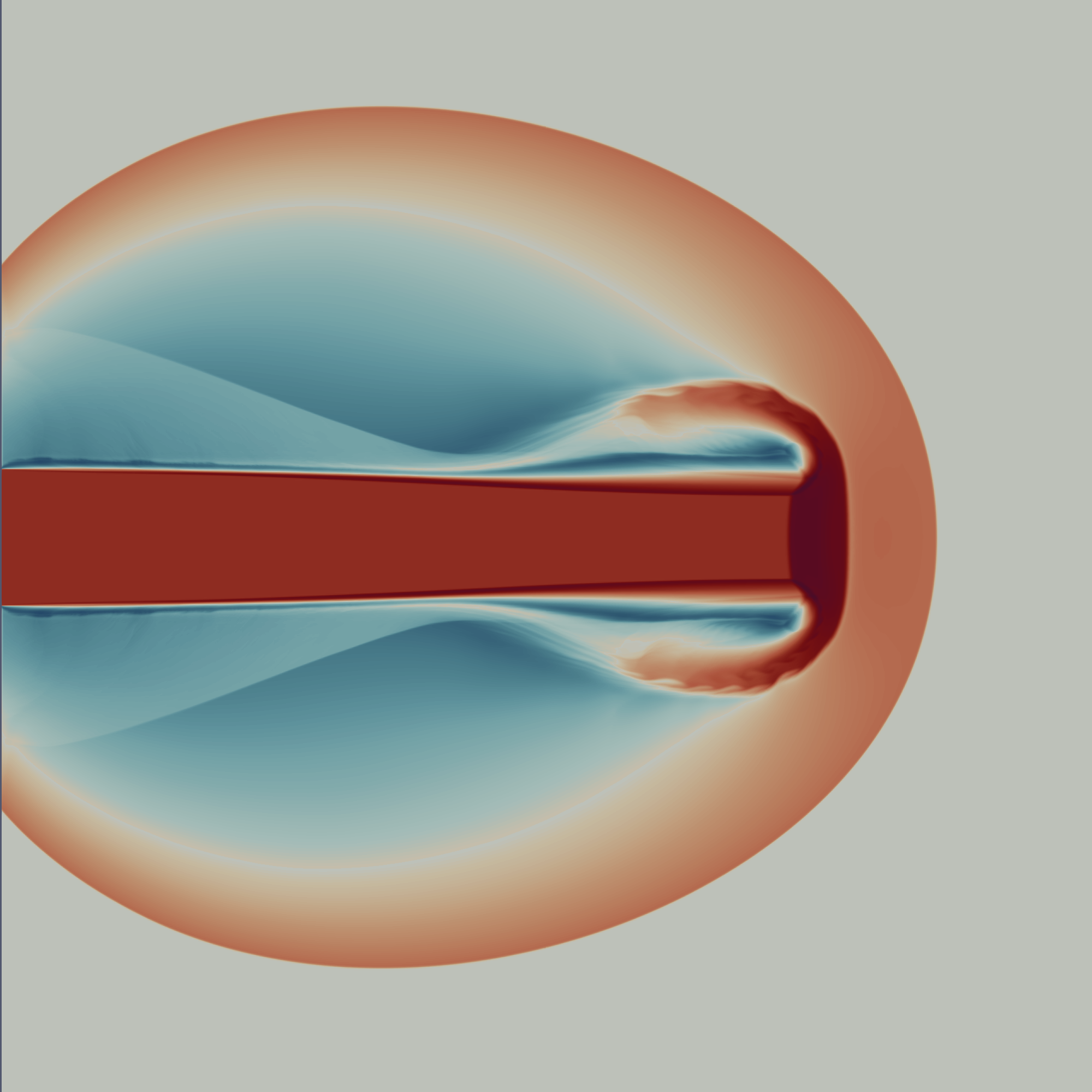} & \qquad \includegraphics[width=0.45\textwidth]{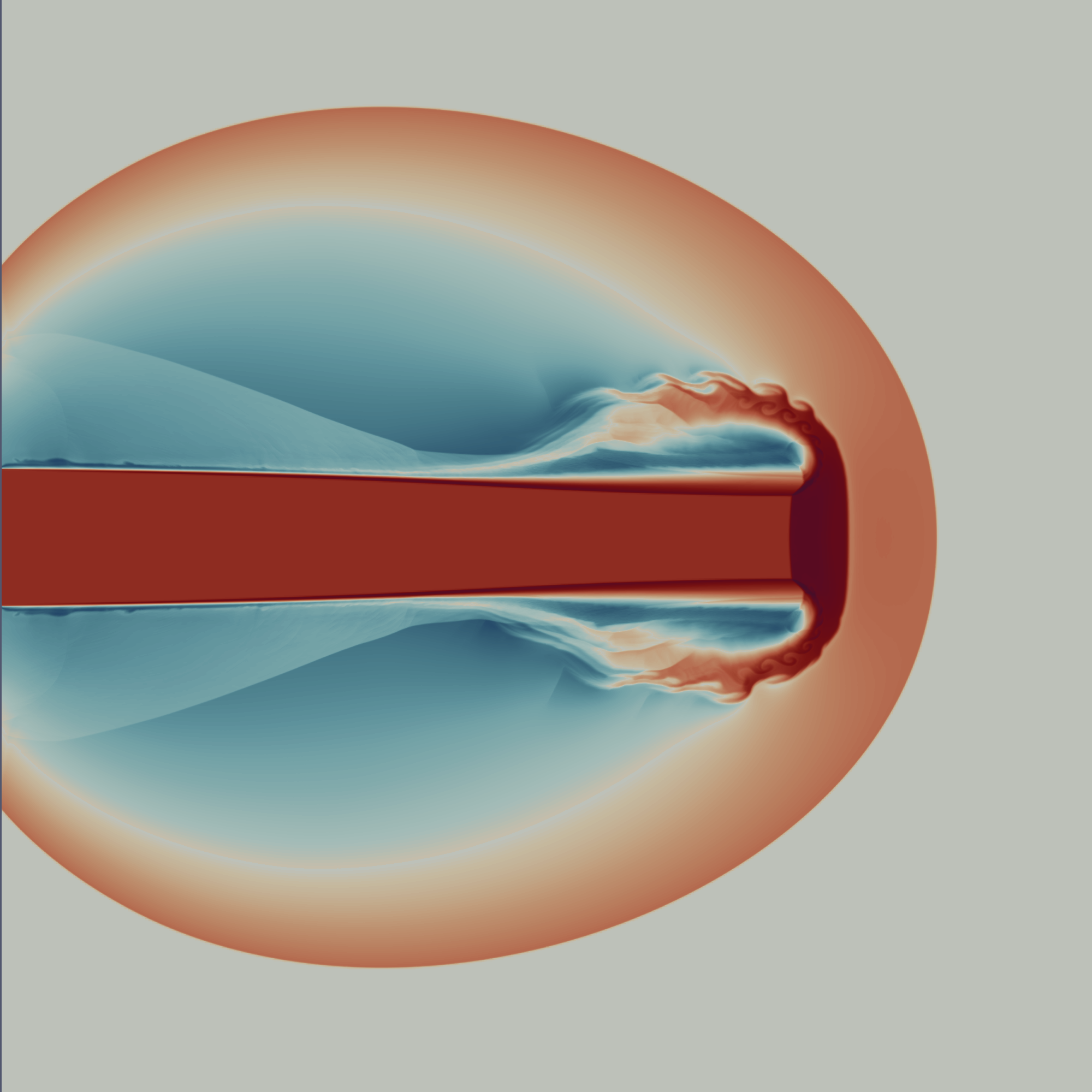} \\
(a) {\tt Trixi.jl} & (b) cRK-MH
\end{tabular}
\caption{Mach 2000 astrophysical jet, density plot of numerical solution in log scales on $400 \times 400$ mesh at time $t=0.001$.}
\label{fig:astrophysical.jet}
\end{figure}

In this test, a hypersonic jet is injected into a uniform medium with a Mach number of 2000 relative to the sound speed in the medium. The very high Mach number makes it a difficult test for admissibility preservation. In addition, the jet flow causes shear structures whose resolution is used as a test for accuracy of the blending limiter. Following~\cite{ha2005, zhang2010c}, the domain is taken to be $[0,1]\times[-0.5,0.5]$, the ambient gas in the interior has state $\boldsymbol u_{a}$ defined in primitive variables as
\[
(\rho,v_1, v_2, p)_a = (0.5, 0, 0, 0.4127)
\]
and jet inflow state $\boldsymbol u_j$ is defined in primitive variables as
\[
(\rho, v_1, v_2, p)_j = (5, 800, 0, 0.4127)
\]
On the left boundary, we impose the boundary conditions
\[
\boldsymbol u_b =
\begin{cases}
\boldsymbol u_a, & \quad \text{if} \quad y \in [-0.05, 0.05],  \\
\boldsymbol u_j, & \quad \text{otherwise},
\end{cases}
\]
and outflow conditions on the right, top and bottom. The HLLC numerical flux was used in the left most elements to distinguish between characteristics entering and exiting the domain (Section~\ref{sec:boundary}). To get better resolution of vortices, we used a smaller time step with $C_s=0.5$ in~\eqref{eq:time.step.2d} and included ghost elements in time step computation to handle the cold start. A more general way to compute time step sizes in problems with such stiff boundary conditions is to use embedded error based time stepping~\cite{babbar2025,Ranocha2021,ranocha2023}. This is a tough test for admissibility preservation because the kinetic energy is substantially higher than the internal energy, making it easy for negative pressure values to occur. At the same time, a Kelvin-Helmholtz instability arises before the bow shock. Thus, it is a good test both for admissibility preservation and  the ability to capture small-scale structures. The simulation gives negative pressures if positivity correction is not applied. While the large scale structures are captured similarly by both the schemes as seen in Figure~\ref{fig:astrophysical.jet}, the cRKFR with MH blending scheme shows more small scales near the front of the jet.

\subsubsection{Detonation shock diffraction}
\begin{figure}
\centering
\begin{tabular}{cc}
\includegraphics[width=0.45\textwidth]{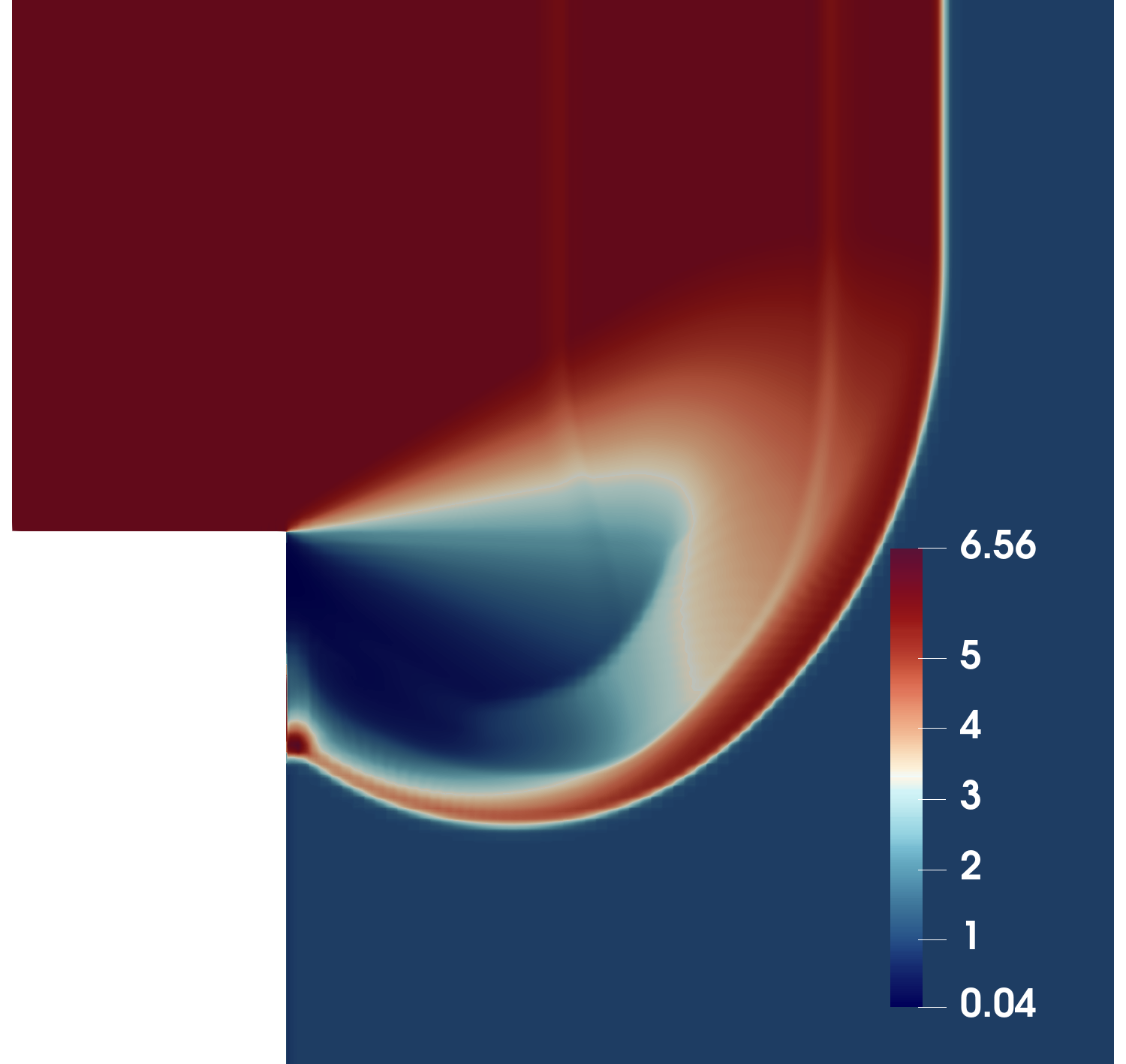} & \qquad \includegraphics[width=0.45\textwidth]{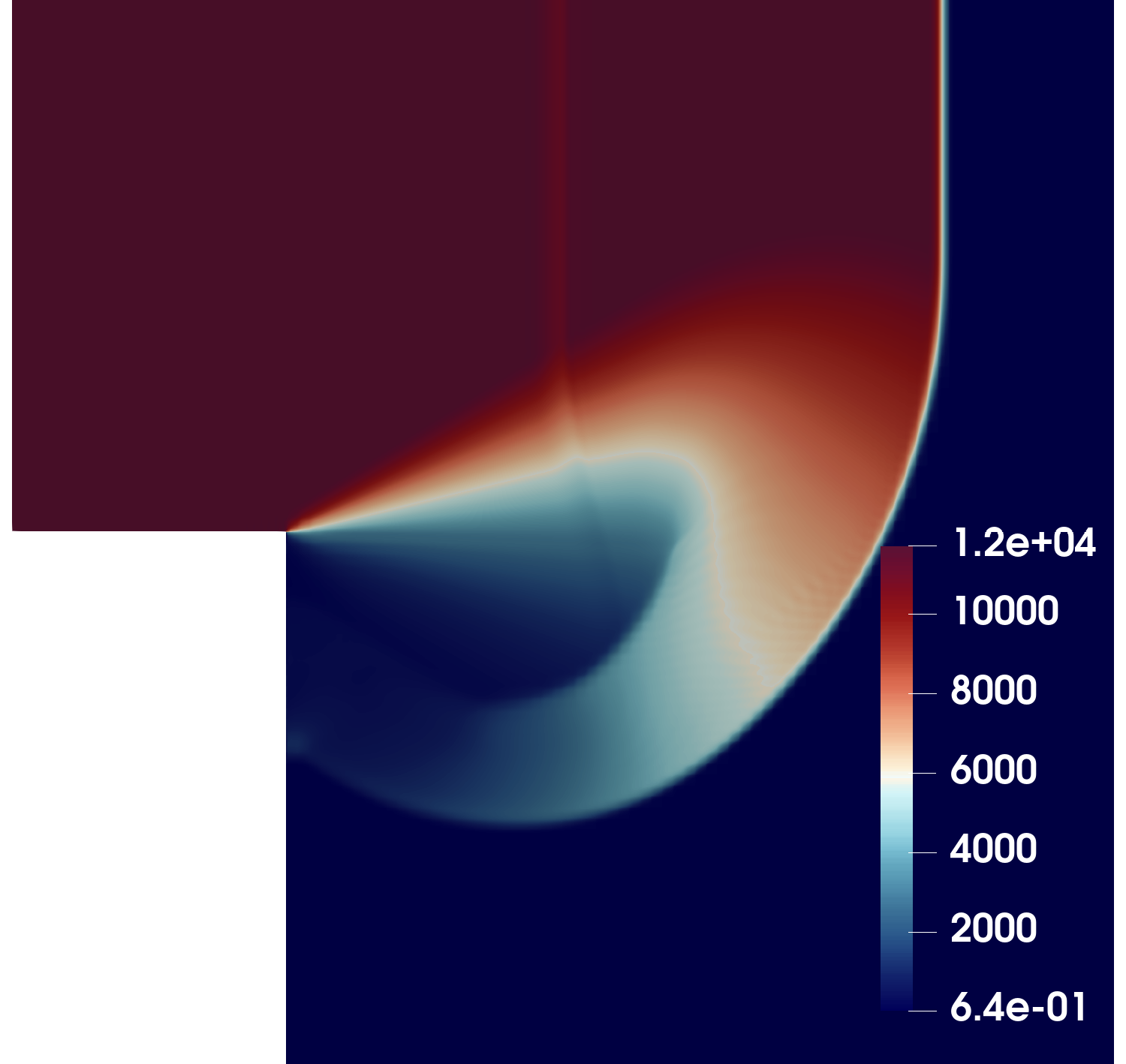} \\
(a) Density & (b) Pressure
\end{tabular}
\caption{Shock diffraction test, numerical solution at $t=0.01$ with degree $N=3$ and grid spacing $\Delta x = \Delta y = 1/100$. (a) Density profile, (b) Pressure profile.}
\label{fig:shock.diffraction}
\end{figure}

This detonation shock diffraction test~\cite{Takayama1991} is another test with a strong (Mach 100) shock, and is thus for validation of the admissibility preservation of our scheme. It involves a planar detonation wave that interacts with a wedge-shaped corner and diffracts around it, resulting in a complicated wave pattern comprising of trasmitted and reflected shocks, as well as rarefaction waves. The computational domain is  $\Omega = [0,2]^2 \backslash ([0,0.5] \times [0,1])$ and following~\cite{hennemann2021}, the simulation is performed by taking the initial condition to be a pure right-moving shock with Mach number of 100 initially located at $x=0.5$ and travelling through a channel of resting gas. The post shock states are computed by normal relations~\cite{naca1951}, so that the initial data is
\begin{alignat*}{2}
\rho(x,y) &= \begin{cases}
5.9970, \qquad & \text{if } x \le 0.5\\
1, & \text{if } x > 0.5
\end{cases}, \qquad
&& u(x,y) = \begin{cases}
98.5914, \qquad & \text{if } x \le 0.5 \\
0, & \text{if } x > 0.5
\end{cases} \\
v(x,y) &= 0, \qquad
&& p(x,y) = \begin{cases}
11666.5, \qquad & \text{if } x \le 0.5 \\
1, & \text{if } x > 0.5
\end{cases} \end{alignat*}
The left boundary is set as inflow and right boundary is set as outflow, all other boundaries are solid walls. The numerical results at $t=0.01$ with polynomial degree $N=3$ on a Cartesian grid consisting of uniformly sized squares with $\Delta x = \Delta y = 1/100$ are shown in Figure~\ref{fig:shock.diffraction}. The results look similar to~\cite{hennemann2021}.

\subsection{Problems with source terms}
To verify the effectiveness of the proposed time-averaging framework in the cRKFR scheme, we demonstrate its capability to handle conservation laws with source terms
\begin{equation} \label{eq:2d.conlaw.source}
\uu_t + \pf_x + \pg_y = \bss(\uu).
\end{equation}
The idea to extend the cRKFR framework to solve~\eqref{eq:2d.conlaw.source} is based on taking time averages of the source terms, and is described in Appendix~\ref{app:source.term}. The usage of source terms for provable admissibility preservation requires additional limiting on the time average sources~\cite{babbar2024generalized}, but we did not need any limiting beyond the flux limiter discussed in Section~\ref{sec:flux.limiter}. The problems also use the same CFL safety factor $C_s = 0.98$ in~\eqref{eq:time.step.2d} as was used in the problems without source terms.

\subsubsection{Euler's equations with gravity: Rayleigh Taylor instability}
\begin{figure}
\centering
\includegraphics[width=0.8\textwidth]{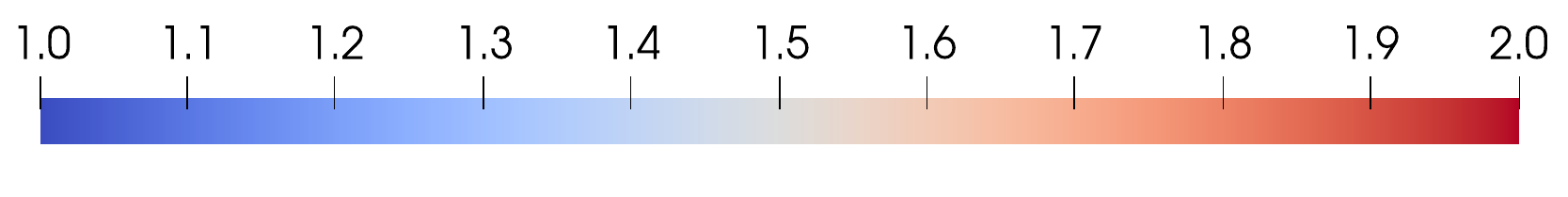} \\
\begin{tabular}{ccccc}
\includegraphics[width=0.173\textwidth]{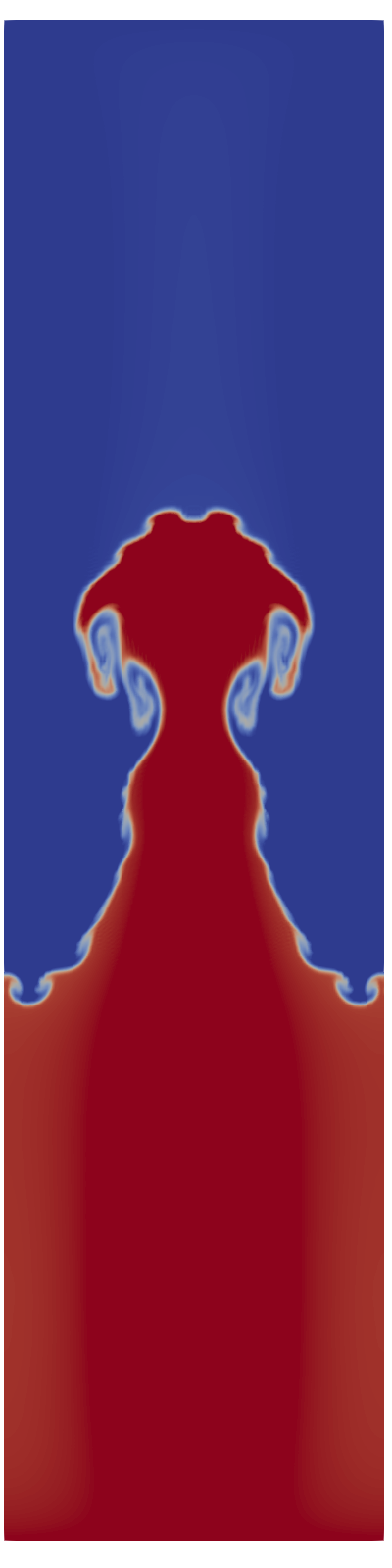} &
\includegraphics[width=0.173\textwidth]{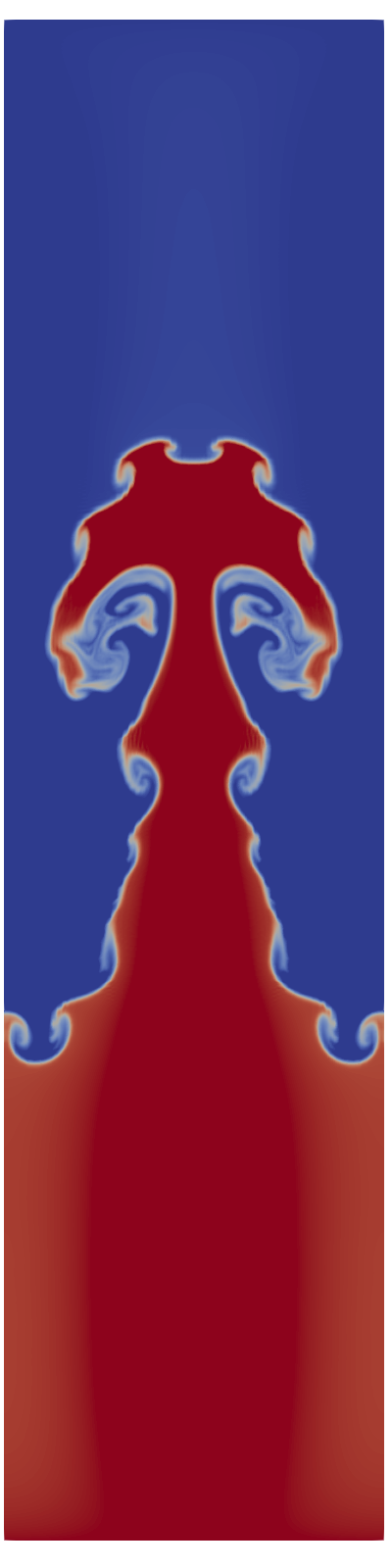} &
\includegraphics[width=0.173\textwidth]{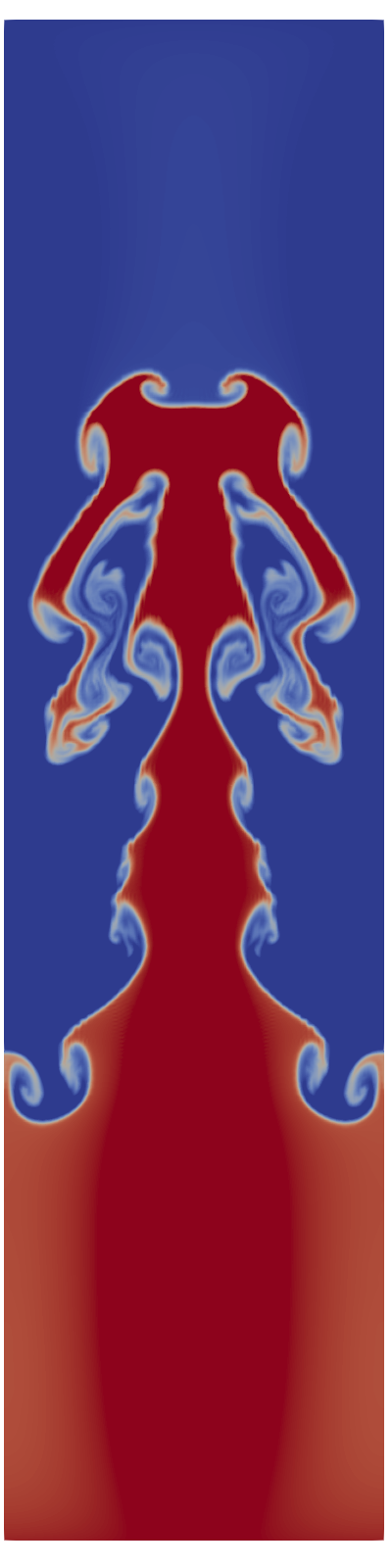} &
\includegraphics[width=0.173\textwidth]{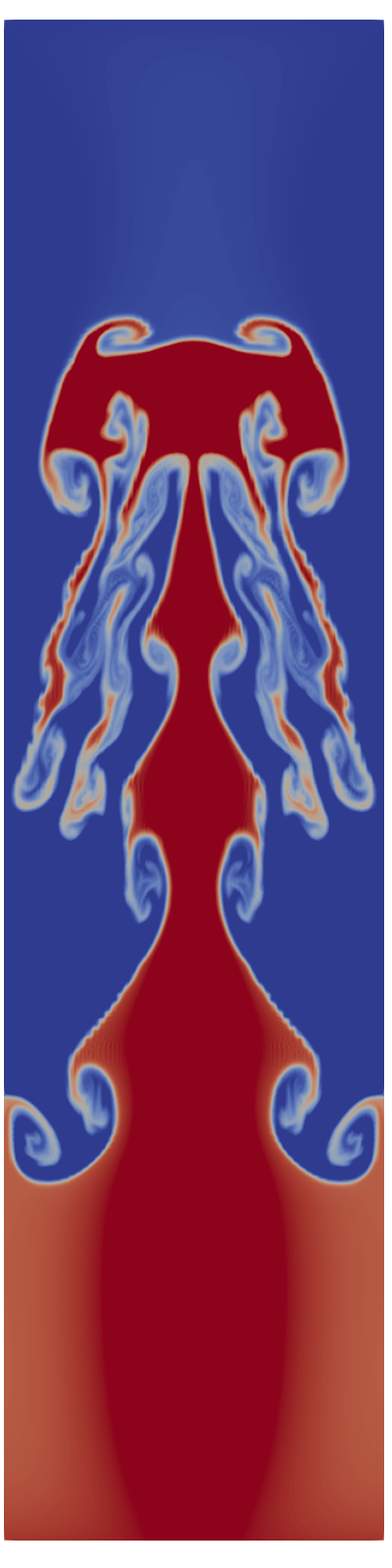} &
\includegraphics[width=0.173\textwidth]{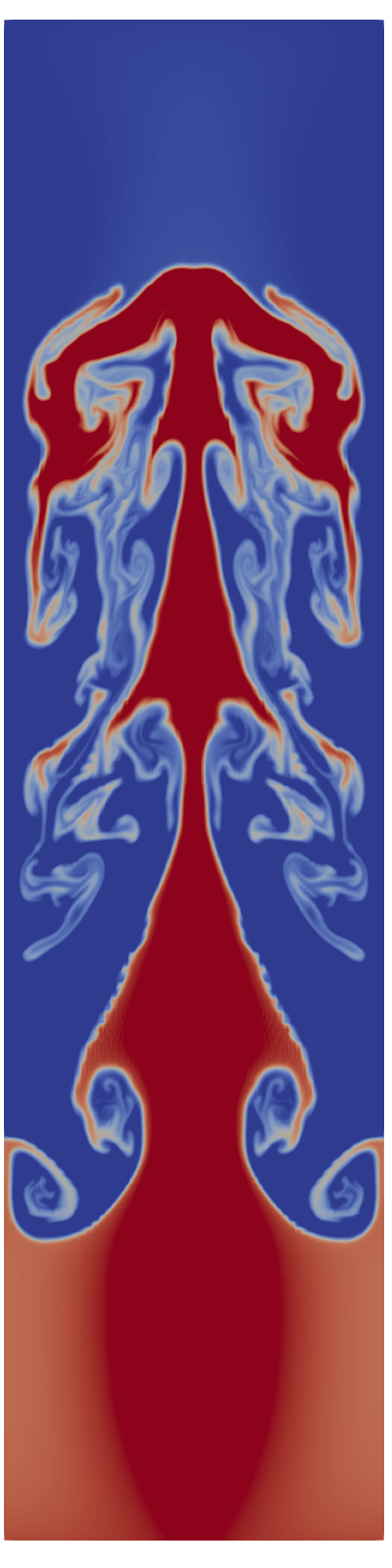} \\
(a) & (b) & (c) & (d) & (e)
\end{tabular}
\caption{Density profile of the Rayleigh Taylor instability problem at times $t=$ (a) 1.5, (b) 1.75, (c) 2, (d) 2.25, (e) 2.5 on a $64 \times 256$ mesh with polynomial degree $N=3$.} \label{fig:rayleigh.taylor}
\end{figure}

The first problem with source terms we consider are the Euler's equations with gravity. In this case, the equations~\eqref{eq:2deuler} are considered with an additional source term~\eqref{eq:2d.conlaw.source} given by $\bss = (0, 0, \rho, \rho v_1)$, where the gravitational acceleration is taken to be $g = 1$ and the direction of gravity is from bottom to top. We test our scheme with the Rayleigh Taylor instability problem, which is a classical problem in fluid dynamics that describes the instability of an interface between two fluids of different densities under the influence of gravity. The initial condition is given by
\[ (\rho, u, v, p) = \begin{cases}
(2, 0, - 0.025 c \cos (8 \pi x), 2 y + 1), \quad & y \leq 0.5,\\
(1, 0, - 0.025 c \cos (8 \pi x), y + 1.5),  & y > 0.5,
\end{cases} \]
where $c=\sqrt{\gamma p / \rho}$ is the sound speed, $\gamma = 5/3$ and the computational domain is $[0,0.25] \times [0,1]$. The left and right boundary conditions are solid walls~\eqref{eq:reflect.bc}. We use HLLC flux to distinguish between inflow and outflow characteristics at the top at bottom (Section~\ref{sec:boundary}). At the top, we specify the inflow values to be $(\rho, v_1, v_2, p) = (1, 0, 0, 2.5)$. At the bottom, we specify the inflow values to be $(\rho, v_1, v_2, p) = (2, 0, 0, 1)$. The simulation is performed using polynomial degree $N=3$ on a mesh of $64 \times 256$ elements. The density profiles at different times are shown in Figure~\ref{fig:rayleigh.taylor}.

\subsubsection{Ten moment problem: two rarefaction problem and a realistic simulation}

The ten moment problem can be seen as a generalization of Euler's equations~\eqref{eq:2deuler} where the energy and pressure are tensor quantities. Specifically, the energy tensor is defined by the ideal equation of state $\ccE = \frac{1}{2}  \ccP + \frac{1}{2} \rho \bv \otimes \bv$ where $\rho$ is the density, $\bv$ is the velocity vector, $\ccR$ is the symmetric pressure tensor. The conservative variables are $\uu =(\rho, \rho v_1, \rho v_2,  \cE_{11},  \cE_{12},  \cE_{22})$, and the fluxes are given by
\begin{equation}
\pf = \left[\begin{array}{c}
\rho v_1\\
\cR_{11} + \rho v_1^2\\
\cR_{12} + \rho v_1 v_2\\
\left( \cE_{11} + \cR_{11} \right) v_1\\
\Eonetwo v_1 + \frac{1}{2}  ( \cR_{11} v_2 + \cR_{12} v_1 )\\
\Etwotwo v_1 + \cR_{12} v_2
\end{array}\right], \quad \pg = \left[\begin{array}{c}
\rho v_2\\
\cR_{12} + \rho v_1 v_2\\
\cR_{22} + \rho v_2^2\\
\cE_{11} v_2 + \cR_{12} v_1\\
\cE_{12} v_2 + \frac{1}{2}  ( \cR_{12} v_2 + \cR_{22} v_1 )\\
\left( \cE_{22} + \cR_{22} \right) v_2
\end{array}\right]. \label{eq:tmp}
\end{equation}
The source terms are given by $\bss = \bss^x + \bss^y$, where
\begin{equation}\label{eq:tenmom.source}
\bss^{x} = \left[\begin{array}{c}
0\\
- \frac{1}{2} \rho \partial_x W\\
0\\
- \frac{1}{2} \rho v_1 \partial_x W\\
- \frac{1}{4} \rho v_2 \partial_x W\\
0
\end{array}\right], \qquad  \bss^{y} = \left[\begin{array}{c}
0\\
0\\
- \frac{1}{2} \rho \partial_y W\\
0\\
- \frac{1}{4} \rho v_1 \partial_y W\\
- \frac{1}{2} \rho v_2 \partial_y W
\end{array}\right],
\end{equation}
where $W (x, y, t)$ is a given function, which models electron quiver energy in the laser~\cite{Berthon2015}. The admissibility set is given by
\[
\Uad = \left\{ \uu \in \re^6 | \rho \left( \uu \right) > 0, \quad \bs{x}^T
\ccP \left( \uu \right)  \bs{x} > 0, \quad \bs{x} \in \mathbb{R}^2 \backslash
\left\{ \bzero \right\} \right\}
\]
These equations are relevant in many
applications, especially related to plasma
flows~(see~{\cite{Berthon_TMP_2006,Berthon2015}} and further
references in~{\cite{meena2017}}), in cases where the \textit{local
	thermodynamic equilibrium} used to close the Euler equations of compressible
flows is not valid, and anisotropic nature of the pressure needs to be
accounted for. The theoretical study of this model, including the hyperbolicity and admissibility properties, can be found in~\cite{Meena_Kumar_Chandrashekar_2017,meena2017,Meena2020}. By casting our scheme in the time average flux and time average sources framework in Appendix~\ref{app:source.term}, we can use the admissibility preserving framework for the ten moment equations using cRKFR schemes as is done for the LWFR schemes in~\cite{babbar2024generalized}. We demonstrate the capability of our scheme by showing two test cases for the ten moment problem.

\begin{figure}
\centering
\begin{tabular}{cc}
\includegraphics[width=0.45\textwidth]{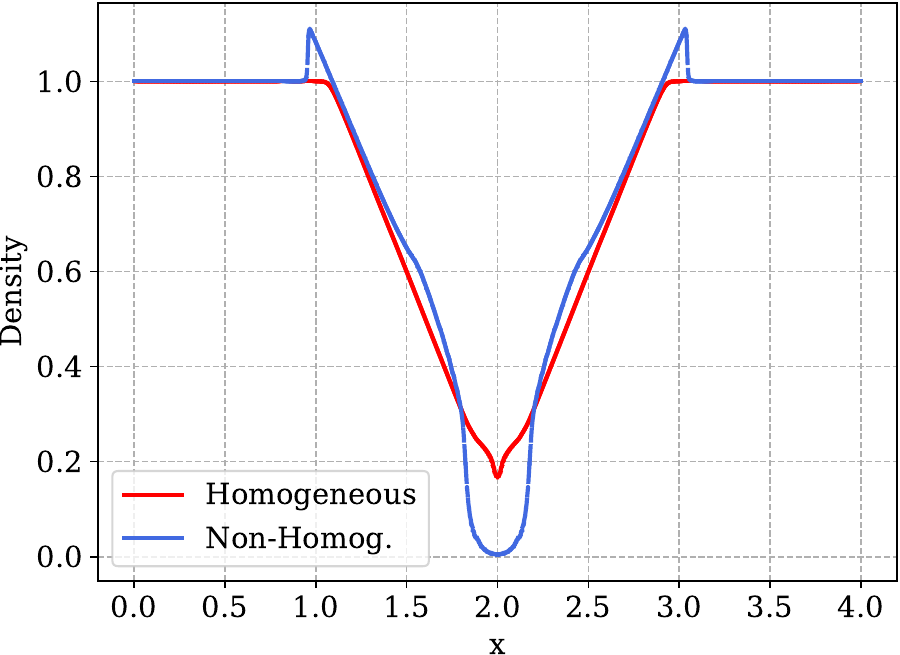} & \qquad \includegraphics[width=0.45\textwidth]{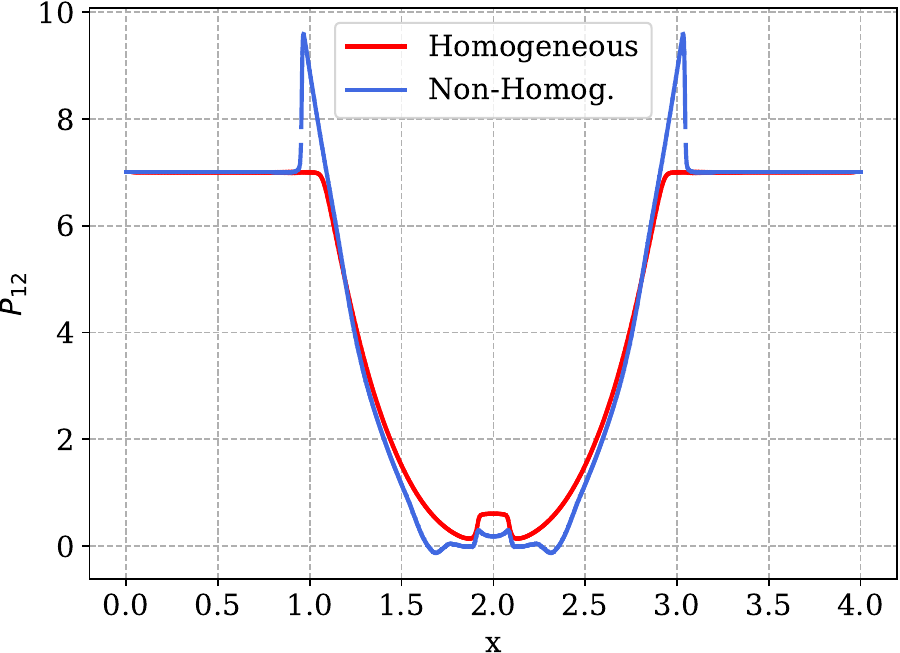} \\
(a) & (b)
\end{tabular}
\caption{Two rarefaction problem for Ten Moment Problem comparing solutions with and without source terms using polynomial degree $N=3$ on a grid of 500 elements. (a) Density, (b) $P_{12}$ profile are shown at time $t=0.1$}
\label{fig:tmp.riemann.problem}
\end{figure}

\paragraph{Two rarefaction problems. } The initial condition for the Riemann problem is given by
\[
\left( \rho, v_1, v_2, \Poneone, \Ponetwo, \Ptwotwo \right) =
\left\{\begin{array}{ll}
(1, - 4, 0, 9, 7, 9), \qquad & x < 0,\\
(1, 4, 0, 9, 7, 9), & x > 0,
\end{array}\right.
\]
with source terms as in~\eqref{eq:tenmom.source} with $W (x, y, t) = 25
\exp (- 200 (x - 2)^2)$. We show the numerical solutions with degree $N = 3$
and 500 elements at $t = 0.1$ in Figure~\ref{fig:tmp.riemann.problem} with and
without the source terms using the blending limiter. The solution with source
terms has a near vacuum state at the centre. Thus, this is a test where low
density is caused by the presence of source terms verifying that our scheme is
able to capture admissibility even in the presence of source terms. The
positivity limiter had to be used to maintain admissibility of the solution.

\begin{figure}
\centering
\begin{tabular}{cc}
\includegraphics[width=0.33\textwidth]{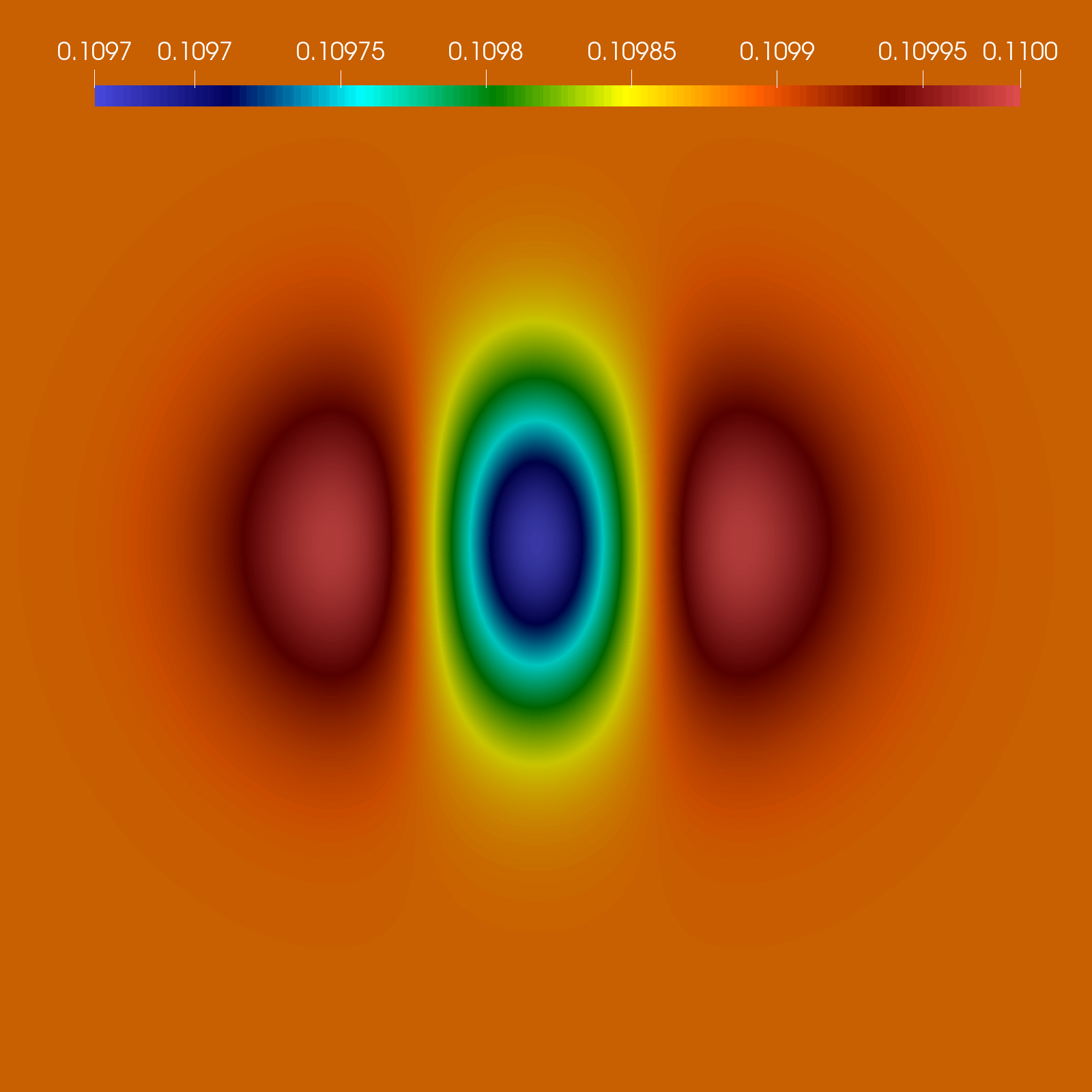} & \qquad \includegraphics[width=0.45\textwidth]{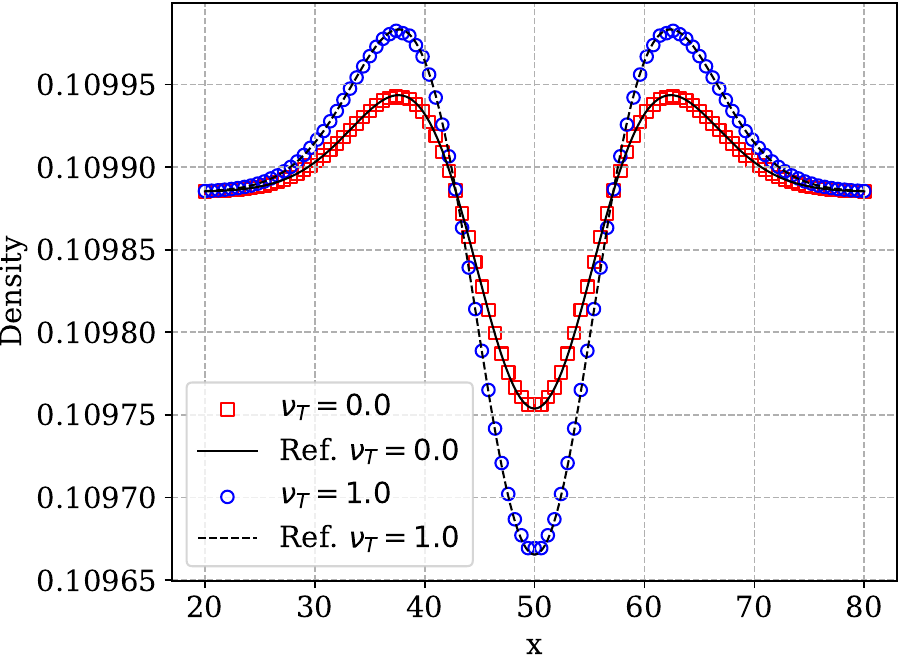} \\
(a) & (b)
\end{tabular}
\caption{Realistic simulation with different absorption coefficients, density profile obtained using a polynomial degree $N=3$ on a grid of 30 elements. (a) Pseudocolor plot, (b) Line cut along $y=50$.}
\label{fig:tmp.realistic}
\end{figure}

\paragraph{Realistic simulation with inverse bremsstrahlung.} Consider the domain $\Omega = [0, 100]^2$ with outflow boundary conditions (Section~\ref{sec:boundary}). The uniform initial condition is taken to be
\[
\rho = 0.109885, \quad \vone = \vtwo = 0, \quad \Poneone = \Ptwotwo = 1,
   \quad \Ponetwo = 0
\]
with the electron quiver energy $W (x, y, t) = \exp (- 0.01 ((x - 50)^2 + (y -
50)^2))$. The source term is taken from~{\cite{Berthon2015}}, and only has the
$x$ component, i.e., $\bss^y \left( \uu \right) = \bzero$, even though $W$
continues to depend on $x$ and $y$. An additional source corresponding to
energy components $\bss_E = (0, 0, 0, \nu_T \rho W, 0, 0)$ is also added where
$\nu_T$ is an absorption coefficient. Thus, the source terms are $\bss =
\bss_x + \bss_E$. The simulation is run till $t = 0.5$ on a grid of 30 elements with polynomial degree $N=3$ using the blending limiter. The density plot with a cut at $y=50$ is shown in Figure~\ref{fig:tmp.realistic}.

\section{Summary and conclusions} \label{sec:conclusion}
The compact Runge-Kutta (cRK) method of~\cite{chen2024}, originally introduced for Discontinuous Galerkin (DG) methods, has been extended to the Flux Reconstruction (FR) framework of~\cite{Huynh2007}. A key feature of the proposed cRKFR method is the interpretation of cRK as a procedure to approximate time-averaged fluxes. With this interpretation, three dissipation models were discussed for the computation of the \textit{time averaged numerical flux}. The original cRKDG scheme of~\cite{chen2024}, equivalent to our D-CSX dissipation model, required $s$ numerical fluxes to be computed at each time step, where $s$ is the number of inner stages of the RK method. The D1 and D2 dissipation models (termed so following~\cite{babbar2022}) required computation of a single numerical flux computation per time step. While the D1 dissipation model has lower CFL numbers than~\cite{chen2024}, the D2 dissipation is equivalent to the D-CSX dissipation of~\cite{chen2024} for linear problems, and thus has the same CFL numbers. The accuracy of the dissipation models is also seen to nearly be the same. The proposed approach thus lowers the total volume of interelement communication required without a loss in accuracy. In addition to the discussion on the dissipation models, we also discussed the {\evaluate} and {\extrapolate} schemes to compute trace values of the time averaged flux needed for the central part of the numerical flux. The {\evaluate} scheme performs the cRK procedure at the element interfaces, in contrast to the {\extrapolate} scheme where the trace values are obtained by extrapolating the flux. It is shown that the {\evaluate} scheme obtains optimal order of accuracy for all polynomial degrees, while the {\extrapolate} scheme gives suboptimal order of accuracy of $O(h^{N+1/2})$ for some nonlinear problems where the degree $N$ is odd. The procedure to handle boundary conditions like inflow, outflow and solid walls for the proposed cRKFR scheme is also detailed. To handle problems with nonsmooth solutions, we use the subcell-based blending limiter of~\cite{babbar2024admissibility,hennemann2021}. This limiter blends a high-order FR method with a lower-order subcell update using a smoothness indicator. Following~\cite{babbar2024admissibility}, the method aims to enhance accuracy by using Gauss-Legendre solution points and perform MUSCL-Hancock reconstruction on subcells. Another crucial benefit of the time averaged flux viewpoint was the ability to use the flux limiter of~\cite{babbar2024admissibility}. The flux limiter of~\cite{babbar2024admissibility} applied to the time-averaged numerical flux preserves admissibility in means, and when combined with the positivity-preserving scaling limiter of~\cite{Zhang2010b}, results in a fully admissibility-preserving cRKFR scheme. The method is further extended to handle source terms following~\cite{babbar2024generalized} by incorporating their contributions as additional time averages, maintaining high order accuracy. Numerical experiments on Euler's equations and the ten-moment problem validate the robustness, accuracy, and admissibility preservation of the proposed method. The numerical results additionally show that the Wall Clock Time (WCT) performance of the proposed cRKFR scheme is better than Lax-Wendroff Flux Reconstruction scheme of~\cite{babbar2022,babbar2024admissibility}, which is another single stage evolution method based on time averages. This is the first paper of a series of our future works on cRKFR, which include extensions to implicit-explicit time stepping, embedded error based time stepping for curvilinear grids and an integration with the stage-dependent-RKDG (sdRKDG) methods proposed in \cite{chen2025runge}.

\section*{Acknowledgments}

The work of Arpit Babbar is funded by the Mainz Institute of Multiscale Modeling.

\appendix

\section{Time average flux viewpoint} \label{app:time.averaged.flux}
We now justify referring to $\F_h, \uU_h$~(\ref{eq:disc.avg.flux},~\ref{eq:time.avg.sol}) as the time average flux, solution approximations respectively by formal order of accuracy arguments. Although a mathematically rigorous justification for  error analysis is beyond the scope of this work and will be postponed  to our future work, there are existing results on the stability analysis and error estimates for schemes with Lax–Wendroff time discretization. The viewpoint of time averages will then be used to prove the \textit{linear equivalence} (equivalence for linear problems) of the cRKFR scheme and the LWFR scheme, as long as both use linearly equivalent dissipation models (Section~\ref{sec:numflux}). For the conservation law~\eqref{eq:con.law}, we consider an augmented system of time dependent equations given by\footnote{This augmentation is equivalent to solving an additional system of ODEs $\bv_t = \fg(\uu)$.}
\begin{equation}\label{eq:augmented.system}
\tu_t := (\uu, \bv)_t  = (-\pf(\uu)_x, \fg(\uu)) =: \bL(\tu),
\end{equation}
where $\bv$ is an auxiliary variable that will not appear in our main conclusion, and $\fg$ can be chosen to be any non-linear function. We will show that numerical solutions at the inner stages of the cRKFR scheme can be used to approximate the time average of a general $\fg(\uu)$. A Runge-Kutta (RK) discretization will be chosen by specifying its Butcher Tableau~\eqref{eq:butcher}, corresponding to which we will have a cRKFR scheme~\eqref{eq:crkfr}. The cRKFR scheme will then be used to solve $\uu_t = -\pf(\uu)_x$ in~\eqref{eq:augmented.system} to obtain the solution $\uu_h$. The auxiliary variable $\bv$ will also be discretized as $\bv_h$ at the same solution points as $\uu_h$ and the chosen RK scheme used will be used to solve the system of ODE $(\bv_h)_t = \fg(\uu_h)$. The chosen RK and cRKFR method are both formally order $N+1$ accurate. Thus, numerical solution of the system~\eqref{eq:augmented.system}, denoted $\tu_h$, will satisfy, for smooth solutions
\begin{equation}\label{eq:formal.analysis}
\tu_h^n + \Delta t \bL_h(\tu) = \tu^{n+1} + \mathcal{O}(\Delta t^{N+2}) \implies \tu^{n+1} - \tu^n = \Delta t \bL_h(\tu) + \mathcal{O}(\Delta t^{N+2}), \quad \bL_h = (-\partial_x \F_h, \fG_h)
\end{equation}
where $\tu^n, \tu^{n+1}$ denote the exact solution of the system~\eqref{eq:augmented.system} at time level $n, n+1$ respectively, $\F_h$ is as in~\eqref{eq:disc.avg.flux} and $\fG_h$ is similarly defined at each solution point as
\begin{equation}\label{eq:G.defn}
\fG_{e,p} = \sum_{i=1}^s b_i \pg({\uu}_{e,p}^{(i)}),
\end{equation}
$\{b_i\}$ are the coefficients of the final stage of the RK method~\eqref{eq:butcher}, $\uu_{e,p}^{(i)}$ are solutions of the cRKFR scheme at inner stages~\eqref{eq:crkfr}. Integrating the system~\eqref{eq:augmented.system} in time over $[t^n, t^{n+1}]$ gives
\[
\tu^{n+1}-\tu^n = \int_{t^n}^{t^{n+1}} \bL(\tu) \ud t,
\]
so that~\eqref{eq:formal.analysis} implies
\begin{equation} \label{eq:L.time.avg}
\bL_h = \frac{1}{\Delta t } \pdx \int_{t^n}^{t^{n+1}} \bL (\tu) \ud t + \mathcal{O}(\Delta t^{N+1}).
\end{equation}
Thus, taking the $\bv$ component~\eqref{eq:augmented.system} of~\eqref{eq:L.time.avg} gives us
\begin{equation}\label{eq:tG.approx}
\fG_h = \frac{1}{\Delta t } \int_{t^n}^{t^{n+1}} \fg(\uu(x_{e,p}, t)) \ud t + \mathcal{O}(\Delta t^{N+1}).
\end{equation}
Equation~\eqref{eq:tG.approx} shows that $\fG_h$~\eqref{eq:G.defn} is an order $N+1$ approximation of the time average of $\fg$ over $[t^n, t^{n+1}]$. By taking the specific choices $\pg(\uu) = \pf(\uu)$ and $\pg(\uu) = \uu$, we get that $\F_h$~\eqref{eq:disc.avg.flux} is an order $N+1$ approximation of the time average of $\pf$ and $\uU_h$~\eqref{eq:time.avg.sol} is an order $N+1$ approximation of the time average of $\uu$.
Even though this argument needed the augment system~\eqref{eq:augmented.system}, the final result~\eqref{eq:tG.approx} is only about the approximation $\fG_h$~\eqref{eq:tG.approx} defined in terms of the numerical solution $\uu_h$ of the original system~\eqref{eq:con.law} and does not contain the auxiliary variable $\bv$~\eqref{eq:augmented.system}.

We now use the time average viewpoint to see that the LWFR and cRKFR schemes are linearly equivalent for polynomial degree $N$ with the RK scheme~\eqref{eq:butcher} of order $N+1$. To be precise, we will show that the D1 dissipation model of Section~\ref{sec:numflux} is equivalent to the D1 dissipation LWFR scheme of~\cite{babbar2022}, and the D-CSX and D2 dissipation of Section~\ref{sec:numflux} are both equivalent to the D2 dissipation LWFR scheme of~\cite{babbar2022}.  We look at the time average approximations of the cRKFR schemes of order $N+1$, given in~(\ref{eq:crk22},~\ref{eq:crk33},~\ref{eq:crk44}), where $N$ is the degree of the polynomial basis. It is easy to see that, for the linear advection equation $\uu_t + a \uu_x = \bzero$, the time averaged flux and solution are given by
\begin{equation} \label{eq:linear.avg}
\F_h = a \uU_h^n, \qquad \uU_h^n = \sum_{m=0}^{N} \frac{(- a \Delta t)^m}{(m+1)!}\partial_x^m \uu_h^n.
\end{equation}
The time averaged flux of the Lax-Wendroff scheme is an approximation of the truncation of the time average flux given by
\[
\F = \sum_{m=0}^N \frac{\Delta t^m}{(m+1)!}\partial_t^m (a \uu).
\]
In the literature~\cite{Qiu2005b,babbar2022,babbar2024admissibility,babbar2025}, using $\partial_t^m\uu = (- a \partial_x)^m \uu$, the following approximation is used
\[
\F = a \sum_{m=0}^N \frac{(-a \Delta t)^m}{(m+1)!} \partial_x^m \uu_h,
\]
which is the same as the time average flux of the cRK scheme~\eqref{eq:linear.avg}. Thus, the time average flux approximation is the same for LWFR and the cRKFR schemes. Since the evolution of the solution of both schemes is given in terms of the same time average evolutions~\eqref{eq:crkfr.avg} (see equations~(4,~5) of~\cite{babbar2024admissibility}), the schemes will be equivalent as long as the numerical fluxes are computed in the same way. Thus, as long as both schemes use the same dissipation model~(\ref{eq:dcsx},~\ref{eq:D1},~\ref{eq:D2}), they will be equivalent. These arguments apply the cRK schemes described in~(\ref{eq:crk22},~\ref{eq:crk33},~\ref{eq:crk44}) since we can explicitly obtain the expression~\eqref{eq:linear.avg}. A heuristic argument for the general cRKFR scheme using an RK method order $N+1$ with polynomial degree $N$ can be given as follows. In this case, the cRKFR is formally of order $N+1$, and thus its time average flux will satisfy
\begin{equation} \label{eq:linear.avg.general}
\F = a \sum_{m=0}^N \frac{(-a \Delta t)^m}{(m+1)!} \partial_x^m \uu_h + \mathcal{O}(\Delta t^{N+1}).
\end{equation}
However, the $\mathcal{O}(\Delta t^{N+1})$ will consist only of derivatives of order greater than $N+1$ of the degree $N$ polynomial $\uu_h$, and will thus be zero. Therefore, the time average flux will be the same as~\eqref{eq:linear.avg}, completing the proof of linear equivalence. Equation~\eqref{eq:linear.avg.general} is a formal step, but it also can verified for particular RK methods.

\section{Source term treatment} \label{app:source.term}
Consider a conservation law with source term
\begin{equation}\label{eq:con.law.source.term}
\uu_t + \pf_x = \bss,
\end{equation}
where $\bss = \bss(\uu, x, t)$ is the source term. Similar to~\cite{babbar2024generalized} for the LWFR scheme, the idea to solve this equation with a cRK scheme is to use a time averaged source term $\bS$. For the scheme with source terms, the evolution will be given by
\begin{equation}\label{eq:crkfr.source}
\uu_h^{n+1} = \uu_h^n - \Delta t \pdx \F_h + \Delta t \bS_h^\delta,
\end{equation}
where $\F_h$ is the continuous time average flux~\eqref{eq:crkfr.avg} and $\bS_h^\delta$ is a local approximation to the time averaged source term $\int_{t^n}^{t^{n+1}} \bss \ud t$. We describe the local compact Runge-Kutta procedure~(Section~\ref{sec:crkfr}) for solving~\eqref{eq:con.law.source.term} with the cRK33 scheme~\eqref{eq:crk33}. The other schemes will be similar. The procedure to compute the local time average flux and source terms is given by
\begin{equation} \label{eq:crk33.source}
\uU = \frac{1}{4}  \uu_h^n + \frac{3}{4}  \uu_h^{(3)}, \qquad
\F = \frac{1}{4}  \pf ( \uu_h^n ) + \frac{3}{4}  \pf( \uu_h^{(3)}), \qquad \bS = \frac{1}{4}  \bss ( \uu_h^n ) + \frac{3}{4}  \bss( \uu_h^{(3)}),
\end{equation}
where
\begin{align*}
\uu^{(2)}_h & = \uu_h^n - \frac{\mathLaplace t}{3} \dlocx
\pf ( \uu_h^n ) + \frac{\mathLaplace t}{3}
\bss ( \uu_h^n ), \\
\uu_h^{(3)} & =  \uu_h^n
- \frac{2}{3} \mathLaplace t \dlocx  \pf (\uu_h^{(2)})
+ \frac{2}{3} \mathLaplace t \bss (\uu_h^{(2)}).
\end{align*}
As justified in Appendix~\ref{app:time.averaged.flux}, $\bS_h^\delta$ is an order $N+1$ approximation of the time averaged source term. Unlike the time average flux, the approximation $\bS_h^\delta$ used in~\eqref{eq:crkfr.source} is completely specified by~\eqref{eq:crk33.source}. In particular, it does not require any interelement correction.
\bibliographystyle{siam}
\bibliography{references}
\end{document}